\documentclass[12pt,a4paper]{amsart}
\usepackage{amsfonts}
\usepackage{amsmath}
\usepackage{amssymb}

\def\gF{\mathfrak{F}}

\def\gH{\mathfrak{H}}

\def\gM{\mathfrak{M}}
\def\gN{\mathfrak{N}}

\begin{document}

\title[INCLUSIONS AND QUANTUM GROUPOIDS III] {Inclusions of Von Neumann Algebras and Quantum Groupo\"ids III}
\author{Michel Enock}
\address{Institut de Math\'ematiques de Jussieu, Unit\'{e} Mixte Paris 6 / Paris 7 /
CNRS de Recherche 7586 \\175, rue du Chevaleret, Plateau 7E, F-75013 Paris}
 \email{enock@math.jussieu.fr}
\date{may 04}
\maketitle
\begin{abstract}
In a former article, in collaboration with Jean-Michel Vallin, we have constructed two "quantum
groupo\"{\i}ds" dual to each other, from a depth 2 inclusion of von Neumann algebras
$M_0\subset M_1$, in such a way that the canonical Jones'tower associated to the inclusion can be
described as a tower of successive crossed-products by these two structures. We are now
investigating in greater details these structures in the presence of an appropriate modular theory on the basis $M'_0\cap M_1$, and we show how these examples fit with Lesieur's "measured quantum
groupo\"{\i}ds".  
\end{abstract}

\maketitle
\newpage
\section{Introduction}
\label{intro}
\subsection{}
 In two articles ([V1], [V2]), J.-M. Vallin has introduced two notions (pseudo-multiplicative
unitary, Hopf-bimodule), in order to generalize, up to the groupo\"{\i}d
case, the classical notions of multiplicative unitary [BS] and of Hopf-von Neumann algebras [ES]
which were introduced to describe and explain duality of groups, and leaded to appropriate notions
of quantum groups ([ES], [W1], [W2], [BS], [MN], [W3], [KV]). 
\\ In a former article [EV], J.-M. Vallin and the author have constructed, from a depth 2 inclusion of
von Neumann algebras $M_0\subset M_1$, with an operator-valued weight $T_1$ verifying a regularity
condition, a pseudo-multiplicative unitary, which leaded to two structures of Hopf bimodules, dual
to each other. Moreover, we have then
construct an action of one of these structures on the algebra $M_1$ such that $M_0$
is the fixed point subalgebra, the algebra $M_2$ given by the basic construction being then
isomorphic to the crossed-product. We construct on $M_2$ an action of the other structure, which
can be considered as the dual action.
\\  If the inclusion
$M_0\subset M_1$ is irreducible, we recovered quantum groups, as proved and studied in former papers
([EN], [E1]).
\\ Therefore, this construction leads to a notion of "quantum groupo\"{\i}d", and a construction of a
duality within "quantum groupo\"{\i}ds". 
\subsection{}
In a finite-dimensional setting, this construction can be
mostly simplified, and is studied in [NV1], [V3], [BSz1],
[BSz2], [Sz], and examples are described. In [NV2], the link between these "finite quantum
groupo\"{\i}ds" and depth 2 inclusions of
$II_1$ factors is given. 
\subsection{}
In a second article [E2], we went on studying these structures, in order to construct, from the pseudo-multiplicative unitary, the analog of an antipod; we showed that this antipod bears a polar decomposition which leads to a co-inverse and a deformation one-parameter group of the structure. We also got a left invariant operator-valued weight (and, using the co-inverse, a right-invariant one). All these results were claimed when adding a modular hypothesis on the basis $M'_0\cap M_1$, namely that there exists on this basis a normal semi-finite faithful weight $\chi$ whose modular automorphism group $\sigma_t^{\chi}$ is equal to the modular automorphism group of the  operator-valued weight $T_1$. 
\\ It appears that this property is much too restrictive; in particular, it does not go throught climbing up the tower and the inclusion $M_1\subset M_2$ does not bear the same property, unless the basis $M'_0\cap M_1$ is semi-finite. There is a mistake in ([E2], 4.1), and all the results claimed in that article are therefore proved only if the basis $M'_0\cap M_1$ is semi-finite. I am indebted to Franck Lesieur who pointed out this mistake. 
\newline
We consider now a wider hypothesis, namely that there exists on the basis $M'_0\cap M_1$ a normal semi-finite faithful weight which is invariant under the modular automorphism group of the operator-valued weight $T_1$. 
\subsection{}
Franck Lesieur introduced in his thesis [L] a notion of "measured quantum groupo\"{\i}ds", in which a modular hypothesis is required. We discuss when our construction fits with his axioms.   
 \subsection{}
The paper is organized as follows : in chapter \ref{pr}, we recall all the preliminaries needed for that theory, mostly Connes-Sauvageot relative tensor product and Hopf-bimodules; in chapter \ref{pmu} is recalled the notion of a pseudo-multiplicative unitary, and how, in [EV], was associated such an object to a depth 2 inclusion of von Neumann algebras, with the appropriate technical conditions. In chapter \ref{modular} is developped a modular theory on the basis, which generalizes what was done in [E2]; in chapter \ref{bW}, we come back to the multiplicative unitary and obtain analytical properties which allow us to construct an antipod, its polar decomposition, a right-invariant operator-valued weight and a left-invariant operator-valued weight; in chapter \ref{modular2}, mimicking [L], we obtain a modular theory for the left-invariant weight and the right-invariant weight, and, in chapter \ref{density}, we obtain a density result which will be usefull in the sequel of the theory. Then, in chapter \ref{relcom}, we easily put all these structures at the level of relative commutants constructed from the Jones' tower, and we finish in chapter \ref{adapted} by discussing when the objects found fit with Lesieur's axioms of measured quantum groupoids. 

\subsection{}
The author is mostly indebted to F. Lesieur, S. Vaes, J.-M. Vallin
and L. Va\u{\i}nerman for many fruitful conversations.

\section{Preliminaries}
\label{pr}
In this chapter are mostly recalled definitions and notations about Connes' spatial
theory (\ref{spatial}, \ref{rel}) and the fiber product construction (\ref{fiber}, \ref{slice})
which are the main technical tools of that theory. The definition of Hopf-bimodules is given
(\ref{Hbimod}). 
 
\subsection{Spatial theory [C1], [S2], [T]}
 \label{spatial}
 Let $N$ be a von Neumann algebra, 
and let $\psi$ be a faithful semi-finite normal weight on $N$; let $\gN_{\psi}$, 
$\gM_{\psi}$, $H_{\psi}$, $\pi_{\psi}$, $\Lambda_{\psi}$,$J_{\psi}$, 
$\Delta_{\psi}$,... be the canonical objects of the Tomita-Takesaki construction 
associated to the weight $\psi$. Let $\alpha$  be a non-degenerate normal representation of $N$ on a
Hilbert space
$\mathcal{H}$. We may as well consider $\mathcal{H}$ as a left $N$-module, and write it then
$_\alpha\mathcal{H}$. Following ([C1], definition 1), we define the set of 
$\psi$-bounded elements of $_\alpha\mathcal{H}$ as :
\[D(_\alpha\mathcal{H}, \psi)= \lbrace \xi \in \mathcal{H};\exists C < \infty ,\| \alpha (y) \xi\|
\leq C \| \Lambda_{\psi}(y)\|,\forall y\in \gN_{\psi}\rbrace\]
Then, for any $\xi$ in $D(_\alpha\mathcal{H}, \psi)$, there exists a bounded operator
$R^{\alpha,\psi}(\xi)$ from $H_\psi$ to $\mathcal{H}$,  defined, for all $y$ in $\gN_\psi$ by :
\[R^{\alpha,\psi}(\xi)\Lambda_\psi (y) = \alpha (y)\xi\]
If there is no ambiguity about the representation $\alpha$, we shall write
$R^{\psi}(\xi)$ instead of $R^{\alpha,\psi}(\xi)$. This operator belongs to $Hom_N (H_\psi , \mathcal{H})$;
therefore, for any
$\xi$, $\eta$ in
$D(_\alpha\mathcal{H}, \psi)$, the operator :
\[\theta^{\alpha,\psi} (\xi ,\eta ) =  R^{\alpha,\psi}(\xi)R^{\alpha,\psi}(\eta)^*\]
belongs to $\alpha (N)'$; moreover, $D(_\alpha\mathcal{H}, \psi)$ is dense ([C1], lemma 2), stable
under $\alpha (N)'$, and the linear span generated by the operators $\theta^{\alpha,\psi} (\xi ,\eta
)$ is a dense ideal in $\alpha (N)'$. With the same hypothesis, the operator :
\[<\xi,\eta>_{\alpha,\psi} = R^{\alpha,\psi}(\eta)^* R^{\alpha,\psi}(\xi)\]
belongs to $\pi_{\psi}(N)'$. Using Tomita-Takesaki's theory, this last algebra is equal to
$J_\psi
\pi_{\psi}(N)J_\psi$, and therefore anti-isomorphic to $N$ (or isomorphic to the opposite von
Neumann algebra $N^o$). We shall consider now $<\xi,\eta>_{\alpha,\psi}$ as an element of $N^o$, and
the linear span generated by these operators is a dense ideal in $N^o$.
\newline
There exists ([C], prop.3) a family  $(e_i)_{i\in I}$ of 
$\psi$-bounded elements of $_\alpha\mathcal H$, such that
\[\sum_i\theta^{\alpha, \psi} (e_i ,e_i )=1\]
Such a family will be called a $(N,\psi )$-basis of $_\alpha\mathcal H$. It is possible ([EN] 2.2) to construct 
a $(N,\psi )$-basis of $_\alpha\mathcal H$, $(e_i)_{i\in I}$, such that the
operators $R^{\alpha, \psi}(e_i)$ are partial isometries with final supports 
$\theta^{\alpha, \psi}(e_i ,e_i )$ 2 by 2 orthogonal, and such that, if $i\neq j$, then 
$<e_i ,e_j>_{\alpha, \psi}=0$. Such a family will be called a $\alpha$-orthogonal basis of $\mathcal H$. 
\newline
Let $\beta$ be a normal non-degenerate
anti-representation
 of
$N$ on
$\mathcal{H}$. We may then as well consider $\mathcal{H}$ as a right $N$-module, and write it $\mathcal{H}_\beta$, or
consider
$\beta$ as a normal non-degenerate representation of the opposite von Neumann algebra $N^o$, and
consider
$\mathcal{H}$ as a left $N^o$-module. 
\newline
We can then define on $N^o$ the opposite faithful
 semi-finite normal weight $\psi ^o$; we have $\gN_{\psi ^o}=\gN_\psi^*$, and 
 the Hilbert space $H_{\psi ^o}$ will be, as usual, identified with $H_\psi$, 
by the identification, for all $x$ in $\gN_\psi$, of 
$\Lambda_{\psi ^o}(x^*)$ with $J_\psi \Lambda_\psi (x)$.
\newline
From these remarks, we infer that the set of 
$\psi ^o$-bounded elements of
$\mathcal{H}_\beta$ is :
\[D(\mathcal{H}_\beta, \psi^o) = \lbrace\xi\in\mathcal{H} ;\exists C < \infty ,
\|\beta (y^*)\xi\|
\le C \| \Lambda_{\psi}(y)\|,\forall y\in \gN_{\psi}\rbrace\]
and, for any $\xi$ in $D(\mathcal{H}_\beta, \psi^o)$ and $y$ in $\gN_\psi$,
 the bounded operator $R^{\beta,\psi^{o}}(\xi )$ is given by the formula :
\[R^{\beta,\psi{^o}}(\xi)J_{\psi}\Lambda_\psi (y) = \beta (y^*)\xi\]
This operator belongs to $Hom_{N^{o}}(H_\psi ,\mathcal{H} )$.
Moreover, $D(\mathcal{H}_\beta, \psi^o)$ is dense, stable under $\beta( N)'=P$, 
and, for all $y$ in $P$, we have :
\[R^{\beta,\psi{^o}}(y\xi)= yR^{\beta,\psi{^o}}(\xi)\]
Then, for any $\xi$, $\eta$ in $D(\mathcal{H}_\beta, \psi^o)$, the operator 
\[\theta^{\beta, \psi^{o}}(\xi ,\eta )=R^{\beta,\psi{^o}}(\xi)R^{\beta,\psi{^o}}(\eta)^*\] 
belongs to $P$, and the linear span generated by these operators is a dense ideal in $P$;
moreover, the operator-valued product 
$<\xi ,\eta >_{\beta,\psi^o}= R^{\beta,\psi{^o}}(\eta)^*R^{\beta,\psi{^o}}(\xi)$
belongs to $\pi_\psi (N)$; we shall consider now, for simplification, that $<\xi ,\eta
>_{\beta,\psi^o}$ belongs to $N$, and the linear span generated by these operators is a dense ideal
in $N$. More precisely, $<\xi ,\eta
>_{\beta,\psi^o}$ belongs to $\gM_\psi$ ([C], lemma 4) and we have ([S1], lemme 1.5) 
\[\Lambda_\psi(<\xi, \eta>_{\beta,\psi^o})=R^{\beta,\psi{^o}}(\eta)^*\xi\]
A $(N^{o},\psi^o)$-basis of $\mathcal H_\beta$ is a family 
$(e_i )_{i\in I}$ of $\psi^o$-bounded elements of $\mathcal H_\beta$, such that 
\[\sum_i\theta^{\beta, \psi^{o}}(e_i ,e_i )=1\]
We have then, for all $\xi$ in $D(\mathcal H_\beta)$ :
\[\xi=\sum_i R^{\beta,\psi^o}(e_i)\Lambda_\psi(<\xi, e_i>_{\beta,\psi^o})\]
It is possible to choose the $(e_i )_{i\in I}$ such that
the $R^{\beta, \psi{^o}}(e_i)$ are partial isometries, with final supports
$\theta^{\beta, \psi^{o}}(e_i ,e_i )$ 2 by 2 orthogonal, and $<e_i, e_j>_{\beta, \psi^o}=0$ if $i\neq j$; such
a family will be then called a $\beta$-orthogonal basis of $\mathcal H$. We have then 
\[R^{\beta, \psi{^o}}(e_i)=\theta^{\beta, \psi^{o}}(e_i ,e_i )R^{\beta, \psi{^o}}(e_i)=
R^{\beta, \psi{^o}}(e_i)<e_i, e_i>_{\beta, \psi^o}\]
 
\subsection{Jones' basic construction}
\label{basic}
Let $M_0\subset M_1$ be an inclusion of von Neumann algebras  (for simplification, these algebras will be supposed to be $\sigma$-finite), equipped with a normal faithful semi-finite operator-valued weight $T_1$ from $M_1$ to $M_0$ (to be more precise, from $M_1^{+}$ to the extended positive elements of $M_0$ (cf. [T] IX.4.12)). Let $\psi_0$ be a normal faithful semi-finite weight on $M_0$, and $\psi_1=\psi_0\circ T_1$; for $i=0,1$, let $H_i=H_{\psi_i}$, $J_i=J_{\psi_i}$, $\Delta_i=\Delta_{\psi_i}$ be the usual objects constructed by the Tomita-Takesaki theory associated to these weights. Following ([J], 3.1.5(i)), the von Neumann algebra $M_2=J_1M'_0J_1$ defined on the Hilbert space $H_1$ will be called the basic construction made from the inclusion $M_0\subset M_1$. We have $M_1\subset M_2$, and we shall say that the inclusion $M_0\subset M_1\subset M_2$ is standard.  Using then Haagerup's construction ([T], IX.4.24), it is possible to construct a normal semi-finite faithful operator-valued weight $T_2$ from $M_2$ to $M_1$ ([EN], 10.7), which will be called the basic construction made from $T_1$. Repeating this construction, we obtain by recurrence successive basic constructions, which lead to Jones' tower $(M_i)_{i\in \Bbb{N}}$ of von Neumann algebras, which is the inclusion 
\[M_0\subset M_1\subset M_2\subset M_3\subset M_4\subset ...\]
which is equipped (for $i\geq 1$)  with normal faithful semi-finite operator-valued weights $T_i$ from $M_i$ to $M_{i-1}$. We define then, by recurrence, the weight $\psi_i=\psi_{i-1}\circ T_i$ on $M_i$, and we shall write $H_i$, $J_i$, $\Delta_i$ instead of $H_{\psi_i}$, etc. We shall define the mirroring $j_i$ on $\mathcal L(H_i)$ by $j_i(x)=J_ix^*J_i$, for all $x$ in $\mathcal L(H_i)$.
\newline
Following ([EN] 10.6), for $x$ in $\gN_{T_i}$, we shall define $\Lambda_{T_i}(x)$ by the following formula, for all $z$ in $\gN_{\psi_{i-1}}$ :
\[\Lambda_{T_i}(x)\Lambda_{\psi_{i-1}}(z)=\Lambda_{\psi_i}(xz)\]
Then, $\Lambda_{T_i}(x)$ belongs to $Hom_{M_{i-1}^o}(H_{i-1}, H_i)$; if $x$, $y$ belong to $\gN_{T_i}$, then $\Lambda_{T_i}(x)^*\Lambda_{T_i}(y)=T_i(x^*y)$, and $\Lambda_{T_i}(x)\Lambda_{T_i}(y)^*$ belongs to $M_{i+1}$; more precisely, it belongs to $\gM_{T_{i+1}}$, and $T_{i+1}(\Lambda_{T_i}(x)\Lambda_{T_i}(y)^*)=xy^*$. 
\newline
By Tomita-Takesaki theory, the Hilbert space $H_1$ bears a natural structure of $M_1-M_1^o$-bimodule, and, therefore, by restriction, of $M_0-M_0^o$-bimodule. Let us write $r$ for the canonical representation of $M_0$ on $H_1$, and $s$ for the canonical antirepresentation given, for all $x$ in $M_0$, by $s(x)=J_1r(x)^*J_1$. Let us have now a closer look to the subspaces $D(H_{1s}, \psi_0^o)$ and $D(_rH_1, \psi_0)$.

\subsection{Proposition}
\label{prop}
 {\it (i) Let $x$ be in $\gN_{T_1}\cap\gN_{\psi_1}$; then $\Lambda_{\psi_1}(x)$ belongs to $D(H_{1s}, \psi_0^o)$, and $R^{s, \psi_0^o}(\Lambda_{\psi_1}(x))=\Lambda_{T_1}(x)$. Moreover, $J_{\psi_1}\Lambda_{\psi_1}(x)$ belongs to $D(_rH_1, \psi_0)$, and $R^{r, \psi_0}(J_{\psi_1}\Lambda_{\psi_1}(x))=J_{\psi_1}\Lambda_{T_1}(x)J_{\psi_0}$. Conversely, let $\xi$ in $D(H_{1s}, \psi_0^o)$; then, there exists a sequence $x_n$ in $\gN_{T_1}\cap\gN_{\psi_1}$ such that $\Lambda_{\psi_1}(x_n)$ is converging to $\xi$, and $\Lambda_{T_1}(x_n)$ is weakly converging to $R^{s, \psi_0^o}(\xi)$. }
 \newline
 {\it (ii) Let $a$ in $\gN_{\psi_1}\cap \gN_{T_1}\cap\gN_{\psi_1}^*\cap \gN_{T_1}^*$, analytic with respect to $\psi_1$, such that, for all $z$ in $\mathbb{C}$, $\sigma_z(a)$ belongs to $\gN_{\psi_1}\cap \gN_{T_1}\cap\gN_{\psi_1}^*\cap \gN_{T_1}^*$ (we shall denote $T_{\psi_1, T_1}$ the set of such elements); then $\Lambda_{\psi_1}(a)$ belongs to $D(H_{1s}, \psi_0^o)\cap D(_rH_1, \psi_0)$}
\newline
{\it (iii) The subspace $D(H_{1s}, \psi_0^o)\cap D(_rH_1, \psi_0)$ is dense in $H_1$. Moreover, if $\xi$ belongs to $D(H_{1s}, \psi_0^o))$, there exists a sequence $\xi_n$ in $D(H_{1s}, \psi_0^o)\cap D(_rH_1, \psi_0)$ such that $R^{s, \psi_0^o}(\xi_n)$ is weakly converging to $R^{s, \psi_0^o}(\xi)$; if $\xi'$ belongs to $D(_rH_1, \psi_0)$, there exists a sequence $\xi'_n$ in $D(H_{1s}, \psi_0^o)\cap D(_rH_1, \psi_0)$ such that $R^{r, \psi_0}(\xi'_n)$ is weakly converging to $R^{r, \psi_0}(\xi')$.}

 \begin{proof}
The first results of (i) are just standard calculation. Let us prove the converse part. Let us consider the basic construction $M_1\subset M_2$, and  the normal semi-finite faithful operator valued weight $T_2$ from $M_2$ to $M_1$. Then, the operator $T_2(\theta^{s, \psi_0^o}(\xi, \xi))$ is a positive self-adjoint operator affiliated to $M_1$, and let us write :
 \[T_2(\theta^{s, \psi_0^o}(\xi, \xi))=\int_0^{\infty}\lambda de_{\lambda}\]
 Let us put $p_n=\int_0^n de_{\lambda}$; then $p_n$ belongs to $M_1$, and $p_n\xi$ belongs to $D(H_{1s}, \psi_0^o)$. Moreover, we get that :
 \[T_2(\theta^{s, \psi_0^o}(p_n\xi, p_n\xi))=\int_0^n\lambda de_{\lambda}\in M_1\]
 from which we get that $\theta^{s, \psi_0^o}(p_n\xi, p_n\xi)$ belongs to $\gM_{T_2}^+$; therefore, we get that there exist $x_n$ in $\gN_{\psi_1}$ such that $p_n\xi=\Lambda_{\psi_1}(x_n)$, and  
 $\theta^{s, \psi^o}(p_n\xi, p_n\xi)=x_nx_n^*$. 
 \\ Let us take now $x$ in $M_0$, in the Tomita algebra associated to the weight $\psi_o$. As we have :
 \[R^{s, \psi_0^o}(p_n\xi)J_{\psi_0}\Lambda_{\psi_0}(x)=s(x^*)p_n\xi=s(x^*)\Lambda_{\psi_1}(x_n)=
 \Lambda_{\psi_1}(x_n\sigma_{-i/2}^{\psi_0}(x^*))\]
 we get :
 \begin{eqnarray*}
 \psi_1(\sigma_{-i/2}^{\psi_0}(x^*)^*x_n^*x_n\sigma_{-i/2}^{\psi_0}(x^*))
 &\leq&\|R^{s, \psi_0^o}(p_n\xi)\|^2\|\Lambda_{\psi_0}(x)\|^2\\
 &\leq&\|R^{s, \psi_0^o}(\xi)\|^2\|\Lambda_{\psi_0}(x)\|^2
 \end{eqnarray*}
 from which we infer that the element $T_1(x_n^*x_n)$ of the extended positive part of $M_0$ satisfies :
 \[(T_1(x_n^*x_n)J_{\psi_0}\Lambda_{\psi_0}(x)|J_{\psi_0}\Lambda_{\psi_0}(x))\leq\|R^{s, \psi_0^o}(\xi)\|^2\|\Lambda_{\psi_0}(x)\|^2\]
 Therefore, $T_1(x_n^*x_n)$ is bounded, i.e. $x_n$ belongs to $\gN_{T_1}$, and, more precisely, for all $n$ in $\mathbb{N}$,  we have $\|T_1(x_n^*x_n)\|\leq\|R^{s, \psi_0^o}(\xi)\|^2$. 
 \\For any $x$ in $\gN_{\psi_0}$, we have :
 \begin{multline*}
 \Lambda_{T_1}(x_n)J_{\psi_0}\Lambda_{\psi_0}(x)=R^{s, \psi_0^o}(\Lambda_{\psi_1}(x_n))J_{\psi_0}\Lambda_{\psi_0}(x)\\
 =s(x^*)\Lambda_{\psi_1}(x_n)=s(x^*)p_n\xi
 \end{multline*}
 which is converging to $s(x^*)\xi=R^{s, \psi_0^o}(\xi)J_{\psi_0}\Lambda_{\psi_0}(x)$. As all the norms $\|\Lambda_{T_1}(x_n)\|$ are uniformly bounded, we get (i). 
 \newline
As $D(_rH_1, \psi_0)=J_{\psi_1}D(H_{1s}, \psi_0^o)$ and $J_{\psi_1}\Lambda_{\psi_1}(a)=\Lambda_{\psi_1}(\sigma_{-i/2}^{\psi_1}(a^*))$, the result (ii) is clear. 
\newline
From which we get the density of $D(H_{1s}, \psi_0^o)\cap D(_rH_1, \psi_0)$, by ([EN], 10.12). More precisely, if $x$ belongs to $\gN_{T_1}\cap\gN_{\psi_1}$, there exists, by ([EN], 10.12) a sequence $x_n$ in $\gN_{T_1}\cap\gN_{\psi_1}$, with $\|x_n\|\leq\|x\|$, $\|T_1(x_n^*x_n)\|\leq\|T_1(x^*x)\|$, such that $\Lambda_{\psi_1}(x_n)$ is converging to $\Lambda_{\psi_1}(x)$ and, using (i), such that $\Lambda_{\psi_1}(x_n)$ belongs to $D(H_{1s}, \psi_0^o)\cap D(_rH_1, \psi)$. Using same arguments as in (i), we find also that $\Lambda_{T_1}(x_n)$ is weakly converging to $\Lambda_{T_1}(x)$. With the help again of (i), starting from $\xi$ in $D(H_{1s}, \psi_0^o)$, we get a sequence $\xi_n$ in $D(H_{1s}, \psi_0^o)\cap D(_rH_1, \psi_0)$. 
\\Again, as $D(_rH_1, \psi_0)=J_{\psi_1}D(H_{1s}, \psi_0^o)$, we finish the proof. \end{proof}

\subsection{Proposition}
\label{propT}
{\it Let $a$, $b$ in $\gN_{\psi_1}$; then $T_1(a^*a)$ and $T_1(b^*b)$ are positive self-adjoint closed operators which verify :}
\[<T_1(b^*b), \omega_{J_{\psi_1}\Lambda_{\psi_1}(a)}>=<T_1(a^*a), \omega_{J_{\psi_1}\Lambda_{\psi_1}(b)}>\]

\begin{proof}
Let us suppose first that $b$ belongs to $\gN_{\psi_1}\cap\gN_{T_1}$. Let $z$ in the Tomita algebra of $\psi_0$; we have :
\begin{eqnarray*}
(zJ_{\psi_1}\Lambda_{\psi_1}(b)|J_{\psi_1}\Lambda_{\psi_1}(b))
&=&\psi_1(b^*b\sigma_{-i/2}^{\psi_1}(z^*))^-\\
&=&\psi_0(T_1(b^*b)\sigma_{-i/2}^{\psi_0}(z^*))^-\\
&=&(J_{\psi_0}\Lambda_{\psi_0}(z)|\Lambda_{\psi_0}(T_1(b^*b)))^-\\
&=&(\Lambda_{\psi_0}(z)|J_{\psi_0}\Lambda_{\psi_0}(T_1(b^*b)))
\end{eqnarray*}
which, by density, remains true for any $z$ in $\gN_{\psi_0}$. Therefore, if $a$ belongs to $\gN_{\psi_1}\cap\gN_{T_1}$, we get :
\[(T_1(a^*a)J_{\psi_1}\Lambda_{\psi_1}(b)|J_{\psi_1}\Lambda_{\psi_1}(b))=(T_1(b^*b)J_{\psi_1}\Lambda_{\psi_1}(a)|J_{\psi_1}\Lambda_{\psi_1}(a))\]
Let us suppose now that $a$ belongs only to $\gN_{\psi_1}$. It is then well known (see for instance [EN], 10.6) that $T_1(a^*a)$ is a positive self-adjoint closed operator which can be written $\int_0^\infty \lambda dp_\lambda$, and that, for all $n$, $ap_n$ belongs to $\gN_{T_1}\cap\gN_{\psi_1}$, and that the sequence $\Lambda_{\psi_1}(ap_n)$ is converging to $\Lambda_{\psi_1}(a)$. In that situation, we get that :
\begin{eqnarray*}
<T_1(a^*a), \omega_{J_{\psi_1}\Lambda_{\psi_1}(b)}>
&=&lim_n<T_1(p_na^*ap_n), \omega_{J_{\psi_1}\Lambda_{\psi_1}(b)}>\\
&=&lim_n(T_1(b^*b)J_{\psi_1}\Lambda_{\psi_1}(ap_n)|J_{\psi_1}\Lambda_{\psi_1}(ap_n))\\
&=&(T_1(b^*b)J_{\psi_1}\Lambda_{\psi_1}(a)|J_{\psi_1}\Lambda_{\psi_1}(a))
\end{eqnarray*}
If now $b$ belongs only to $\gN_{\psi_1}$, using the same trick and the fact we are dealing with closed operators, we obtain the result. \end {proof}


\subsection{Relative tensor product [C1], [S2], [T]}
\label{rel}
Using the notations of \ref{spatial}, let now $\mathcal{K}$ be another Hilbert space on which there exists
a non-degenerate representation
$\gamma$ of
$N$. Following J.-L. Sauvageot ([S2], 2.1), we define
the relative tensor product $\mathcal{H}\underset{\psi}{_\beta\otimes_\gamma}\mathcal{K}$ as the
Hilbert space obtained from the algebraic tensor product $D(\mathcal{H}_\beta ,\psi^o )\odot
\mathcal{K} $ equipped with the scalar product defined, for $\xi_1$, $\xi_2$ in $D(\mathcal{H}_\beta
,\psi^o )$,
$\eta_1$, $\eta_2$ in $\mathcal{K}$, by 
\[(\xi_1\odot\eta_1 |\xi_2\odot\eta_2 )=(\gamma(<\xi_1 ,\xi_2 >_{\beta,\psi^o})\eta_1 |\eta_2 )\]
where we have identified $N$ with $\pi_\psi (N)$ to simplifly the notations.
\newline
The image of $\xi\odot\eta$ in $\mathcal{H}\underset{\psi}{_\beta\otimes_\gamma}\mathcal{K}$ will be
denoted by
$\xi\underset{\psi}{_\beta\otimes_\gamma}\eta$. We shall use intensively this construction; one
should bear in mind that, if we start from another faithful semi-finite normal weight $\psi '$, we
get another Hilbert space $\mathcal{H}\underset{\psi'}{_\beta\otimes_\gamma}\mathcal{K}$; there exists an isomorphism $U^{\psi, \psi'}_{\beta, \gamma}$ from $\mathcal{H}\underset{\psi}{_\beta\otimes_\gamma}\mathcal{K}$ to $\mathcal{H}\underset{\psi'}{_\beta\otimes_\gamma}\mathcal{K}$, which is unique up to some functorial property ([S2], 2.6) (but this isomorphism
does not send  $\xi\underset{\psi}{_\beta\otimes_\gamma}\eta$ on
$\xi\underset{\psi'}{_\beta\otimes_\gamma}\eta$ !). 
\newline
When no confusion is possible about the representation and the anti-representation, we shall write
$\mathcal{H}\otimes_{\psi}\mathcal{K}$ instead of
$\mathcal{H}\underset{\psi}{_\beta\otimes_\gamma}\mathcal{K}$, and $\xi\otimes_\psi\eta$ instead
of
$\xi\underset{\psi}{_\beta\otimes_\gamma}\eta$.
\newline
For any $\xi$ in $D(\mathcal{H}_\beta,
\psi^o)$, we define the bounded linear application $\lambda_\xi^{\beta, \gamma}$ from $\mathcal
K$ to
$\mathcal{H}\underset{\psi}{_\beta\otimes_\gamma}\mathcal{K}$ by, for all $\eta$ in $\mathcal K$,
$\lambda_\xi^{\beta, \gamma} (\eta)=\xi\underset{\psi}{_\beta\otimes_\gamma}\eta$. We shall write
$\lambda_\xi$ if no confusion is possible. We get ([EN], 3.10) :
\[\lambda_\xi^{\beta, \gamma}=R^{\beta, \psi^o}(\xi)\otimes_\psi 1_\mathcal K\]
where we recall the canonical identification (as left $N$-modules) of 
$L^2(N)\otimes_\psi\mathcal K$ with $\mathcal K$. We have :
\[(\lambda_\xi^{\beta, \gamma})^*\lambda_\xi^{\beta, \gamma}=\gamma(<\xi, \xi>_{\beta, \psi^o})\]
\newline
In ([S1] 2.1), the relative tensor product
$\mathcal{H}\underset{\psi}{_\beta\otimes_\gamma}\mathcal{K}$ is defined also, if
$\xi_1$, $\xi_2$ are in $\mathcal{H}$, $\eta_1$, $\eta_2$ are in $D(_\gamma\mathcal{K},\psi)$, by the
following formula :
\[(\xi_1\odot\eta_1 |\xi_2\odot\eta_2 )= (\beta(<\eta_1, \eta_2>_{\gamma,\psi})\xi_1 |\xi_2)\]
which leads to the the definition of a relative flip $\sigma_\psi$ which will be an isomorphism from
$\mathcal{H}\underset{\psi}{_\beta\otimes_\gamma}\mathcal{K}$ onto
$\mathcal{K}\underset{\psi^o}{_\gamma\otimes _\beta}\mathcal{H}$, defined, for any 
$\xi$ in $D(\mathcal{H}_\beta ,\psi^o )$, $\eta$ in $D(_\gamma \mathcal{K},\psi)$, by :
\[\sigma_\psi (\xi\otimes_\psi\eta)=\eta\otimes_{\psi^o}\xi\]
This allows us to define a relative flip $\varsigma_\psi$ from
$\mathcal{L}(\mathcal{H}\underset{\psi}{_\beta\otimes_\gamma}\mathcal{K})$ to $\mathcal{L}
(\mathcal{K}\underset{\psi^o}{_\gamma\otimes_\beta}\mathcal{H})$ which sends
$X$ in 
$\mathcal{L}(\mathcal{H}\underset{\psi}{_\beta\otimes_\gamma}\mathcal{K})$ onto
$\varsigma_\psi(X)=\sigma_\psi X\sigma_\psi^*$. Starting from another faithful semi-finite normal
weight $\psi'$, we get a von Neumann algebra
$\mathcal{L}(\mathcal{H}\underset{\psi'}{_\beta\otimes_\gamma}\mathcal{K})$ which is isomorphic to
$\mathcal{L}(\mathcal{H}\underset{\psi}{_\beta\otimes_\gamma}\mathcal{K})$, and a von Neumann
algebra $\mathcal{L} (\mathcal{K}\underset{\psi^{'o}}{_\gamma\otimes_\beta}\mathcal{H})$ which is
isomorphic to
$\mathcal{L} (\mathcal{K}\underset{\psi^o}{_\gamma\otimes_\beta}\mathcal{H})$; as we get that :
 \[\sigma_{\psi'}\circ U^{\psi, \psi'}_{\beta, \gamma}=U^{\psi^o, \psi'^o}_{\gamma, \beta}\]
 we see that these isomorphisms exchange $\varsigma_\psi$ and
$\varsigma_{\psi '}$. Therefore, the homomorphism $\varsigma_{\psi}$ can be denoted $\varsigma_N$
without any reference to a specific weight.
\newline
We may define, for any $\eta$ in $D(_\gamma\mathcal{K}, \psi)$, an application
$\rho_\eta^{\beta, \gamma}$ from $\mathcal H$ to
$\mathcal{H}\underset{\psi}{_\beta\otimes_\gamma}\mathcal{K}$ by
$\rho_\eta^{\beta, \gamma} (\xi)=\xi\underset{\psi}{_\beta\otimes_\gamma}\eta$. We shall write
$\rho_\eta$ if no confusion is possible. We get that :
\[(\rho_\eta^{\beta, \gamma})^*\rho_\eta^{\beta, \gamma}=\beta(<\eta, \eta>_{\gamma, \psi})\]
\newline
We recall, following
([S2], 2.2b) that, for all $\xi$ in $\mathcal{H}$, $\eta$ in $D(_\gamma\mathcal{K},\psi)$, $y$ in
$N$, analytic with respect to $\psi$, we have :
\[\beta (y)\xi \otimes_\psi\eta=\xi\otimes_\psi\gamma(\sigma^\psi_{-i/2}(y))\eta\]
Let $x$ be an element
of
$\mathcal{L}(\mathcal{H})$, commuting with the right action of $N$ on $\mathcal{H}_\beta$ (i.e. $x\in\beta(N)'$). It
is possible to define an operator $x\underset{\psi}{_\beta\otimes_\gamma} 1_{\mathcal{K}}$ on
$\mathcal{H}\underset{\psi}{_\beta\otimes_\gamma}
\mathcal{K}$. By the same way, if $y$ commutes with the left action of $N$ on
$_\gamma\mathcal{K}$ (i.e. $y\in\gamma(N)'$), it is possible to define
$1_{\mathcal{H}}\underset{\psi}{_\beta\otimes_\gamma}y$ on
$\mathcal{H}\underset{\psi}{_\beta\otimes_\gamma} \mathcal{K}$, and by composition, it is possible
to define then
$x\underset{\psi}{_\beta\otimes_\gamma} y$. If we start from another faithful semi-finite normal
weight $\psi '$, the canonical isomorphism $U^{\psi, \psi'}_{\beta, \gamma}$ from $\mathcal{H}\underset{\psi}{_\beta\otimes_\gamma}
\mathcal{K}$ to $\mathcal{H}\underset{\psi'}{_\beta\otimes_\gamma} \mathcal{K}$ sends
$x\underset{\psi}{_\beta\otimes_\gamma} y$ on $x\underset{\psi'}{_\beta\otimes_\gamma} y$ ([S2],
2.3 and 2.6); therefore, this operator can be denoted $x\underset{N}{_\beta\otimes_\gamma} y$
without any reference to a specific weight\vspace{2mm}.  
\newline
Let us suppose now that $\mathcal{K}$ is a $N-P$ bimodule; that means that there exists a von
Neumann algebra $P$, and a non-degenerate normal anti-representation $\epsilon$ of $P$ on
$\mathcal{K}$, such that
$\epsilon (P)\subset\gamma (N)'$. We shall write then $_\gamma\mathcal{K}_\epsilon$. If $y$ is in $P$, we
have seen that it is possible to define then the operator
$1_{\mathcal{H}}\underset{\psi}{_\beta\otimes_\gamma}\epsilon (y)$ on
$\mathcal{H}\underset{\psi}{_\beta\otimes_\gamma}\mathcal{K}$, and we define this way a
non-degenerate normal antirepresentation of $P$ on
$\mathcal{H}\underset{\psi}{_\beta\otimes_\gamma}\mathcal{K}$, we shall call again $\epsilon$ for
simplification. If $\mathcal H$ is a $Q-N$ bimodule, then
$\mathcal{H}\underset{\psi}{_\beta\otimes_\gamma}\mathcal{K}$ becomes a $Q-P$ bimodule (Connes'
fusion of bimodules).
\newline
Taking a faithful semi-finite normal weight
$\nu$  on $P$, and a left $P$-module $_{\zeta}\mathcal{L}$ (i.e. a Hilbert space $\mathcal{L}$ and a normal
non-degenerate representation $\zeta$ of $P$ on $\mathcal{L}$), it is possible then to define
$(\mathcal{H}\underset{\psi}{_\beta\otimes_\gamma}\mathcal{K})\underset{\nu}{_\epsilon\otimes_\zeta}\mathcal{L}$.
Of course, it is possible also to consider the Hilbert space
$\mathcal{H}\underset{\psi}{_\beta\otimes_\gamma}(\mathcal{K}\underset{\nu}{_\epsilon\otimes_\zeta}\mathcal{L})$.
It can be shown that these two Hilbert spaces are isomorphics as $\beta (N)'-\zeta
(P)^{'o}$-bimodules. (In ([V1] 2.1.3), the proof, given for $N=P$ abelian can be used, without
modification, in that wider hypothesis). We shall write then
$\mathcal{H}\underset{\psi}{_\beta\otimes_\gamma}\mathcal{K}\underset{\nu}{_\epsilon\otimes_\zeta}\mathcal{L}$
without parenthesis, to emphazise this coassociativity property of the relative tensor
product.
 \newline
 If $\pi$ denotes the canonical left representation of $N$ on the Hilbert space $L^2(N)$, then it is straightforward to verify that the application which sends, for all $\xi$ in $\mathcal {H}$ and $x$ in $\gN_\chi$, the vector $\xi{}_\beta\underset{\chi}{\otimes}{}_\pi J_\chi\Lambda_\chi (x)$ on $\beta(x^*)\xi$, gives an isomorphism of $\mathcal{H}{}_\beta\underset{\chi}{\otimes}{}_\pi L^2(N)$ on $\mathcal{H}$, which will send the antirepresentation of $N$ given by $n\rightarrow 1_{\mathcal{H}}{}_\beta\underset{\chi}{\otimes}{}_\pi J_\chi n^*J_\chi$ on $\beta$\vspace{2mm}. 
\newline
If $\mathcal H$ and $\mathcal K$ are finite-dimensional Hilbert spaces, the relative tensor product 
$\mathcal{H}\underset{\psi}{_\beta\otimes_\gamma}
\mathcal{K}$ can be identified with a subspace of the tensor Hilbert space 
$\mathcal{H}\otimes
\mathcal{K}$ ([EV] 2.4), the projection on which belonging to $\beta (N)\otimes\gamma (N)$.

\subsection{Fiber product [V1], [EV]} 
\label{fiber}
Let us follow the notations of \ref{rel}; let now
$M_1$ be a von Neumann algebra on $\mathcal{H}$, such that $\beta (N)\subset
M_1$, and $M_2$ be a von Neumann algebra on $\mathcal{K}$, such that $\gamma (N)\subset
M_2$. The von Neumann algebra generated by all elements $x\underset{N}{_\beta\otimes_\gamma} y$,
where
$x$ belongs to $M'_1$, and $y$ belongs $M'_2$ will be denoted
$M'_1\underset{N}{_\beta\otimes_\gamma} M'_2$ (or $M'_1\otimes_N M'_2$ if no confusion if
possible), and will be called the relative tensor product of
$M'_1$ and $M'_2$ over $N$. The commutant of this algebra will be denoted 
$M_1\underset{N}{_\beta *_\gamma} M_2$ (or $M_1*_N M_2$ if no confusion is possible) and called the
fiber product of $M_1$ and
$M_2$, over
$N$. It is straightforward to verify that, if $P_1$ and $P_2$ are two other von Neumann
algebras satisfying the same relations with $N$, we have 
\[M_1*_N M_2\cap P_1*_N P_2=(M_1\cap P_1)*_N (M_2\cap P_2)\]
Moreover, we get that $\varsigma_N (M_1\underset{N}{_\beta *_\gamma}
M_2)=M_2\underset{N^o}{_\gamma *_\beta}M_1$.
\newline
In particular, we have :
\[(M_1\cap \beta (N)')\underset{N}{_\beta\otimes_\gamma} (M_2\cap \gamma (N)')\subset
M_1\underset{N}{_\beta *_\gamma} M_2\] and :
\[M_1\underset{N}{_\beta *_\gamma} \gamma(N)=(M_1\cap\beta (N)')\underset{N}{_\beta\otimes_\gamma} 1\]
More generally, if
$\beta$ is a non-degenerate normal involutive antihomomorphism from
$N$ into a von Neumann algebra
$M_1$, and
$\gamma$ a non-degenerate normal involutive homomorphism from $N$ into a von Neumann
algebra
$M_2$, it is possible
to define, without any reference to a specific Hilbert space, a von Neumann algebra
$M_1\underset{N}{_\beta *_ \gamma}M_2$. 
\newline
Moreover, if now $\beta '$ is a non-degenerate normal involutive antihomomorphism from $N$ into
another von Neumann algebra
$P_1$,
$\gamma '$ a non-degenerate normal involutive homomorphism from $N$ into another
von Neumann algebra $P_2$, $\Phi$ a normal involutive homomorphism from $M_1$ into $P_1$ such that
$\Phi\circ\beta =\beta '$, and $\Psi$ a normal involutive homomorphism from $M_2$ into $P_2$ such that
$\Psi\circ\gamma=\gamma'$, it is possible then to define a normal involutive homomorphism (the proof
given in ([S1] 1.2.4) in the case when $N$ is abelian can be extended without modification in the
general case) :
\[\Phi\underset{N}{_\beta *_\gamma}\Psi 
: M_1\underset{N}{_\beta
*_\gamma}M_2\rightarrow P_1\underset{N}{_{\beta '}*_{\gamma '}}P_2\]
In the case when $_\gamma\mathcal{K}_\epsilon$ is a $N-P^o$ bimodule as explained in \ref{rel} and
$_\zeta\mathcal{L}$ a $P$-module, if
$\gamma (N)\subset M_2$ and $\epsilon (P)\subset M_2$, and if $\zeta (P)\subset M_3$, where $M_3$ is
a von Neumann algebra on $\mathcal{L}$, it is possible to consider then $(M_1\underset{N}{_\beta
*_\gamma}M_2)\underset{P}{_\epsilon *_\zeta}M_3$ and $M_1\underset{N}{_\beta
*_\gamma}(M_2\underset{P}{_\epsilon *_\zeta}M_3)$. The coassociativity property for relative tensor
products leads then to the isomorphism of these von Neumann algebra we shall write now 
$M_1\underset{N}{_\beta
*_\gamma}M_2\underset{P}{_\epsilon *_\zeta}M_3$ without parenthesis.
\newline
If $M_1$ and $M_2$ are finite-dimensional, the fiber product $M_1\underset{N}{_\beta *_ \gamma}M_2$
can be identified to a reduced algebra of $M_1\otimes M_2$ (reduced by a projector which belongs to
$\beta (N)\otimes \gamma (N)$). ([EV] 2.4)

\subsection{Slice maps [E2]}
\label{slice}
Let $A$ be in $M_1\underset{N}{_\beta *_\gamma}M_2$, and let $\xi_1$, $\xi_2$ be in
$D(\mathcal{H}_\beta,\psi^o)$; let us define :
\[(\omega_{\xi_1, \xi_2}*id)(A)=\lambda_{\xi_2}^*A\lambda_{\xi_1}\]
We define this way a $(\omega_{\xi_1, \xi_2}*id)(A)$ as a bounded operator on $\mathcal{K}$,
which belongs to $M_2$, such that :
\[((\omega_{\xi_1, \xi_2}*id)(A)\eta_1|\eta_2)=
(A(\xi_1\underset{\psi}{_\beta\otimes_\gamma}\eta_1)|
\xi_2\underset{\psi}{_\beta\otimes_\gamma}\eta_2)\]
One should note that $(\omega_{\xi_1, \xi_2}*id)(1)=\gamma (<\xi_1, \xi_2 >_{\beta, \psi^o})$. 
\newline
Let us define the same way, for any $\eta_1$, $\eta_2$ in
$D(_\gamma\mathcal{K}, \psi)$:
\[(id*\omega_{\eta_1, \eta_2})(A)=\rho_{\eta_2}^*A\rho_{\eta_1}\]
which belongs to $M_1$. 
\newline
We therefore have a Fubini formula for these slice maps : for any $\xi_1$, $\xi_2$ in
$D(\mathcal{H}_\beta,\psi^o)$, $\eta_1$, $\eta_2$ in $D(_\gamma\mathcal{K}, \psi)$, we have :
\[<(\omega_{\xi_1, \xi_2}*id)(A), \omega_{\eta_1, \eta_2}>=<(id*\omega_{\eta_1,
\eta_2})(A),\omega_{\xi_1, \xi_2}>\]
Equivalently ([E2] 3.3), for any $\omega_1$
in $M_{1*}^+$ such that there exists $k_1$ in $\mathbb{R}^+$ such that $\omega_1\circ\beta\leq
k_1\psi^o$, and $\omega_2$ in $M_{2*}^+$ such that there exists $k_2$ in $\mathbb{R}^+$ such that
$\omega_2\circ\gamma\leq k_2\psi$, we have :
\[<(\omega_1*id)(A), \omega_2>=<(id*\omega_2)(A), \omega_1>\]
Let $\phi_1$ be a normal semi-finite weight on
$M_1^+$, and $A$ be a positive element of the fiber product
$M_1\underset{N}{_\beta*_\gamma}M_2$, then we may define an element of the extended positive part
of $M_2$, denoted
$(\phi_1*id)(A)$, such that, for all $\eta$ in $D(_\gamma L^2(M_2), \psi)$, we have :
\[\|(\phi_1*id)(A)^{1/2}\eta\|^2=\phi_1(id*\omega_\eta)(A)\]
Moreover, then, if $\phi_2$ is a normal semi-finite weight on $M_2^+$, we have :
\[\phi_2(\phi_1*id)(A)=\phi_1(id*\phi_2)(A)\]
and if $\omega_i$ be in $M_{1*}$ such that $\phi_1=sup_i\omega_i$, we
have $(\phi_1*id)(A)=sup_i(\omega_i*id)(A)$.
\newline
Let now $P_1$ be a von Neuman algebra such that :
\[\beta(N)\subset P_1\subset M_1\]
and let $\Phi_i$ ($i=1,2$)
be a normal faithful semi-finite operator valued weight from $M_i$ to $P_i$; for any positive
operator $A$ in the fiber product
$M_1\underset{N}{_\beta*_\gamma}M_2$, there exists an element $(\Phi_1*id)(A)$
of the extended positive part
of $P_1\underset{N}{_\gamma*_\beta}M_2$, such that ([E2], 3.5), for all
$\eta$ in $D(_\gamma L^2(M_2), \psi)$, and $\xi$ in $D(L^2(P_1)_\beta, \psi^o)$, we have :
\[\|(\Phi_1*id)(A)^{1/2}(\xi\underset{\psi}{_\beta\otimes_\gamma}\eta)\|^2=
\|\Phi_1(id*\omega_\eta)(A)^{1/2}\xi\|^2\]
If $\phi$ is a normal semi-finite weight on $P$, we have :
\[(\phi\circ\Phi_1*id)(A)=(\phi*id)(\Phi_1*id)(A)\]
We define the same way an element $(id*\Phi_2)(A)$ of the extended positive part
of
$M_1\underset{N}{_\gamma*_\beta}P_2$, and we have :
\[(id*\Phi_2)((\Phi_1*id)(A))=(\Phi_1*id)((id*\Phi_2)(A))\]
 Considering now an element $x$ of $M_1{}_\beta\underset{\psi}{*}{}_\pi \pi(N)$, which can be identified 
(\ref{fiber}) to $M_1\cap\beta(N)'$, we get that, for $e$ in $\gN_\psi$, we have  \[(id_\beta\underset{\psi}{*}{}_\pi\omega_{J_\psi \Lambda_{\psi}(e)})(x)=\beta(ee^*)x\]
 Therefore, by increasing limits, we get that $(id_\beta\underset{\psi}{*}{}_\pi\psi)$ is the injection of $M_1\cap\beta(N)'$ into $M_1$.  More precisely, if $x$ belongs to $M_1\cap\beta(N)'$, we have :
 \[(id_\beta\underset{\psi}{*}{}_\pi\psi)(x{}_\beta\underset{\psi}{\otimes}{}_\pi 1)=x\]
 \newline
 Therefore, if $\Phi_2$ is a normal faithful semi-finite operator-valued weight from $M_2$ onto $\gamma(N)$, we get that, for all $A$ positive in $M_1\underset{N}{_\beta*_\gamma}M_2$, we have :
 \[(id_\beta\underset{\psi}{*}{}_\gamma\psi\circ\Phi_2)(A){}_\beta\underset{\psi}{\otimes}{}_\gamma 1=
 (id_\beta\underset{\psi}{*}{}_\gamma\Phi_2)(A)\]
 Let now $\xi$ be $D(L^2(M_1)_\beta, \psi^o)$; then, the operator $(\omega_\xi{}_\beta\underset{\psi}{*}{}_\gamma id)(A)$ is bounded and we get :
 \begin{multline*}
 \Phi_2((\omega_\xi{}_\beta\underset{\psi}{*}{}_\gamma id)(A))
 =(\omega_\xi{}_\beta\underset{\psi}{*}{}_\gamma id)(id_\beta\underset{\psi}{*}{}_\gamma\Phi_2)(A)\\
 =\gamma(<id_\beta\underset{\psi}{*}{}_\gamma\psi\circ\Phi_2)(A)\xi, \xi>_{\beta, \psi^o})
 \end{multline*}

\subsection{Hopf-bimodules}
\label{Hbimod}
A quadruplet $(N, M, r, s, \Gamma)$ will be called a Hopf-bimodule, following ([Val1], [EV] 6.5), if
$N$,
$M$ are von Neumann algebras, $r$ a faithful non-degenerate representation of $N$ into $M$, $s$ a
faithful non-degenerate anti-representation of
$N$ into $M$, with commuting ranges, and $\Gamma$ an injective involutive homomorphism from $M$
into
$M\underset{N}{_s *_r}M$ such that, for all $X$ in $N$ :
\newline
(i) $\Gamma (s(X))=1\underset{N}{_s\otimes_r}s(X)$
\newline
(ii) $\Gamma (r(X))=r(X)\underset{N}{_s\otimes_r}1$ 
\newline
(iii) $\Gamma$ satisfies the co-assosiativity relation :
\[(\Gamma \underset{N}{_s *_r} id)\Gamma =(id \underset{N}{_s *_r}\Gamma)\Gamma\]
This last formula makes sense, thanks to the two preceeding ones and
\ref{fiber}\vspace{5mm}.\newline
If $(N, M, r, s, \Gamma)$ is a Hopf-bimodule, it is clear that
$(N^o, M, s, r,
\varsigma_N\circ\Gamma)$ is another Hopf-bimodule, we shall call the symmetrized of the first
one. (Recall that $\varsigma_N\circ\Gamma$ is a homomorphism from $M$ to
$M\underset{N^o}{_r*_s}M$).
\newline
If $N$ is abelian, $r=s$, $\Gamma=\varsigma_N\circ\Gamma$, then the quadruplet $(N, M, r, r,
\Gamma)$ is equal to its symmetrized Hopf-bimodule, and we shall say that it is a symmetric
Hopf-bimodule\vspace{5mm}.\newline
Let $\mathcal G$ be a groupo\"{\i}d, with $\mathcal G^{(0)}$ as its set of units, and let us denote
by $r$ and $s$ the range and source applications from $\mathcal G$ to $\mathcal G^{(0)}$, given by
$xx^{-1}=r(x)$ and $x^{-1}x=s(x)$. As usual, we shall denote by $\mathcal G^{(2)}$ (or $\mathcal
G^{(2)}_{s,r}$) the set of composable elements, i.e. 
\[\mathcal G^{(2)}=\{(x,y)\in \mathcal G^2; s(x)=r(y)\}\]
In [Val1] was associated to a locally compact groupo\"{\i}d $\mathcal G$, equipped with a Haar system (see [R1],
[R2], [C2] II.5 and [AR] for more details, precise definitions and examples of groupo\"{\i}ds) two
Hopf-bimodules : 
\newline
The first one is $(L^\infty (\mathcal G^{(0)}), L^\infty (\mathcal G), r_{\mathcal G}, s_{\mathcal G}, \Gamma_{\mathcal
G})$, where, for $g$ in $L^\infty (\mathcal G^{(0)})$, we put $r_{\mathcal G}(g)=g\circ r$,
$s_{\mathcal G}(g)=g\circ s$, and
$\Gamma_{\mathcal G}$ is defined the following way : for $f$ in $L^\infty (\mathcal G)$, $\Gamma_{\mathcal
G}(f)$ will be the function defined on $\mathcal G^{(2)}$ by $(s,t)\rightarrow f(st)$;
$\Gamma_{\mathcal G}$ is then an involutive homomorphism from $L^\infty (\mathcal G)$ into $L^\infty
(\mathcal G^2_{s,r})$ (which can be identified to
$L^\infty (\mathcal G){_s*_r}L^\infty (\mathcal G)$).
\newline
The second one is symmetric; it is $(L^\infty (\mathcal G^{(0)}), \mathcal L(\mathcal G), r_{\mathcal G}, r_{\mathcal G},
\widehat{\Gamma_{\mathcal G}})$, where
$\mathcal L(\mathcal G)$ is the von Neumann algebra generated by the convolution algebra associated to the
groupo\"{\i}d
$\mathcal G$, and $\widehat{\Gamma_{\mathcal G}}$ has been defined in [Y1] and
[Val1]\vspace{5mm}.\newline
If $(N,M,r,s,\Gamma)$ be a Hopf-bimodule with a finite-dimensional algebra $M$, then, the
identification of $M\underset{N}{_s*_r}M$ with a reduced algebra $(M\otimes M)_e$ (\ref{fiber})
leads to an injective homomorphism $\widetilde{\Gamma}$ from $M$ to $M\otimes M$ such that
$\widetilde{\Gamma}(1)=e\not= 1$ and $(\widetilde{\Gamma}\otimes
id)\widetilde{\Gamma}=(id\otimes\widetilde{\Gamma})\widetilde{\Gamma}$ ([EV] 6.5). Then $(M,
\widetilde{\Gamma})$ is a weak Hopf
$\mathbb{C}^*$-algebra in the sense of ([BSz1], [BSz2], [Sz]).


\section{Pseudo-multiplicative unitary}
\label{pmu}
In this chapter, we recall (\ref{defmult}) the definition of a pseudo-multiplicative unitary, give the fundamental example given by groupo\"{\i}ds (\ref{gd}), and construct the Hopf-bimodules "generated by the left (resp. right) leg" of a pseudo-multiplicative unitary (\ref{AW}). Then, we recall the definition of a depth 2 inclusions (\ref{basic}), and the precise hypothesis under which, in [EV],  a pseudo-multiplicative unitary was then associated (\ref{exW}).
\subsection{Definition}
\label{defmult}
Let $N$ be a von Neumann algebra; let
$\gH$ be a Hilbert space on which $N$ has a non-degenerate normal representation $\alpha$ and two
non-degenerate normal anti-representations $\beta$ and $\hat{\beta}$. These 3 applications
are supposed to be injective, and to commute two by two.  Let $\nu$ be a normal semi-finite faithful weight on
$N$; we can therefore construct the Hilbert spaces
$\gH\underset{\nu}{_{\hat{\beta}}\otimes_\alpha}\gH$ and
$\gH\underset{\nu^o}{_\alpha\otimes_{\beta}}\gH$. A unitary $W$ from
$\gH\underset{\nu}{_{\hat{\beta}}\otimes_\alpha}\gH$ onto
$\gH\underset{\nu^o}{_\alpha\otimes_{\beta}}\gH$
will be called a pseudo-multiplicative unitary over the basis $N$, with respect to the
representation $\alpha$, and the anti-representations $\beta$ and $\hat{\beta}$, if :
\newline
(i) $W$ intertwines $\alpha$, $\beta$, $\hat{\beta}$ in the following way :
\[W(\alpha
(X)\underset{N}{_{\hat{\beta}}\otimes_\alpha}1)=
(1\underset{N^o}{_\alpha\otimes_{\beta}}\alpha(X))W\]
\[W(1\underset{N}{_{\hat{\beta}}\otimes_\alpha}\hat{\beta}
(X))=(1\underset{N^o}{_\alpha\otimes_{\beta}}\hat{\beta} (X))W\]
\[W(\beta(X) \underset{N}{_{\hat{\beta}}\otimes_\alpha}1)=
(\beta(X)\underset{N^o}{_\alpha\otimes_{\beta}}1)W\]
\[W(1\underset{N}{_{\hat{\beta}}\otimes_\alpha}\beta(X))=
(\hat{\beta}(X)\underset{N^o}{_\alpha\otimes_{\beta}}1)W\]
(ii) The operator satisfies :
\begin{multline*}
(1_{\gH}\underset{N^o}{_\alpha\otimes_{\beta}}W)
(W\underset{N}{_{\hat{\beta}}\otimes_\alpha}1_{\gH})=\\
=(W\underset{N^o}{_\alpha\otimes_{\beta}}1_{\gH})
(\sigma_{\nu^o}\underset{N^o}{_\alpha\otimes_{\beta}}1_{\gH})
(1_{\gH}\underset{N^o}{_\alpha\otimes_{\beta}}W)
\sigma_{2\nu}
(1_{\gH}\underset{N}{_{\hat{\beta}}\otimes_\alpha}\sigma_{\nu^o})
(1_{\gH}\underset{N}{_{\hat{\beta}}\otimes_\alpha}W)
\end{multline*}
In that formula, the first $\sigma_{\nu^o}$ is the relative flip defined in \ref{rel} from $\gH\underset{\nu^o}{_\alpha\otimes_{\hat{\beta}}}\gH$
to $\gH\underset{\nu}{_{\hat{\beta}}\otimes_\alpha}\gH$, and the second is the relative flip from
$\gH\underset{\nu^o}{_\alpha\otimes_{\beta}}\gH$ to $\gH\underset{\nu}{_\beta\otimes_\alpha}\gH$; while $\sigma_{2\nu}$ is the relative flip
from $\gH\underset{\nu}{_{\hat{\beta}}\otimes_\alpha}\gH\underset{\nu}{_\beta\otimes_\alpha}\gH$ to
$\gH\underset{\nu^o}{_\alpha\otimes_\beta}(\gH\underset{\nu}{_{\hat{\beta}}\otimes_\alpha}\gH)$. The index $2$ is written to recall that the
flip "turns" around the second relative tensor product, and, in such a formula, the parenthesis are written to recall that, in such a situation,
associativity rules does not occur because the anti-representation $\beta$ is here acting in the second leg of
$\gH\underset{\nu}{_{\hat{\beta}}\otimes_\alpha}\gH$. 
\newline
All the properties supposed in (i) allow us to write such a formula, which will be called the
"pentagonal relation". 
\newline
If we start from another normal semi-finite faithful weight $\nu'$ on $N$, we may define, using \ref{rel}, another unitary $W^{\nu'}=U^{\nu^o,
\nu^{'o}}_{\alpha, \beta}WU^{\nu', \nu}_{\hat{\beta}, \alpha}$ from $\gH\underset{\nu'}{_{\hat{\beta}}\otimes_\alpha}\gH$ onto
$\gH\underset{\nu^{'o}}{_\alpha\otimes_{\beta}}\gH$. The formulae which link these isomorphims between relative product Hilbert spaces and the
relative flips allow us to check that this operator $W^{\nu'}$ is also pseudo-multiplicative; which can be resumed in saying that a
pseudo-multiplicative unitary does not depend on the choice of the weight on $N$.

\subsection{Fundamental example}
\label{gd}
Let $\mathcal G$ be a measured groupo\"{\i}d, with $\mathcal G^{(0)}$ as space
of units, and $r$ and $s$ the range and source functions from $\mathcal G$ to $\mathcal G^{(0)}$. Let us
note :
\[\mathcal G^2_{r,r}=\{(x,y)\in \mathcal G^2, r(x)=r(y)\}\]
\[\mathcal G^2_{s,r}=\{(x,y)\in \mathcal G^2, s(x)=r(y)\}\]
Then, it has been shown [Val1] that the formula $W_{\mathcal G}f(x,y)=f(x,x^{-1}y)$, where $x$, $y$ are
in
$\mathcal G$, such that $r(y)=r(x)$, and $f$ belongs to $L^2(\mathcal G^2_{s,r})$ (with respect to an
appropriate measure), is a unitary from $L^2(\mathcal G^2_{s,r})$ to $L^2(\mathcal G^2_{r,r})$. Moreover, this
unitary can be interpreted [Val2] as a pseudo-multiplicative unitary over the basis
$L^\infty (\mathcal G^{(0)})$, with respect to the representation $r_{\mathcal G}$, and anti-representation
$s_{\mathcal G}$ and
$r_{\mathcal G}$ (as here the basis is abelian, the notions of representation and anti-representations are
the same, and the commutation property is fulfilled), where $r_{\mathcal G}$ and $s_{\mathcal G}$ are defined,
from
$L^\infty (\mathcal G^{(0)})$ to $L^\infty (\mathcal G)$ (and then considered as representations on $\mathcal L(L^2(\mathcal
G))$, for any
$f$ in
$L^\infty (\mathcal G^{(0)})$, by
$r_{\mathcal G}(f)=f\circ r$ and $s_{\mathcal G}(f)=f\circ s$.

\subsection{Hopf-bimodules associated to a pseudo-multiplicative unitary}
\label{AW}
For $\xi_1$ in $D(\gH_{\hat{\beta}},\nu^o)$, $\eta_1$ in $D(_\alpha\gH, \nu)$, the operator $(\lambda_{\eta_1}^{\alpha,
\beta})^*W\lambda_{\xi_1}^{\hat{\beta}, \alpha}$ will be written $(\omega_{\xi_1,
\eta_1}*id)(W)$ for ; we have,
therefore, for all
$\xi_2$,
$\eta_2$ in
$\gH$ :
\[((\omega_{\xi_1,\eta_1}*id)(W)\xi_2|\eta_2)=(W(\xi_1\underset{\nu}{_{\hat{\beta}}\otimes_\alpha}\xi_2)|
\eta_1\underset{\nu^o}{_\alpha\otimes_\beta}\eta_2)\]
and, using the intertwining property of $W$ with $\hat{\beta}$, we easily get that $(\omega_{\xi_1,
\eta_1}*id)(W)$ belongs to $\hat{\beta}(N)'$. If $\xi$ belongs to $D(\gH_{\hat{\beta}}, \nu^o)\cap D(_\alpha\gH, \nu)$, we shall write
$(\omega_\xi*id)(W)$ instead of $(\omega_{\xi, \xi}*id)(W)$. 
\newline
Following ([EV] 6.1 and 6.5), we shall write $\mathcal A(W)$ (or $\mathcal A$) the von Neumann algebra generated by  these operators. We then have $\mathcal A(W)\subset\hat{\beta}(N)'$. 
\newline
For $\xi_2$ in $D(_\alpha\gH, \nu)$, $\eta_2$ in $D(\gH_\beta, \nu^o)$, the operator $(\rho_{\eta_2}^{\alpha,
\beta})^*W\rho_{\xi_2}^{\hat{\beta}, \alpha}$ will be written $(id*\omega_{\xi_2, \eta_2})(W)$; we have, therefore, for all
$\xi_1$, $\eta_1$ in $\gH$ :
\[((id*\omega_{\xi_2, \eta_2})(W)\xi_1|\eta_1)=(W(\xi_1\underset{\nu}{_{\hat{\beta}}\otimes_\alpha}\xi_2)|
\eta_1\underset{\nu^o}{_\alpha\otimes_\beta}\eta_2)\]
and, using the intertwining property of $W$ with $\beta$, we easily get that $(id*\omega_{\xi_2, \eta_2})(W)$ belongs
to $\beta (N)'$. If $\xi$ belongs to $D(_\alpha\gH, \nu)\cap D(\gH_\beta, \nu^o)$, we shall write $(id*\omega_\xi)(W)$ instead of
$(id*\omega_{\xi, \xi})(W)$.
\newline
Following ([EV] 6.1 and 6.5), we shall write $\widehat{\mathcal A}(W)$ (or $\widehat{\mathcal A}$) the von Neumann algebra generated by  these operators.
We then have $\widehat{\mathcal A}(W)\subset\beta(N)'$.
\newline
In ([EV] 6.3 and 6.5), using the pentagonal equation, we got
that
$(N^o,\mathcal A,\beta,\alpha,\Gamma)$, and
$(N,\widehat{\mathcal A}, \alpha,\hat{\beta},  \widehat{\Gamma})$ are Hopf-bimodules, where $\Gamma$ and
$\widehat{\Gamma}$ are defined, for any $x$ in $\mathcal A$ and $y$ in $\widehat{\mathcal
A}$, by :
\[\Gamma(x)=W(x\underset{N}{_{\hat{\beta}}\otimes_\alpha}1)W^*\]
\[\widehat{\Gamma}(y)=W^*(1\underset{N^o}{_\alpha\otimes_\beta}y)W\]
In ([EV] 6.1(iv)), we had obtained that $x$ in $\mathcal L(\gH)$ belongs to $\mathcal A'$ if and only if $x$ belongs to $\alpha(N)'\cap
\beta(N)'$ and verify $(1\underset{N^o}{_\alpha\otimes_\beta}x)W=W(1\underset{N}{_{\hat{\beta}}\otimes_\alpha}x)$. We obtain the same way
that $y$ in $\mathcal L(\gH)$ belongs to $\widehat{\mathcal A}'$ if and only if $y$ belongs to $\alpha(N)'\cap
\hat{\beta}(N)'$ and verify $(y\underset{N^o}{_\alpha\otimes_\beta}1)W=W(y\underset{N}{_{\hat{\beta}}\otimes_\alpha}1)$. 
\newline
Moreover, we get that $\alpha(N)\subset\mathcal A\cap\widehat{\mathcal A}$, $\beta(N)\subset\mathcal A$,
$\hat{\beta}(N)\subset\widehat{\mathcal A}$, and, for all $x$ in $N$ :
\[\Gamma (\alpha (x))=1\underset{N^o}{_\alpha\otimes_\beta}\alpha (x)\]
\[\Gamma (\beta (x))=\beta (x)\underset{N^o}{_\alpha\otimes_\beta}1\]
\[\widehat{\Gamma}(\alpha(x))=\alpha (x)\underset{N}{_{\hat{\beta}}\otimes_\alpha}1\]
\[\widehat{\Gamma}(\hat{\beta}(x))=1\underset{N}{_{\hat{\beta}}\otimes_\alpha}\hat{\beta}(x)\]

Let us take the notations of \ref{gd}; the von Neumann algebra $\mathcal A(W_{\mathcal G})$ is equal to the von Neumann algebra $\mathcal
L(\mathcal G)$ ([Val2], 3.2.6 and 3.2.7); using ([Val2] 3.1.1), we get that the Hopf-bimodule homomorphism $\Gamma$ defined on $\mathcal L(\mathcal
G)$ by
$W_{\mathcal G}$ is the usual Hopf-bimodule homomrphism $\widehat{\Gamma_{\mathcal G}}$ studied in [Y1] and [Val1]. The
von Neumann algebra
$\widehat{\mathcal A(W_{\mathcal G})}$ is equal to the von Neumann algebra $L^{\infty}(\mathcal G, \nu)$ ([Val2], 3.2.6 and 3.2.7); using ([Val2]
3.1.1), we get that the Hopf-bimodule homomorphism
$\widehat{\Gamma}$ defined on
$L^{\infty}(\mathcal G, \nu)$ by $W_{\mathcal G}$ is equal to the usual Hopf-bimodule homomorphism $\Gamma_{\mathcal G}$ studied in [Val1], and
recalled in
\ref{gd}.

\subsection{Pseudo-multiplicative unitary associated to a depth 2 inclusion}
\label{exW}
Let $M_0\subset M_1$ be an inclusion of von Neumann algebras, equipped with a normal semi-finite faithful operator-valued weight $T_1$ from $M_1$ to $M_0$; with the notations of \ref{basic}, following ([GHJ] 4.6.4), we shall say that the inclusion $M_0\subset M_1$ is depth 2 if the inclusion $M'_0\cap M_1\subset M'_0\cap M_2\subset M'_0\cap M_3$ is standard, and, following ([EN], 11.12), we shall say that the operator-valued weight $T_1$ is regular if both restrictions $T_{2|M'_0\cap M_2}$ and $T_{3|M'_1\cap M_3}$ are semi-finite. 
\newline
In the sequel, we shall consider a depth 2 inclusion $M_0\subset M_1$ of $\sigma$-finite von Neumann algebras, equipped with a regular normal semi-finite faithful operator-valued weight $T_1$. In such a case, we had constructed and studied in ([EV]) a pseudo-multiplicative
unitary, the construction of which is here recalled. 
 \newline
We consider the set of double intertwinners
$Hom_{M_0, M_1^o}(H_1, H_2)$, and, more precisely, if $\chi$ is a normal semi-finite faithful weight on $M'_0\cap M_1$, the subset $Hom_\chi$
defined this way :
\[Hom_\chi=\{x\in Hom_{M_0, M_1^o}(H_1, H_2) / \chi (x^*x)<\infty\}\]
This set is clearly a pre-Hilbert space, and we shall denote $\gH$ its completion, and $\Lambda_\chi$ the canonical injection of $Hom_\chi$ into $\gH$. 
\newline
 For all $a\in \gN_{T_2}\cap M'_0$, 
$\Lambda_{T_2}(a)$ belongs to $Hom_{M_0, M_1^o}(H_1, H_2)$, and, for any $e$ in $\gN_\chi$,
$\Lambda_{T_2}(a)e$ belongs to $Hom_\chi$. As, by hypothesis, $\gN_{T_2}\cap M'_0$ is not reduced to $\{0\}$, $Hom_\chi$ (and $\gH$) are not reduced to
$\{0\}$. 
\newline
By hypothesis, the restriction $\tilde{T_2}$ of $T_2$ to $M'_0\cap M_2$ is semi-finite, and if we write $\chi_2=\chi\circ\tilde{T_2}$, which is a normal
semi-finite faithful weight on $M'_0\cap M_2$, we obtain this way an injection $I$ from $L^2(M'_0\cap M_2)$ into $\gH$, defined, for all $a$ in
$\gN_{\chi_2}$, by :
\[I\Lambda_{\chi_2}(a)=\Lambda_\chi(\Lambda_{T_2}(a))\]
Then, we can prove ([EV], 3.8) that this isometry $I$ is surjective, and we shall
identify $\gH$ with $L^2(M'_0\cap M_2)$. 
 \newline
By hypothesis, the inclusion $M'_0\cap M_1\subset M'_0\cap
M_2\subset M'_0\cap M_3$ is a basic construction, which means that there exist a normal faithful representation $\pi$ of $M'_0\cap M_3$ on the
Hilbert space $L^2(M'_0\cap M_2)$, such that $\pi(M'_0\cap M_3)=J_{\chi_2}\pi_{\chi_2}(M'_0\cap M_1)'J_{\chi_2}$. The restriction of $\pi$ to
$M'_0\cap M_2$ is $\pi_{\chi_2}$; moreover, $\pi$ can be easily described using the identification of $L^2(M'_0\cap M_2)$ with
$\gH$ ([EV], 32(ii), 3.9, 3.10); we have, for all
$X$ in $M'_0\cap M_3$ and $x$ in $Hom_\chi$ :
\[\pi(X)\Lambda_\chi (x)=\Lambda_\chi (Xx)\]
 Moreover, we get a normal semi-finite faithful operator-valued weight $\tilde{T_3}$ from $M'_0\cap M_3$ to $M'_0\cap M_2$, and a normal semi-finite faithful weight $\chi_3$ on $M'_0\cap M_3$, such that, for all $x$ in $M'_0\cap M_3$, we have $\pi(\sigma_t^{\chi_3}(x))=\Delta_{\chi_2}^{it}\pi(x)\Delta_{\chi_2}^{-it}$.
 \newline
In that situation, we had constructed in ([EV]) a pseudo-multiplicative
unitary over the basis
$(M'_0\cap M_1)^o$, with respect to a representation $s$, and two antirepresentations $r$ and $\hat{r}$ of $(M'_0\cap
M_1)^o$ on the Hilbert space $\gH=L^2(M'_0\cap M_2)$.
\newline
Here, $r$ is the restriction of $\pi_{\chi_2}$ to $M'_0\cap M_1$, $\hat{r}$ is
the isomorphism of $M'_0\cap M_1$ onto $M'_2\cap M_3$ given by $j_2\circ j_1$, composed with the restriction of $\pi$ to
$M'_2\cap M_3$, and $s$ is the antirepresentation of $M'_0\cap M_1$ given, for all $x$ in
$M'_0\cap M_1$, by $s(x)=J_{\chi_2}x^*J_{\chi_2}$, which sends $M'_0\cap M_1$ onto $J_{\chi_2}\pi_{\chi_2}(M'_0\cap
M_1)J_{\chi_2}$, which, by the depth 2 hypothesis, is equal to $\pi(M'_0\cap M_3)'$. 
\newline
The operator $W$ can be defined the following way ([EV], 5.3) : for any $x$ in $Hom_\chi$, $\Lambda_\chi (x)$ belongs to $D(\gH_s, \chi^o)$
([EV],4.2), and for $x$, $y$, $z_\alpha$, $t_\alpha$ in $Hom_\chi$ such that :
\[W(\Lambda_\chi (x)\underset{\chi^o}{_r\otimes_s}\Lambda_\chi (y))=\sum_\alpha (\Lambda_\chi
(z_\alpha)\underset{\chi}{_s\otimes_{\hat{r}}}\Lambda_\chi (t_\alpha))\]
means that :
\[(y\otimes_{M_0}1)x=\sum_\alpha (1\otimes_{M_0}z_\alpha)t_\alpha\]
the sum being strongly convergent in $Hom(H_1, H_3)$ (let us recall, by \ref{basic}, that $y\otimes_{M_0}1\in Hom(H_2, H_3)$, and that
$1\otimes_{M_0}z_\alpha\in Hom(H_2, H_3)$).
 \newline
 In [EV] is shown that the von Neumann algebra $\mathcal A(W)$ generated by the operators $(\omega_{\xi, \eta}*id)(W)$,  with $\xi$ in $D({}_r\gH, \chi)$, and $\eta$ in $D(\gH_s, \chi^o)$, is then equal to
$\pi(M'_0\cap M_2)'$, and that the bimodule homomorphism defined for $x$ in $\pi(M'_0\cap M_2)'$ by :
 \[\Gamma(x)=W(x{}_r\underset{\chi^o}{\otimes}{}_s1)W^*\]
  sends this algebra to $\pi(M'_0\cap M_2)'\underset{M'_0\cap
M_1}{_s*_{\hat{r}}}\pi(M'_0\cap M_2)'$. 
  \\It is shown also that the von Neumann algebra
$\widehat{\mathcal A(W)}$ generated by the operators $(id*\omega_{\xi', \eta'})(W)$, with $\xi'$ in $D(\gH_s, \chi^o)$ and $\eta'$ in $D({}_{\hat{r}}\gH, \chi)$, is then equal to $\pi(M'_1\cap M_3)'$, and the bimodule homomorphism defined for $x$ in $\pi(M'_1\cap M_3)'$ by :
\[\widehat{\Gamma}(y)=W^*(1 {}_s\underset{\chi}{\otimes}{}_{\hat{r}}y)W\]
  sends this algebra to $\pi(M'_1\cap M_3)'\underset{(M'_0\cap M_1)^o}{_r*_s}\pi(M'_1\cap M_3)'$.

\section{Modular theory on the basis}
\label{modular}
For the construction of the pseudo-multiplicative unitary made in \ref{exW}, we had made a technical use of some auxilliary weight on the basis; we now suppose a modular property on this weight, in order to get more analytical properties. 
\subsection{Invariant property for a weight on the basis}
\label{inv}
Let $M_0\subset M_1$ be a depth 2 inclusion of $\sigma$-finite von Neumann algebras, equipped with a regular normal semi-finite faithful operator-valued weight $T_1$. We shall use all the notations described in \ref{basic} and in \ref{exW}. We recall that, by definition, the modular automorphism group $\sigma_t^{T_1}$ is the restriction on $M'_0\cap M_1$ of $\sigma_t^{\psi_1}$; let now $\chi$ be a normal semi-finite faithful weight on $M'_0\cap M_1$; let us climb one step in the tower: we define the normal faithful semi-finite weight $\chi^2$ on $M'_1\cap M_2$ by $\chi^2=\chi\circ j_1$. We have, for all $x$ in $M'_0\cap M_1$ :
\[\sigma_s^{\chi^2}(j_1(x))=j_1\circ\sigma_{-s}^{\chi}(x)\]
and :
\begin{multline*}
\sigma_t^{T_2}(j_1(x))=\sigma_t^{\psi_2}(j_1(x))=\Delta_1^{it}J_1x^*J_1\Delta_1^{-it}\\
=J_1\Delta_1^{it}x^*\Delta_1^{-it}J_1=j_1(\sigma_t^{\psi_1}(x))=j_1(\sigma_t^{T_1}(x))
\end{multline*}
We shall now suppose that there exists some normal semi-finite faithful weight $\chi$ on $M'_0\cap M_1$ which is invariant under the modular automorphism group of $T_1$; we then get that $\chi^2$ is invariant under $\sigma_t^{T_2}$. 

\subsection{Proposition}
\label{propinv}
{\it Let $M_0\subset M_1$ be a depth 2 inclusion of $\sigma$-finite von Neumann algebras, equipped with a regular normal semi-finite faithful operator-valued weight $T_1$. Let us use all the notations of \ref{basic} and \ref{exW}; let $\chi$ be a normal semi-finite faithful weight on $M'_0\cap M_1$, invariant under the modular automorphism group of $T_1$. Then :
\newline 
(i) the weight $\chi_2$ on $M'_0\cap M_2$ is invariant under the modular automorphism group $\sigma_t^{T_1\circ T_2}$, and there exists a positive invertible operator $h$ on $\gH$, such that, for all $a$ in $\gN_{\chi_2}$ and $x$ in $M'_0\cap M_3$ :
\[h^{it}\Lambda_{\chi_2}(a)=\Lambda_{\chi_2}(\sigma_t^{\psi_2}(a))\]
\[h^{it}\pi(x)h^{-it}=\sigma_t^{\psi_3}(x)\]
The one-parameter unitary group $h^{it}$ is the canonical implementation of $\sigma_t^{T_1\circ T_2}$, and, therefore, the operator $h$ verifies $J_{\chi_2}hJ_{\chi_2}=h^{-1}$, and $h$ commutes with $\Delta_{\chi_2}$; for any $x$ in $\mathcal A(W)=\pi(M'_0\cap M_2)'$, let us put $\rho_t(x)=h^{it}xh^{-it}$; then $\rho_t$ is a one-parameter group of automorphisms of $\mathcal A (W)$; if $\chi'_2$ denotes the canonical normal faithful semi-finite weight on $\mathcal A(W)$ constructed by the Tomita-Takesaki theory from the weight $\chi_2$ on the commutant $\pi(M'_0\cap M_2)$, then, for all $s$, $t$ in $\Bbb{R}$, we have $\rho_t\circ\sigma_s^{\chi'_2}=\sigma_s^{\chi'_2}\circ\rho_t$.
\newline
(ii) for all $t$ in $\Bbb{R}$, the application which sends, for all $a$ in $\gN_{\chi_2}\cap\gN_{T_2}$ and $b$ in $\gN_{\psi_1}$, $\Lambda_{\chi_2}(a)_s\underset{\chi}{\otimes}{}_r\Lambda_{\psi_1}(b)$ on $h^{it}\Lambda_{\chi_2}(a)_s\underset{\chi}{\otimes}{}_r\Delta_{\psi_1}^{it}\Lambda_{\psi_1}(b)$ is well defined, and can be extended to a unitary on the Hilbert space $\gH_s\underset{\chi}{\otimes}{}_rH_1$ that we shall denote by 
$h^{it}_s\underset{\chi}{\otimes}{}_r\Delta_{\psi_1}^{it}$; we obtain this way a one-parameter group of unitaries on this Hilbert space, whose generator will be denoted by $h_s\underset{\chi}{\otimes}{}_r\Delta_{\psi_1}$; this group of unitaries satisfies} 
\[U_1(h^{it}_s\underset{\chi}{\otimes}{}_r\Delta_{\psi_1}^{it})=\Delta_{\psi_2}^{it}U_1\]

\begin{proof}
The proof of [E2], 4.4 (i) and (iii) remains valid in that context. \end{proof}

\subsection{Definitions and notations ([EN], 10.11, 10.13, 10.14)}
Let $X$ be in $End_{M_0^o}(H_0, H_1)$ such that $XX^*$ belongs to $\gM_{T_2}^+$; then, there exists a unique element $\Phi_1(X)$ in $M_1$ such that $T_2(X\Lambda_{T_1}(a)^*)=\Phi_1(X)a^*$, for all $a$ in $\gN_{T_1}$. We define this way an injective morphism $\Phi_1$ of $(M_1, M_0)$-bimodules such that :
\[\Phi_1(X)\Phi_1(X)^*\leq T_2(XX^*)\]
and we have clearly $\Phi_1(\Lambda_{T_1}(a))=a$, for all $a$ in $\gN_{T_1}$. 
\newline
If $X$ belongs to $\gN_{T_2}$ and $e$ to $\gN_{T_1}\cap\gN_{\psi_1}$, then $\Phi_1(X^*\Lambda_{T_1}(e))$ belongs to $\gN_{\psi_1}$, and we get :
\[X^*\Lambda_{\psi_1}(e)=\Lambda_{\psi_1}(\Phi_1(X^*\Lambda_{T_1}(e)))\]
If $X$ belongs to $\gN_{T_2}$ and $e$ to $\gN_{T_1}$, then $\Phi_1(X^*\Lambda_{T_1}(e))$ belongs to $\gN_{T_1}$, and we get :
\[X^*\Lambda_{T_1}(e)=\Lambda_{T_1}(\Phi_1(X^*\Lambda_{T_1}(e)))\]
Starting from the inclusion $M_1\subset M_2$, we get an application $\Phi_2$.

\subsection{Lemma}
\label{lemma1}
{\it Let $x$ be in $M'_1\cap \gN_{T_3}$, analytical with respect to $\sigma_t^{\chi_3}$, $a$ in $\gN_{T_2}\cap \gN_{T_2}^*\cap\gN_{\chi_2}\cap \gN_{\chi_2}^*$; then $\Phi_2(x^*\Lambda_{T_2}(a^*))$ belongs to $\gN_{\chi_2}\cap\gN_{\chi_2}^*$, and we have :}
\[\Lambda_{\chi_2}(\Phi_2(x^*\Lambda_{T_2}(a^*)))=\pi(x^*)\Lambda_{\chi_2}(a^*)\]
\[\Lambda_{\chi_2}(\Phi_2(x^*\Lambda_{T_2}(a^*))^*)=J_{\chi_2}\pi(\sigma_{-i/2}^{\chi_3}(x^*))J_{\chi_2}\Lambda_{\chi_2}(a)\]
\begin{proof}
The first part had been obtained in ([E2], 5.2), and the second part is easy, using the fact that, for all $x$ in $M'_0\cap M_3$, we have $\pi(\sigma_t^{\chi_3}(x))=\Delta_{\chi_2}^{it}\pi(x)\Delta_{\chi_2}^{-it}$ (\ref{exW}).  \end{proof}

\subsection{Lemma}
\label{lemma2}
{\it Let $x$ be in $M'_1\cap \gN_{T_3}$, analytical with respect to $\sigma_t^{\psi_3}$, $a$ in $\gN_{T_2}\cap \gN_{T_2}^*\cap\gN_{\psi_2}\cap \gN_{\psi_2}^*$, $b$ in $\gN_{T_2}\cap \gN_{T_2}^*$; then :
\newline
(i) $\Phi_2(x^*\Lambda_{T_2}(a^*))$ belongs to $\gN_{\psi_2}\cap\gN_{\psi_2}^*$, and we have :
\[\Lambda_{\psi_2}(\Phi_2(x^*\Lambda_{T_2}(a^*))^*)=J_2\sigma_{-i/2}^{\psi_3}(x^*)J_2\Lambda_{\psi_2}(a)\]
(ii) $\Phi_2(x^*\Lambda_{T_2}(b^*))$ belongs to $\gN_{T_2}\cap\gN_{T_2}^*$, and we have :}
\[\Lambda_{T_2}(\Phi_2(x^*\Lambda_{T_2}(b^*))^*)=J_2\sigma_{-i/2}^{\psi_3}(x^*)J_2\Lambda_{T_2}(b)\]
\begin{proof}
The proof had been given in ([E2], 5.3). \end{proof}

\subsection{Lemma}
\label{lemma3}
{\it Let $x$ be in $\gM_{T_3}\cap M'_1$; let us define the elements $x_n$ by :
\[x_n=\frac{n}{\pi}\int_{-\infty}^{\infty}\int_{-\infty}^{\infty} e^{-n(t^2+s^2)}\sigma_t^{\psi_3}\circ\sigma_s^{\chi_3}(x)dsdt\]
Then, for all $n$ in $\Bbb{N}$, $x_n$ belongs to $\gM_{T_3}\cap M'_1$, are analytic with respect both to $\sigma_t^{\psi_3}$ and $\sigma_s^{\chi_3}$, and the sequence $x_n$ is strongly converging to $x$.  Moreover, for all $z$ in $\Bbb{C}$, $\sigma_z^{\psi_3}(x_n)$ belongs to $\gM_{T_3}$ and is analytic with respect to $\sigma_t^{\chi_3}$.}
\begin{proof}
Let us first remark that, as the restriction of $\sigma_t^{\psi_3}$ to $M_1$ is equal to $\sigma_t^{\psi_1}$, it is clear that $\sigma_t^{\psi_3}$ leaves $M'_1\cap M_3$ globally invariant; using \ref{propinv}, we have, for all $x$ in $M'_1\cap M_3$ :
\[\pi(\sigma_t^{\psi_3}(x))=h^{it}\pi(x)h^{-it}\]
and, using \ref{exW}, we have :
\[\pi(\sigma_t^{\chi_3}(x))=\Delta_{\chi_2}^{it}\pi(x)\Delta_{\chi_2}^{-it}\]
As, by \ref{propinv}, $h$ and $\Delta_{\chi_2}$ commute, we see that the two automorphism groups $\sigma_t^{\psi_3}$ and $\sigma_s^{\chi_3}$ of $M'_1\cap M_3$ commute. 
\newline
Let us suppose first that $x$ is positive; for any $u$ in $\Bbb{R}$, we get that :
\[\sigma_u^{\psi_3}(x_n)=\frac{n}{\pi}\int\int 
e^{-n((t-u)^2+s^2)}\sigma_t^{\psi_3}\circ\sigma_s^{\chi_3}(x)dsdt\]
which is the restriction of the analytic function 
\[z\rightarrow\frac{n}{\pi}\int\int 
e^{-n((t-z)^2+s^2)}\sigma_t^{\psi_3}\circ\sigma_s^{\chi_3}(x)dsdt\]
So, $x_n$ is analytic with respect to 
$\sigma_t^{\psi_3}$. Using the commutation property of the automorphism groups, we get the same way that $x_n$ is analytic with respect to $\sigma_s^{\chi_3}$. The fact that $\sigma_z^{\psi_3}(x_n)$ is analytic with respect to $\sigma_s^{\chi_3}$ is easy to get by the same type of calculations. We obtain that it is a linear combination of four elements in $\gM_{T_3}^+$, using same arguments than in ([EN], 10.12). We obtain that $x_n$ is strongly converging to $x$ by same arguments as in ([EN], 10.12); by linearity, all the results remain true for any $x$ be in $\gM_{T_3}\cap M'_1$. \end{proof}

\subsection{Proposition}
\label{propj}
{\it Let $M_0\subset M_1$ be a depth 2 inclusion of $\sigma$-finite von Neumann algebras, equipped with a regular normal semi-finite faithful operator-valued weight $T_1$. Let us use all the notations of \ref{basic} and \ref{exW}; let $\chi$ be a normal semi-finite faithful weight on $M'_0\cap M_1$, invariant under the modular automorphism group of $T_1$. Then :
\newline
(i) for all $x$ in $M'_1\cap M_3$, we have :
\[\pi(j_2(x))=J_{\chi_2}\pi(x)^*J_{\chi_2}\]
(ii)
the operator $h$ defined in \ref{propinv} is such that the positive invertible operator $h^{-1}\Delta_{\chi_2}$ is affiliated to $\pi(M'_1\cap M_3)'$. }
\begin{proof}
Let us take $x$ in $\gN_{T_3}\cap M'_1$, analytical with respect both to $\sigma_t^{\psi_3}$ and $\sigma_s^{\chi_3}$,and such that  $\sigma_z^{\psi_3}(x)$ belongs to $\gN_{T_3}\cap M'_1$ and is analytic with respect to $\chi_3$ (does exist by \ref{lemma3}). 
\newline
Then, if $a$ belongs to $\gN_{\tilde{T_2}}\cap \gN_{\tilde{T_2}}^*\cap\gN_{\chi_2}\cap\gN_{\chi_2}^*$, we have, by \ref{lemma1} :
\[\Lambda_{\chi_2}(\Phi_2(x^*\Lambda_{T_2}(a^*))^*)=J_{\chi_2}\pi(\sigma_{-i/2}^{\chi_3}(x^*))J_{\chi_2}\Lambda_{\chi_2}(a)\]
On the other hand, using \ref{lemma2}, we see that $\Phi_2(x^*\Lambda_{T_2}(a^*))$ belongs to $\gN_{T_2}\cap \gN_{T_2}^*$, and we get that :
\begin{eqnarray*}
\Lambda_{\chi_2}(\Phi_2(x^*\Lambda_{T_2}(a^*))^*)&=
&\Lambda_\chi(\Lambda_{T_2}(\Phi_2(x^*\Lambda_{T_2}(a^*))^*)\\
&=&\Lambda_\chi(J_2\sigma_{-i/2}^{\psi_3}(x^*)J_2\Lambda_{T_2}(a))\\
&=&\pi(j_2(\sigma_{i/2}^{\psi_3}(x)))\Lambda_\chi(\Lambda_{T_2}(a))\\
&=&\pi(j_2(\sigma_{i/2}^{\psi_3}(x)))\Lambda_{\chi_2}(a)
\end{eqnarray*}
from which we get that :
\[J_{\chi_2}\pi(\sigma_{-i/2}^{\chi_3}(x^*))J_{\chi_2}\Lambda_{\chi_2}(a)=
\pi(j_2(\sigma_{i/2}^{\psi_3}(x)))\Lambda_{\chi_2}(a)\]
which, by density, leads to :
\[J_{\chi_2}\pi(\sigma_{-i/2}^{\chi_3}(x^*))J_{\chi_2}=
\pi(j_2(\sigma_{i/2}^{\psi_3}(x)))\]
Using again \ref{lemma3}, it is possible to find elements such that we can apply this formula to the element $y=\sigma_{i/2}^{\psi_3}(x)$, which leads to :
\[J_{\chi_2}\pi(j_2(y^*))J_{\chi_2}=\pi(\sigma_{i/2}^{\chi_3}\circ\sigma_{-i/2}^{\psi_3}(y))\]
from which we get that :
\begin{eqnarray*}
\pi(\sigma_{-i/2}^{\chi_3}\circ\sigma_{i/2}^{\psi_3}(y^*))
&=&\pi(\sigma_{i/2}^{\chi_3}\circ\sigma_{-i/2}^{\psi_3}(y))^*\\
&=&J_{\chi_2}\pi(j_2(y))J_{\chi_2}\\
&=&\pi(\sigma_{i/2}^{\chi_3}\circ\sigma_{-i/2}^{\psi_3}(y^*))
\end{eqnarray*}
which leads to the fact that $h^{-1}\Delta_{\chi_2}$ commutes with all the operators $\pi(y)$, with $y$ in $\gM_{T_3}\cap M'_1$, analytic with respect both to $\sigma_t^{\chi_3}$ and $\sigma_s^{\psi_3}$, such that the elements $\sigma_z^{\psi_3}(y)$ are analytic with respect to $\sigma_t^{\chi_3}$, and belong to $\gM_{T_3}$; using again \ref{lemma3}, we get that it will commute with all elements in $\gM_{T_3}\cap M'_1$, which is, by hypothesis, a dense subset of $M_3\cap M'_1$. This finishes the proof of (ii). Then, we get 
\[\pi(j_2(x))=J_{\chi_2}\pi(x)^*J_{\chi_2}\]
for all $x$ in $\gM_{T_3}\cap M'_1$, analytic with respect both to $\sigma_t^{\chi_3}$ and $\sigma_s^{\psi_3}$, such that the elements $\sigma_z^{\psi_3}(x)$ are analytic with respect to $\sigma_t^{\chi_3}$, and belong to $\gM_{T_3}$, using \ref{lemma3} again, by density, it is true for all $x$ in $\gM_{T_3}\cap M'_1$, and, by density again, we obtain (i). \end{proof}

\subsection{Corollaries}
\label{corj}
{\it With the hypothesis of \ref{propj}, we have :
\newline
(i) for any $x$ in $\pi (M'_1\cap M_3)'$ and $t$ in $\mathbb{R}$, let us put $\widehat{\rho_t}(x)=h^{it}xh^{-it}$, and $\hat{j}(x)=J_{\chi_2}x^*J_{\chi_2}$; then $\widehat{\rho_t}$ is a one-parameter automorphism group of $\pi (M'_1\cap M_3)'$, and $\hat{j}$ is an anti-isomorphism of $\pi (M'_1\cap M_3)'$, such that $\hat{j}\circ r=s$, and :
\[\hat{j}\circ\widehat{\rho_t}=\widehat{\rho_t}\circ\hat{j}\]
(ii) We have, with the notations of \ref{exW} :
\[\pi (M'_0\cap M_2)'\cap\pi (M'_1\cap M_3)'=\pi(M'_0\cap M_3)'=s(M'_0\cap M_1)\]
\[\pi (M'_0\cap M_2)\cap\pi (M'_1\cap M_3)'=\pi(M'_0\cap M_1)=r(M'_0\cap M_1)\]
\[\pi (M'_0\cap M_2)'\cap\pi (M'_1\cap M_3)=\pi(M'_2\cap M_3)=\hat{r}(M'_0\cap M_1)\]
\[\pi (M'_0\cap M_2)\cap\pi (M'_1\cap M_3)=\pi(M'_1\cap M_2)=\pi\circ j_1(M'_0\cap M_1)\]
(iii) We get :
\[M'_0\cap M_1=M'_0\cap M_2\cap (M'_1\cap M_3)'\]
\[M'_1\cap M_2=(M'_0\cap M_2)\cap (M'_1\cap M_3)\]
\[M'_2\cap M_3=M'_1\cap M_3\cap(M'_0\cap M_2)'\]
\[((M'_0\cap M_2)\cup(M'_1\cap M_3))''=M'_0\cap M_3\]
(iv) For all $t$ in $\mathbb{R}$, we have $\sigma_t^{T_2}=\sigma^{\tilde{T_2}}_{t|(M'_1\cap M_2)}$, where $\tilde{T_2}$ is the restriction of $T_2$ to $M'_0\cap M_2$, which is, by definition, a normal semi-finite faithful operator-valued weight from $M'_0\cap M_2$ to $M'_0\cap M_1$, whose modular group $\sigma_t^{\tilde{T_2}}$ is a one-parameter automorphism group of $M'_0\cap M_2\cap(M'_0\cap M_1)'$ (which contains clearly $M'_1\cap M_2$). }
\begin{proof}
For (i), (ii) and (iii), we can follow ([E2], 5.5). For any $x$ in $M'_1\cap M_2$, we have :
\[\pi(\sigma_t^{T_2}(x))=\pi(\sigma_t^{\psi_2}(x))=h^{it}\pi(x)h^{-it}\]
\[\pi(\sigma_t^{\tilde{T_2}}(x))=\pi(\sigma_t^{\chi_2}(x))=\Delta_{\chi_2}^{it}\pi(x)\Delta_{\chi_2}^{-it}\]
which gives the result, by \ref{propj}(ii). \end{proof}

\section{Back to W}
\label{bW}
In this section, following what remains of ([E2] 6 and 7) under the wider hypothesis of the existence of a normal faithful semi-finite weight $\chi$ on the basis, invariant under the modular group $\sigma_t^T$, we construct a co-inverse, a right-invariant operator-valued weight (and, therefore, a left-invariant operator-valued weight) on the Hopf-bimodule $\mathcal A(W)$. All these results are written down in \ref{thbimod}, and give (with a wider hypothesis) the proof of the results claimed in [E2]. 
\subsection{Proposition}
\label{propbeta}
{\it For any $x$ in $M'_0\cap M_2$, let us define :
\[\beta (\pi (x))=W^*(\pi (x)_s\underset{M'_0\cap M_1}{\otimes}{}_r1)W\]
Then :
\newline
(i) we have : $\beta(\pi (x))=\varsigma(id_s\underset{M'_0\cap M_1}{*}{}_r)(U^*j_2j_1(x)U)$
\newline
(ii) $\beta(\pi (x))$ belongs to $\pi(M'_0\cap M_2)_r\underset{(M'_0\cap M_1)^o}{*}{}_s \pi(M'_0\cap M_2)'$
\newline
(iii) we have :
\[(\beta_r\underset{(M'_0\cap M_1)^o}{*}{}_s id)\beta(\pi (x))=
(id_r\underset{(M'_0\cap M_1)^o}{*}{}_s \varsigma\Gamma)\beta(\pi (x))\]
\newline
(iv) for all $t$ in $\mathbb{R}$, we have :
\[\beta(\pi(\sigma_t^{\psi_2}(x)))=(\sigma_t^{\psi_2}{}_s\underset{M'_0\cap M_1}{*}{}_r Adh^{it})\beta(\pi (x))\]
\newline
(v) if $x$ is positive, we have :
\[\hat{r}(\tilde{T_2}(x))=(\chi_2 {}_r\underset{(M'_0\cap M_1)^o}{*}{}_s id)\beta(\pi (x))=
(\tilde{T_2}{}_r\underset{(M'_0\cap M_1)^o}{*}{}_s id)\beta(\pi (x))\]
\newline
(vi) if $x$ belongs to $\gN_{\tilde{T_2}}$ (resp. $\gN_{\chi_2}$), and $\xi$, $\eta$ to $D(\gH_s, \chi^o)$, the operator $(id*\omega_{\xi, \eta})\beta(\pi (x))$ belongs to $\gN_{\tilde{T_2}}$ (resp. $\gN_{\chi_2}$).  
\newline
(vii) the subspace $D(_{\hat{r}}\gH, \chi)\cap D(\gH_s, \chi^o)$ is dense in $\gH$.
\newline
(viii) if $x$ belongs to $\gN_{\chi_2}$), $\xi$ to $D(_{\hat{r}}\gH, \chi)\cap D(\gH_s, \chi^o)$ and $\eta$ to $D(\gH_s, \chi^o)$, we have :}
\[\Lambda_{\chi_2}((id*\omega_{\xi, \eta})\beta(\pi (x)))=(id*\omega_{\xi, \eta})(W^*)\Lambda_{\chi_2}(x)\]

\begin{proof}
The proof is identical to ([E2], 6.1 to 6.6). \end{proof}

\subsection{Proposition}
\label{hW}
{\it It is possible to define a one-parameter unitary group, denoted $h^{it}_r\underset{\chi^o}{\otimes}{}_sh^{it}$, on $\gH_r\underset{\chi^o}{\otimes}{}_s\gH$, with natural values on elementary tensors, and one-parameter unitary group, denoted $h^{it}_s\underset{\chi}{\otimes}{}_{\hat{r}}h^{it}$, on $\gH_s\underset{\chi}{\otimes}{}_{\hat{r}}\gH$, with natural values on elementary tensors, and they verify, for all $t$ in $\mathbb{R}$ :}
\[W(h^{it}_r\underset{\chi^o}{\otimes}{}_sh^{it})=(h^{it}_s\underset{\chi}{\otimes}{}_{\hat{r}}h^{it})W\]

\begin{proof}
As, for any $x$ in $M'_0\cap M_1$, we have, using \ref{propinv} :
\[h^{it}r(x)h^{-it}=r(\sigma_t^{\psi_1}(x))\]
\[h^{it}s(x)h^{-it}=s(\sigma_t^{\psi_1}(x))\]
\[h^{it}\hat{r}(x)h^{-it}=\sigma_t^{\psi_3}(\hat{r}(x))=\hat{r}(\sigma_t^{\psi_1}(x))\]
which, by standard computations, will give a meaning to $h^{it}_r\underset{\chi^o}{\otimes}{}_sh^{it}$ and $h^{it}_s\underset{\chi}{\otimes}{}_{\hat{r}}h^{it}$. Moreover, using again \ref{propinv} and \ref{propbeta}, we have, for all $a$ in $\gN_{\tilde{T_2}}\cap\gN_{\chi_2}$, $\xi$ in $D(_{\hat{r}}\gH, \chi)\cap D(\gH_s, \chi^o)$ and $\eta$ in $D(\gH_s, \chi^o)$ :
\begin{eqnarray*}
(id*\omega_{\xi, \eta})(W^*)h^{it}\Lambda_{\chi_2}(a)
&=&(id*\omega_{\xi, \eta})(W^*)\Lambda_{\chi_2}(\sigma_t^{\psi_2}(a))\\
&=&\Lambda_{\chi_2}((id*\omega_{\xi, \eta})\beta(\pi (\sigma_t^{\psi_2}(a))))\\
&=&\Lambda_{\chi_2}(\sigma_t^{\psi_2}((id*\omega_{h^{-it}\xi, h^{-it}\eta})\beta(\pi (a))))\\
&=&h^{it}(id*\omega_{h^{-it}\xi, h^{-it}\eta})(W^*)\Lambda_{\chi_2}(a)
\end{eqnarray*}
from which we get the result. \end{proof}

\subsection{Proposition}
\label{tau}
{\it For all $t$ in $\mathbb{R}$ and $x$ in $\pi(M'_1\cap M_3)'$, let us define $\widehat{\tau_t}(x)=\Delta_{\chi_2}^{it}x\Delta_{\chi_2}^{-it}$; then, $\widehat{\tau_t}$ is a one-parameter automorphism group of $\pi(M'_1\cap M_3)'$, which commutes with the anti-isomorphism $\hat{j}$  defined in \ref{corj} (i). Moreover, we have, for all $a$ in $M'_0\cap M_1$ :}
\[\widehat{\tau_t}(r(a))=r(\sigma_t^{\chi}(a))\]
\[\widehat{\tau_t}(s(a))=s(\sigma_t^{\chi}(a))\]
\begin{proof}
Clearly, $\widehat{\tau_t}$ is the composition of the one-parameter automorphism group $\widehat{\rho_t}$ of $\pi(M'_1\cap M_3)'$, defined in \ref{corj} (i), with the inner automorphism group implemented by $h^{-it}\Delta_{\chi_2}^{it}$. The commutation relation between $J_{\chi_2}$ and $\Delta_{\chi_2}^{it}$ gives all the other results. \end{proof}

\subsection{Theorem}
\label{thtau}
{\it Let $\xi$,$\eta$ be in $D(_{\hat{r}}\gH, \chi)\cap D(\gH_s, \chi^o)$. Then, the operator $(id*\omega_{\xi, \eta})(W^*)$ belongs to $D(\widehat{\tau_{-i/2}})$, and we have :}
\[\hat{j}\circ\widehat{\tau_{-i/2}}((id*\omega_{\xi, \eta})(W^*))=(id*\omega_{\xi, \eta})(W)\]

\begin{proof}
We can repeat what was written in ([E2], 6.8). \end{proof}

\subsection{Theorem}
\label{thtau2}
{\it There exists an involutive antilinear isomorphism $\hat{J}$ on $\gH$, which implements an anti-isomorphism $j$ of $\pi(M'_0\cap M_2)'$ such that $j\circ\hat{r}=s$, and a positive self-adjoint invertible operator $\widehat{\Delta}$ on $\gH$, which implements a one-parameter automorphism group $\tau_t$ of $\pi(M'_0\cap M_2)'$ commuting with $j$, such that, for all $\xi$,$\eta$ in $D(_r\gH, \chi)\cap D(\gH_s, \chi^o)$, the operator $(\omega_{\xi, \eta}*id)(W)$ belongs to $D(\tau_{-i/2})$, and we have : }
\[j\circ\tau_{-i/2}((\omega_{\xi, \eta}*id)(W))=(\omega_{\xi, \eta}*id)(W^*)\]
{\it Moreover, for all $a$ in $M'_0\cap M_1$, we have :}
\[\tau_t(s(a))=s(\sigma_{t}^{\chi}(a))\]
\[\tau_t(\hat{r}(a))=\hat{r}(\sigma_{t}^{\chi}(a))\]
\begin{proof}
As we had obtained \ref{thtau} under the hypothesis that there exists a weight $\chi$ invariant under $\sigma_t^{T_1}$, we can (\ref{inv}) apply this result to the inclusion $M_1\subset M_2$.  
\\ Let us precise the notations: we start from the set of double intertwiners $Hom_{M_1, M_2^o}(H_2, H_3)$, and the normal semi-finite faithful weight $\chi^2=\chi\circ j_1$ on $M'_1\cap M_2$, from which we construct (\ref{exW}) an Hilbert space $\gH_2$; more precisely, if $X$ belongs to $Hom_{M_0, M_1^o}(H_1, H_2)$, the application $X\rightarrow J_3(1_{H_1}\otimes_{\psi_0}X)J_2$ defines an antilinear isomorphism from $Hom_{M_0, M_1^o}(H_1, H_2)$ onto $Hom_{M_1, M_2^o}(H_2, H_3)$, which sends the subset $Hom_{\chi}$ onto the subset $Hom_{\chi^2}$, and extends to an antilinear isomorphism $\gF$ from $\gH$ onto the Hilbert space $\gH_2$. This Hilbert space can be identified with the standard Hilbert space $L^2(M'_1\cap M_3)$, using the normal semi-finte faithful weight $\chi^2_3=\chi^2\circ T_{3|M'_1\cap M_3}$ on $M'_1\cap M_3$, and we can construct a normal faithful representation $\pi_2$ of $M'_1\cap M_4$ on $\gH_2$, which will verify ([EV] 3.3), for all $X$ in $M'_0\cap M_3$ :
\[\gF\pi(X)^*\gF^*=\pi_2(j_2(X))\]
Then, we can construct a pseudo-multiplicative unitary $W_2$ over the basis $(M'_1\cap M_2)^o$, with respect to the representations $r_2$ and $\hat{r}_2$ of $M'_1\cap M_2$ and the anti-representation $s_2$ of $M'_1\cap M_2$. Then, the antilinear isomorphism $\gF$ from $\gH$ to $\gH_2$ intertwines $r$ with $\hat{r}_2$, $\hat{r}$ with $r_2$, $s$ with $s_2$, and we have ([EV], 5.4) :
\[W_2=(\gF_s\underset{\chi}{\otimes}{}_r\gF)\sigma_{\chi^o}W^*\sigma_{\chi^o}(\gF^*_{r_2}\underset{\chi^{2 o}}{\otimes}{}_{s_2}\gF^*)\]
Writing $\hat{J}=\gF^*J_{\chi^2_3}\gF$, we construct an antilinear involutive isomorphism $\hat{J}$ of $\gH$ which, using \ref{propj} applied to $M_1\subset M_2$ and the properties of $\gF$,  will verify, for all $x$ in $M'_0\cap M_2$ :
\[\hat{J}\pi (x^*)\hat{J}=\pi\circ j_1(x)\]
and, therefore, the application $x\rightarrow \hat{J}x^*\hat{J}$ is an anti-automorphism $j$ of $\pi (M'_0\cap M_2)'$, such that $j\circ \hat{r}=s$. 
\\Writing $\widehat{\Delta}=\gF^*\Delta_{\chi^2_3}^{-1}\gF$, we construct a positive invertible self-adjoint operator on $\gH$, which commutes with $\hat{J}$, and a one-parameter group of unitaries $\widehat{\Delta}^{it}=\gF^*\Delta_{\chi^2_3}^{it}\gF$, which, using \ref{tau} applied to $M_1\cap M_2$ and the properties of $\gF$, implements on $\pi(M'_0\cap M_2)'$ a one-parameter automorphism group $\tau_t$ which verify, for all $a$ in $M'_0\cap M_1$ :
\[\tau_t(\hat{r}(a))=\hat{r}(\sigma_{t}^{\chi}(a))\]
\[\tau_t(s(a))=s(\sigma_{t}^{\chi}(a))\]
and we get, using \ref{thtau} applied to $M_1\subset M_2$ and the properties of $\gF$, that, for all $\xi$, $\eta$ in $D(_r\gH, \chi)\cap D(\gH_s, \chi^o)$, the operator $(\omega_{\xi, \eta}*id)(W)$ belongs to $D(\tau_{i/2})$ and that :
\[j\circ\tau_{-i/2}((\omega_{\xi, \eta}*id)(W))=(\omega_{\xi, \eta}*id)(W^*)\]
\end{proof}

\subsection{Theorem}
\label{thW}
{\it (i) It is possible to define a one-parameter unitary group $\Delta_{\chi_2}^{it}{}_r\underset{\chi^o}{\otimes}{}_s\widehat{\Delta}^{it}$ on $\gH_r\underset{\chi^o}{\otimes}{}_s\gH$, with natural values on elementary tensors, and a one-parameter unitary group $\Delta_{\chi_2}^{it}{}_s\underset{\chi}{\otimes}{}_{\hat{r}}\widehat{\Delta}^{it}$ on  $\gH_s\underset{\chi}{\otimes}{}_{\hat{r}}\gH$, such that :
\[W(\Delta_{\chi_2}^{it}{}_r\underset{\chi^o}{\otimes}{}_s\widehat{\Delta}^{it})=
(\Delta_{\chi_2}^{it}{}_s\underset{\chi}{\otimes}{}_{\hat{r}}\widehat{\Delta}^{it})W\]
(ii) It is possible to define an antilinear bijective isometry $J_{\chi_2}{}_r\underset{\chi^o}{\otimes}{}_s\hat{J}$ from $\gH_r\underset{\chi^o}{\otimes}{}_s\gH$ onto 
$\gH_s\underset{\chi}{\otimes}{}_{\hat{r}}\gH$, with natural values on elementary tensors, whose inverse is the antilinear bijective isometry $J_{\chi_2}{}_s\underset{\chi}{\otimes}{}_{\hat{r}}\hat{J}$, from $\gH_s\underset{\chi}{\otimes}{}_{\hat{r}}\gH$ onto $\gH_r\underset{\chi^o}{\otimes}{}_s\gH$, defined the same way with natural values on elementary tensors. This antilinear bijective isometry verify :}
\[(J_{\chi_2}{}_r\underset{\chi^o}{\otimes}{}_s\hat{J})W^*(J_{\chi_2}{}_r\underset{\chi^o}{\otimes}{}_s\hat{J})=W\]
\begin{proof}
The proof is identical to ([E2], 7.2). \end{proof}

\subsection{Theorem}
\label{thtau3}
{\it (i) The one-parameter automorphism group $\tau_t$ defined in \ref{thtau2} satisfies, for all $t$ in $\mathbb{R}$ :
\[\Gamma\circ\tau_t=(\tau_t{}_s\underset{\chi}{*}{}_{\hat{r}}\tau_t)\circ\Gamma\]
(ii) the anti-isomorphism $j$ defined in \ref{thtau2} satisfies $\Gamma\circ j=\varsigma_{\chi}(j{}_s\underset{\chi}{*}{}_{\hat{r}}j)\Gamma$.}
\begin{proof}
Let $\xi$ in $D(_r\gH, \chi)$ and $\eta$ in $D(\gH_s, \chi^o)$; using the pentagonal equation, we get that :
\begin{multline*}
\Gamma((\omega_{\xi, \eta}*id)(W))=\\(\omega_{\xi, \eta}*id*id)(W{}_s\underset{\chi}{\otimes}{}_{\hat{r}}1_\gH)(\sigma_\chi{}_s\underset{\chi}{\otimes}{}_{\hat{r}}1_\gH)(1_{\gH}{}_s\underset{\chi}{\otimes}{}_{\hat{r}}W)\sigma_{2\chi^o}(1_\gH{}_r\underset{\chi^o}{\otimes}{}_s\sigma_{\chi}))
\end{multline*}
On the other hand, using \ref{thW}(i), we get, for all $t$ in $\mathbb{R}$ :
\[\tau_t((\omega_{\xi, \eta}*id)(W))=\widehat{\Delta}^{it}(\omega_{\xi, \eta}*id)(W)\widehat{\Delta}^{-it}=
(\omega_{\Delta_{\chi_2}^{it}\xi, \Delta_{\chi_2}^{it}\eta}*id)(W)\]
from which we get that :
\begin{multline*}
\Gamma(\tau_t(\omega_{\xi, \eta}*id)(W)))=\\(\omega_{\Delta_{\chi_2}^{it}\xi, \Delta_{\chi_2}^{it}\eta}*id*id)(W{}_s\underset{\chi}{\otimes}{}_{\hat{r}}1_\gH)(\sigma_\chi{}_s\underset{\chi}{\otimes}{}_{\hat{r}}1_\gH)(1_{\gH}{}_s\underset{\chi}{\otimes}{}_{\hat{r}}W)\sigma_{2\chi^o}(1_\gH{}_r\underset{\chi^o}{\otimes}{}_s\sigma_{\chi}))
\end{multline*}
and, using again \ref{thW}(i), we have :
\[\Gamma(\tau_t(\omega_{\xi, \eta}*id)(W)))=
(\tau_t{}_s\underset{\chi}{*}{}_{\hat{r}}\tau_t)\Gamma((\omega_{\xi, \eta}*id)(W))\]
from which we get that $\Gamma\circ\tau_t(x)=(\tau_t{}_s\underset{\chi}{*}{}_{\hat{r}}\tau_t)\circ\Gamma(x)$, for all $x$ in the involutive algebra generated by the elements of the form $(\omega_{\xi, \eta}*id)(W)$, which, by continuity, gives (i). 
\\ On the other hand, using \ref{thW}(ii), we get :
\[(j((\omega_{\xi, \eta}*id)(W))=\hat{J}(\omega_{\xi, \eta}*id)(W)^*\hat{J}=(\omega_{J_{\chi_2}\eta, J_{\chi_2}\xi}*id)(W)\]
from which we get that :
\begin{multline*}
\Gamma(j((\omega_{\xi, \eta}*id)(W))=\\
(\omega_{J_{\chi_2}\eta, J_{\chi_2}\xi}*id*id)(W{}_s\underset{\chi}{\otimes}{}_{\hat{r}}1_\gH)(\sigma_\chi{}_s\underset{\chi}{\otimes}{}_{\hat{r}}1_\gH)(1_{\gH}{}_s\underset{\chi}{\otimes}{}_{\hat{r}}W)\sigma_{2\chi^o}(1_\gH{}_r\underset{\chi^o}{\otimes}{}_s\sigma_{\chi}))
\end{multline*}
and, using again \ref{thW}(ii), we finally have :
\[(j{}_s\underset{\chi}{*}{}_{\hat{r}}j)\Gamma(j((\omega_{\xi, \eta}*id)(W))=\varsigma_{\chi^o}\Gamma((\omega_{\xi, \eta}*id)(W))\]
from which we get that $\Gamma\circ j(x)=\varsigma_{\chi}(j{}_s\underset{\chi}{*}{}_{\hat{r}}j)\Gamma(x)$, for all $x$ in the involutive algebra generated by the elements of the form $(\omega_{\xi, \eta}*id)(W)$, which, by continuity, gives (ii).  \end{proof}

\subsection{Theorem}
\label{thinv}
{\it (i) Let us define $\tilde{T_2}'$ the normal semi-finite faithful operator valued weight from $\pi(M'_0\cap M_2)'=J_{\chi_2}\pi(M'_0\cap M_2)J_{\chi_2}$ to $s(M'_0\cap M_1)=J_{\chi_2}\pi(M'_0\cap M_1)J_{\chi_2}$ given, for all positive $x$ in $\pi(M'_0\cap M_2)'$ by :
\[\tilde{T_2}'(x)=J_{\chi_2}\tilde{T_2}(J_{\chi_2}xJ_{\chi_2})J_{\chi_2}\]
which verify $\chi'_2=\chi\circ s^{-1}\circ \tilde{T_2}'$. We have then, for all positive $x$ in  $\pi(M'_0\cap M_2)'$:
\[\tilde{T_2}'(x)=(\tilde{T_2}'_s\underset{\chi}{*}{}_{\hat{r}} id)\Gamma (x)=(\chi'_2{}_s\underset{\chi}{*}{}_{\hat{r}} id)\Gamma (x)\]
Therefore, $\tilde{T_2}'$ is a right-invariant operator valued weight for the Hopf-bimodule $\mathcal A(W)$.
\newline
(ii) let us consider $j\circ\tilde{T_2}'\circ j$ the normal semi-finite operator valued weight from $\pi(M'_0\cap M_2)'$ to $j\circ s(M'_0\cap M_1)=\hat{r}(M'_0\cap M_1)$; it verifies $\chi'_2\circ j=\chi\circ\hat{r}^{-1}\circ j\circ\tilde{T_2}'\circ j$, and, for all positive $x$ in  $\pi(M'_0\cap M_2)'$, we get :
\[j\circ\tilde{T_2}'\circ j(x)=(id{}_s\underset{\chi}{*}{}_{\hat{r}}j\circ\tilde{T_2}'\circ j)\Gamma (x)=
(id{}_s\underset{\chi}{*}{}_{\hat{r}}\chi'_2\circ j)\Gamma (x)\]
Therefore, $j\circ\tilde{T_2}'\circ j$ is a left-invariant operator valued weight for the Hopf-bimodule $\mathcal A(W)$.}
\begin{proof}
The proof of (i) is identical to ([E2], 7.5). Then (ii) is a straightforward corollary of \ref{thtau3}(ii). \end{proof}

\subsection{Proposition}
\label{propsigma}
{\it For any $x$ in $\pi(M'_0\cap M_2)'$, and $t$ in $\mathbb{R}$, we have :}
\[\Gamma\circ\sigma_t^{\chi'_2}(x)=(\sigma_t^{\chi'_2}{}_s\underset{\chi}{*}{}_{\hat{r}}\tau_{-t})\Gamma (x)\]
\[\Gamma\circ\sigma_t^{\chi'_2\circ j}(x)=(\tau_t{}_s\underset{\chi}{*}{}_{\hat{r}}\sigma_t^{\chi'_2\circ j})\Gamma (x)\]

\begin{proof}
The first formula is a direct corollary of \ref{thW}(i). 
\newline
As $\sigma_t^{\chi'_2\circ j}(x)=j\circ\sigma_{-t}^{\chi'_2}\circ j(x)$, the second formula comes from the first and \ref{thtau3}(ii). \end{proof}

\subsection{Theorem}
\label{thbimod}
{\it Let $M_0\subset M_1$ be a depth 2 inclusion of $\sigma$-finite von Neumann algebras, equipped with a regular normal semi-finite faithful operator-valued weight, in the sense of ([EV], \ref{exW}); let $r$, $\hat{r}$ (resp. $s$), be the representations (resp. the antirepresntation) of $M'_0\cap M_1$ on $L^2(M'_0\cap M_2)$ defined in  ([EV], \ref{exW}) and $W$ be the pseudo-multiplicative associated ([EV], \ref{exW}); let $((M'_0\cap M_1)^o, s, \hat{r}, \pi(M'_0\cap M_2)', \Gamma)$ be the Hopf-bimodule associated in ([EV], \ref{exW}). Let us suppose that there exists on $M'_0\cap M_1$ a normal faithful semi-finite weight $\chi$ invariant under the modular automorphism group $\sigma_t^{T_1}$. Then, the Hopf-bimodule bears the extra structures :
\newline
(i) there exists an anti-automorphism $j$ of $\pi(M'_0\cap M_2)'$ such that $j\circ s=\hat{r}$, and :
\[\Gamma\circ j=\varsigma_{\chi}(j{}_s\underset{\chi}{*}{}_{\hat{r}}j)\Gamma\]
(ii) there exists a one-parameter automorphism group $\tau_t$ of $\pi(M'_0\cap M_2)'$, such that, for all $t$ in $\mathbb{R}$, $j\circ\tau_t=\tau_t\circ j$, and :
\[\Gamma\circ\tau_t=(\tau_t{}_s\underset{\chi}{*}{}_{\hat{r}}\tau_t)\circ\Gamma\]
Moreover, if $\xi$,$\eta$ are in $D(_r\gH, \chi)\cap D(\gH_s, \chi^o)$, the operator $(\omega_{\xi, \eta}*id)(W)$ belongs to $D(\tau_{-i/2})$, and we have : 
\[j\circ\tau_{-i/2}((\omega_{\xi, \eta}*id)(W))=(\omega_{\xi, \eta}*id)(W^*)\]
Moreover, for all $a$ in $M'_0\cap M_1$, we have :
\[\tau_t(s(a))=s(\sigma_{t}^{\chi}(a))\]
\[\tau_t(\hat{r}(a))=\hat{r}(\sigma_{t}^{\chi}(a))\]
(iii) there exists a normal semi-finite faithful operator-valued weight $\tilde{T'_2}$ from $\pi(M'_0\cap M_2)'$ so $s(M'_0\cap M_1)$, which is right-invariant in the sense that, for all positive $x$ in  $\pi(M'_0\cap M_2)'$:
\[\tilde{T_2}'(x)=(\tilde{T_2}'_s\underset{\chi}{*}{}_{\hat{r}} id)\Gamma (x)=(\chi'_2{}_s\underset{\chi}{*}{}_{\hat{r}} id)\Gamma (x)\]
where we have $\chi'_2=\chi\circ s^{-1}\circ \tilde{T_2}'$. Moreover, for all $t$ in $\mathbb{R}$, we have :}
\[\Gamma\circ\sigma_t^{\chi'_2}=(\sigma_t^{\chi'_2}{}_s\underset{\chi}{*}{}_{\hat{r}}\tau_{-t})\Gamma\]

\begin{proof} 
This is a r\'{e}sum\'{e} of \ref{thtau2}, \ref{thtau3}, \ref{thinv} and \ref{propsigma}. \end{proof}

\section{Modular theory on $\mathcal A(W)$.}
\label{modular2}
In this chapter, mostly inspired by the work of Lesieur ([L], chap. 3, 4 and 5), we obtain results about the modular theory of the left-invariant weight and of the right-invariant weight constructed on $\mathcal A(W)$ in \ref{thinv}.  We construct a scaling operator and a modulus (\ref{modulus}), and prove that the modulus is, in a special sense, a cocharacter (\ref{thmodulus}). 

\subsection{Proposition}
\label{propsigma2}
{\it (i) For any $a$ in $M'_0\cap M_1$, we have :
\[\sigma_t^{\chi'_2}(\hat{r}(a))=\hat{r}(\sigma_{-t}^{T_1}(a))\]
\[\sigma_t^{\chi'_2\circ j}(s(a))=s(\sigma_t^{T_1}(a))\]
(ii) For any positive $x$ in $\pi(M'_0\cap M_2)'$, and $t$ in $\mathbb{R}$, we have :}
\[\chi'_2\circ j(x)=\chi'_2\circ j\circ\sigma_t^{\chi'_2}\circ\tau_t(x)\]
\[\chi'_2(x)=\chi'_2\circ\sigma_{-t}^{\chi'_2\circ j}\circ\tau_t(x)\]
\begin{proof}
Using \ref{propj}(ii), we get  :
\begin{multline*}
\sigma_t^{\chi'_2}(\hat{r}(a))=\Delta_{\chi_2}^{-it}\hat{r}(a)\Delta_{\chi_2}^{it}=h^{-it}\hat{r}(a)h^{it}\\
=\sigma_{-t}^{\psi_3}(\hat{r}(a))=\hat{r}(\sigma_{-t}^{\psi_1}(a))=\hat{r}(\sigma_{-t}^{T_1}(a))
\end{multline*}
which is the first result of (i). We then get :
\begin{multline*}
\sigma_t^{\chi'_2\circ j}(s(a))=j\circ\sigma_{-t}^{\chi'_2}\circ j(s(a))=j\circ\sigma_{-t}^{\chi'_2}(\hat{r}(a))\\=
j\circ\hat{r}(\sigma_t^{T_1}(a))=s(\sigma_t^{T_1}(a))
\end{multline*}
which finishes the proof of (i). 
\newline
From \ref{propsigma} and \ref{thtau3}(i), we get :
\[\Gamma\circ\sigma_t^{\chi'_2}\circ\tau_t(x)=(\sigma_t^{\chi'_2}\circ\tau_t{}_s\underset{\chi}{*}{}_{\hat{r}}id)\Gamma (x)\]
and then, using \ref{thinv} (ii), we have :
\begin{eqnarray*}
j\circ\tilde{T_2}'\circ j\circ\sigma_t^{\chi'_2}\circ\tau_t(x)
&=&(id{}_s\underset{\chi}{*}{}_{\hat{r}}\chi'_2\circ j)\Gamma \circ\sigma_t^{\chi'_2}\circ\tau_t(x)\\
&=&\sigma_t^{\chi'_2}\circ\tau_t(j\circ\tilde{T_2}'\circ j(x))
\end{eqnarray*}
If $a$ is in $M'_0\cap M_1$, we have, by \ref{thtau2}, $\tau_t(\hat{r}(a))=\hat{r}(\sigma_t^{\chi}(a))$, and, therefore, using (i), $\sigma_t^{\chi'_2}\circ\tau_t(\hat{r}(a))=\hat{r}(\sigma_{-t}^{T_1}\sigma_t^{\chi}(a))$. So, we get :
\begin{eqnarray*}
\chi'_2\circ j\circ\sigma_t^{\chi'_2}\circ\tau_t(x)
&=&\chi\circ\hat{r}^{-1}\circ j\circ\tilde{T_2}'\circ j\circ\sigma_t^{\chi'_2}\circ\tau_t(x)\\
&=&\chi\circ\sigma_{-t}^{T_1}\circ\sigma_t^{\chi}\circ\hat{r}^{-1}\circ j\circ\tilde{T_2}'\circ j(x)\\
&=&\chi\circ\hat{r}^{-1}\circ j\circ\tilde{T_2}'\circ j(x)\\
&=&\chi'_2\circ j(x)
\end{eqnarray*}
which gives the first formula of (ii). The second formula comes then from the commutation of $j$ and $\tau_t$ (\ref{thtau2}). \end{proof}

\subsection{Theorem}
\label{thsigmatau}
{\it The three automorphism groups $\tau_t$, $\sigma_t^{\chi'_2}$, $\sigma_t^{\chi'_2\circ j}$ are two by two commuting.}

\begin{proof}
Using \ref{propsigma2}, we get that the automorphism group $\sigma_s^{\chi'_2\circ j}$ is commuting, for all $t$ in $\mathbb{R}$, with the automorphism $\sigma_t^{\chi'_2}\circ\tau_t$. Which we can write, for all $s$, $t$ in $\mathbb{R}$ :
\[\sigma_s^{\chi'_2\circ j}\circ\sigma_t^{\chi'_2}\circ\tau_t\circ\sigma_{-s}^{\chi'_2\circ j}=\sigma_t^{\chi'_2}\circ\tau_t\]
Therefore, using \ref{propsigma}, we get :
\begin{eqnarray*}
\Gamma\circ\tau_s\circ\sigma_t^{\chi'_2}\circ\tau_t\circ\tau_{-s}
&=&
(\tau_s{}_s\underset{\chi}{*}{}_{\hat{r}}\tau_s)\Gamma\circ\sigma_t^{\chi'_2}\circ\tau_t\circ\tau_{-s}\\
&=&(\tau_s{}_s\underset{\chi}{*}{}_{\hat{r}}\tau_s)(\sigma_t^{\chi'_2}\circ\tau_t{}_s\underset{\chi}{*}{}_{\hat{r}}id)\Gamma\circ\tau_{-s}\\
&=&(\tau_s{}_s\underset{\chi}{*}{}_{\hat{r}}\tau_s)(\sigma_t^{\chi'_2}\circ\tau_t{}_s\underset{\chi}{*}{}_{\hat{r}}id)(\tau_{-s}{}_s\underset{\chi}{*}{}_{\hat{r}}\tau_{-s})\Gamma\\
&=&(\tau_s{}_s\underset{\chi}{*}{}_{\hat{r}}\sigma_s^{\chi'_2\circ j})(\sigma_t^{\chi'_2}\circ\tau_t{}_s\underset{\chi}{*}{}_{\hat{r}}id)(\tau_{-s}{}_s\underset{\chi}{*}{}_{\hat{r}}\sigma_{-s}^{\chi'_2\circ j})\Gamma\\
&=&(\tau_s{}_s\underset{\chi}{*}{}_{\hat{r}}\sigma_s^{\chi'_2\circ j})(\sigma_t^{\chi'_2}\circ\tau_t{}_s\underset{\chi}{*}{}_{\hat{r}}id)\Gamma\circ\sigma_{-s}^{\chi'_2\circ j}\\
&=&(\tau_s{}_s\underset{\chi}{*}{}_{\hat{r}}\sigma_t^{\chi'_2\circ j})\Gamma\circ\sigma_t^{\chi'_2}\circ\tau_t\circ\sigma_{-s}^{\chi'_2\circ j}\\
&=&\Gamma\circ\sigma_{s}^{\chi'_2\circ j}\circ\sigma_t^{\chi'_2}\circ\tau_t\circ\sigma_{-s}^{\chi'_2\circ j}\\
&=&\Gamma\circ\sigma_t^{\chi'_2}\circ\tau_t
\end{eqnarray*}
and, by the injectivity of $\Gamma$, we get that 
\[\sigma_t^{\chi'_2}\circ\tau_t=\tau_s\circ\sigma_t^{\chi'_2}\circ\tau_t\circ\tau_{-s}
=\tau_s\circ\sigma_t^{\chi'_2}\circ\tau_{-s}\circ\tau_{t}\]
from which we get that $\sigma_t^{\chi'_2}=\tau_s\circ\sigma_t^{\chi'_2}\circ\tau_{-s}$, which gives the commutation of the one-parameter groups of automorphisms $\tau_s$ and $\sigma_t^{\chi'_2}$. 
\newline
It is then straightforward to get the commutation of $\tau_s$ with $\sigma_{t}^{\chi'_2\circ j}$. 
\newline
Finally, we have :
\begin{eqnarray*}
\Gamma\circ\sigma_{t}^{\chi'_2}\circ\sigma_s^{\chi'_2\circ j}
&=&(\sigma_t^{\chi'_2}{}_s\underset{\chi}{*}{}_{\hat{r}}\tau_{-t})\Gamma \circ\sigma_s^{\chi'_2\circ j}\\
&=&(\sigma_t^{\chi'_2}{}_s\underset{\chi}{*}{}_{\hat{r}}\tau_{-t})
(\tau_s{}_s\underset{\chi}{*}{}_{\hat{r}}\sigma_s^{\chi'_2\circ j})\Gamma\\
&=&(\tau_s{}_s\underset{\chi}{*}{}_{\hat{r}}\sigma_s^{\chi'_2\circ j})
(\sigma_t^{\chi'_2}{}_s\underset{\chi}{*}{}_{\hat{r}}\tau_{-t})\Gamma\\
&=&(\tau_s{}_s\underset{\chi}{*}{}_{\hat{r}}\sigma_s^{\chi'_2})
\Gamma\circ\sigma_{t}^{\chi'_2\circ j}\\
&=&\Gamma\circ\sigma_s^{\chi'_2\circ j}\circ\sigma_{t}^{\chi'_2}
\end{eqnarray*}
which, by the injectivity of $\Gamma$, gives the commutation of the one-parameter groups $\sigma_s^{\chi'_2}$ and $\sigma_{t}^{\chi'_2\circ j}$. \end{proof}

\subsection{Theorem}
\label{modulus}
{\it There exists a positive invertible operator $\delta'$ affiliated to $\pi(M'_0\cap M_2)'$, called the modulus, and a positive invertible element $\lambda$ affiliated to $Z(M_0)\cap Z(M_1)$, called the scaling operator, such that, for all $s$, $t$ in $\mathbb{R}$, $x$ positive in $\pi(M'_0\cap M_2)'$ :
\newline
(i) $\sigma_s^{\chi'_2}(\delta'^{it})=\pi(\lambda)^{ist}\delta'^{it}$
\newline
(ii) $(D\chi'_2\circ j:D\chi'_2)_t=\pi(\lambda)^{it^2/2}\delta'^{it}$
\newline
(iii) $\chi'_2\circ\sigma_t^{\chi'_2\circ j}(x)=\chi'_2\circ\tau_t(x)=\chi'_2(\pi(\lambda)^{t/2}x\pi(\lambda)^{t/2})$}

\begin{proof}
Using ([V3], Prop. 5.2), and the fact that the automorphism groups $\sigma_s^{\chi'_2}$ and $\sigma_{t}^{\chi'_2\circ j}$ commute by \ref{thsigmatau}, we get the existence of $\delta'$ and a positive invertible element $\mu$ in $\pi(Z(M'_0\cap M_2))$ such that 
\[(D\chi'_2\circ j:D\chi'_2)_t=\mu^{it^2/2}\delta'^{it}\]
\[\sigma_s^{\chi'_2}(\delta'^{it})=\mu^{ist}\delta'^{it}\]
By the unicity of this decomposition, we easily get that $j(\delta')=\delta'^{-1}$, and $j(\mu)=\mu$. Moreover, a calculation completely similar to ([V3], Prop. 5.5) leads to :
\[\chi'_2\circ\sigma_t^{\chi'_2\circ j}(x)=\chi'_2(\mu^{t/2}x\mu^{t/2})\]
and, by \ref{propsigma2}, we have :
\[\chi'_2\circ\sigma_t^{\chi'_2\circ j}(x)=\chi'_2\circ\tau_t(x)\]
and, therefore, we get that $\mu^{ist}=(D\chi'_2\circ\tau_t:D\chi'_2)_s$. 
\newline
We have, using \ref{tau} and \ref{thinv} :
\[\chi'_2\circ\tau_t=
\chi\circ s^{-1}\circ\tilde{T_2}'\circ\tau_t=\chi\circ s^{-1}\circ\tau_{-t}\circ\tilde{T_2}'\circ\tau_t\]
from which we get that $\mu^{ist}=(D\tau_{-t}\circ\tilde{T_2}'\circ\tau_t:D\tilde{T_2}')_s$ and belongs therefore to $s(M'_0\cap M_3)$; as $\mu$ is affiliated to $\pi(Z(M'_0\cap M_2))$, we have $J_{\chi_2}\mu J_{\chi_2}=\mu^{-1}$, and $\mu$ is therefore affiliated to $r(M'_0\cap M_1)$; as $j(\mu)=\mu$, we get that $\mu$ is also affiliated to $\hat{r}(M'_0\cap M_1)=M'_2\cap M_3$. Let us write $\mu=\pi(\lambda)$; we get that $\lambda^{it}$ belongs to $M'_0\cap M_1\cap Z(M'_0\cap M_2)\cap (M'_0\cap M_3)'\cap M'_2\cap M_3$; so $\lambda^{it}$ belongs to $Z(M_1)$; we have $J_1\lambda J_1=\lambda^{-1}$, and, therefore, as $\lambda^{it}$ belongs to $M'_2$, it belongs also to $M_0$, and, more precisely, to $Z(M_0)$, which finishes the proof. \end{proof}

\subsection{Corollary}
\label{cormodulus}
{\it  The  operator-valued weight $\tilde{T'_2}$ (resp. $j\circ\tilde{T'_2}\circ j$) from $\pi(M'_0\cap M_2)'$ to $s(M'_0\cap M_1)$ (resp. $\hat{r}(M'_0\cap M_1)$) introduced in \ref{thinv}, satisfies, for all $t$ in $\mathbb{R}$ and $x$ positive in  $\pi(M'_0\cap M_2)'$ :}
\[\tilde{T'_2}\circ\sigma_t^{\chi'_2\circ j}(x)=\pi(\lambda)^{t/2}s\circ\sigma_t^{T_1}\circ s^{-1}(\tilde{T'_2}(x))\pi(\lambda)^{t/2}\]
\[ j\circ\tilde{T'_2}\circ j\circ\sigma_t^{\chi'_2}(x)=\pi(\lambda)^{-t/2}\hat{r}\circ\sigma_{-t}^{T_1}\circ\hat{r}^{-1}(j\circ\tilde{T'_2}\circ j(x))\pi(\lambda)^{-t/2}\]

\begin{proof}
Using \ref{thinv} and \ref{propsigma}, we get :
\[\tilde{T'_2}(\sigma^{\chi'_2\circ j}_t(x))=(\chi'_2{}_s\underset{\chi}{*}_{\hat{r}}id)\Gamma(\sigma^{\chi'_2\circ j}_t(x))=
(\chi'_2{}_s\underset{\chi}{*}_{\hat{r}}id)(\tau_t{}_s\underset{\chi}{*}_{\hat{r}}\sigma^{\chi'_2\circ j}_t)\Gamma(x)\]
Using then \ref{modulus}, we get it is equal to :
\[\sigma^{\chi'_2\circ j}_t((\chi'_2{}_s\underset{\chi}{*}{}_{\hat{r}}id)(\pi(\lambda)^{t/2}_s\underset{\chi}{\otimes}{}_{\hat{r}}1)\Gamma(x)(\pi(\lambda)^{t/2}_s\underset{\chi}{\otimes}{}_{\hat{r}}1))\]
and, using the fact that $\pi(\lambda)^{t/2}$ is affiliated 
first, to $\hat{r}(M'_0\cap M_1)$, second, to $s(M'_0\cap M_1)$, and , third, to the center of $\pi(M'_0\cap M_2)'$, we obtain it is equal to :
\begin{eqnarray*}
\sigma^{\chi'_2\circ j}_t((\chi'_2{}_s\underset{\chi}{*}_{\hat{r}}id)\Gamma(\pi(\lambda)^{t/2}x\pi(\lambda)^{t/2}))
&=&
\sigma^{\chi'_2\circ j}_t(\tilde{T'_2}(\pi(\lambda)^{t/2}x\pi(\lambda)^{t/2}))\\
&=&
\sigma^{\chi'_2\circ j}_t(\pi(\lambda)^{t/2}\tilde{T'_2}(x)\pi(\lambda)^{t/2})\\
&=&
\pi(\lambda)^{t/2}\sigma^{\chi'_2\circ j}_t(\tilde{T'_2}(x))\pi(\lambda)^{t/2}\\
\end{eqnarray*}
which, using \ref{propsigma2}, gives the first result. The second formula is then obtained by similar calculations. \end{proof}

\subsection{Corollary}
\label{cormodulus2}
{\it Let us denote by $T_{\chi'_2, \tilde{T'_2}}^{\chi'_2\circ j}$ the subset of $T_{\chi'_2, \tilde{T'_2}}$ (defined in \ref{prop}) made of elements $a$ in $\gN_{\tilde{T'_2}}\cap\gN_{\chi'_2}\cap\gN_{\chi'_2\circ j}$, analytic with respect to both $\chi'_2$ and $\chi'_2\circ j$, such that, for all $z$, $z'$ in $\mathbb{C}$, $\sigma_z^{\chi'_2}\circ\sigma_{z'}^{\chi'_2\circ j}(a)$ belongs to $\gN_{\tilde{T'_2}}\cap\gN_{\chi'_2}\cap\gN_{\chi'_2\circ j}$. This linear space is weakly dense in $\pi(M'_0\cap M_2)'$, and the set of $\Lambda_{\chi'_2}(a)$ (resp. $\Lambda_{\chi'_2\circ j}(a)$) , for $a$ in $T_{\chi'_2, \tilde{T'_2}}^{\chi'_2\circ j}$, is a linear dense subset in $H_{\chi'_2}$ (resp. $H_{\chi'_2\circ j}$); moreover, the subset $J_{\chi_2}\Lambda_{\chi'_2}(T_{\chi'_2, \tilde{T'_2}}^{\chi'_2\circ j})$ is included in the domain of $\delta'^z$, for all $z$ in $\mathbb{C}$, and is an essential domain for $\delta'^z$. }

\begin{proof}
Let us take $x$ in $T_{\chi'_2, \tilde{T'_2}}^+$; let us write $\lambda=\int_0^\infty tde_t$ and let us define $f_p=\int_{1/p}^p de_t$; if we put :
\[x_{q,p}=\pi(f_p)\sqrt{\frac{q}{\pi}}\int_{-\infty}^{+\infty}e^{-qt^2}\sigma_t^{\chi'_2\circ j}(x)dt\]
we, by classical remarks (see [EN], 10.12 for instance for similar calculations) obtain that $x_{q,p}$ is in 
$T_{\chi'_2, \tilde{T'_2}}^+$ and, moreover, is analytical with respect to $\chi'_2\circ j$, and that all $\sigma_z^{\chi'_2\circ j}(x_{q,p})$ belong to $T_{\chi'_2, \tilde{T'_2}}$. Moreover, when $q$, $p$ go to infinity, $x_{q,p}$ is weakly converging to $x$, and $\Lambda_{\chi'_2}(x_{q,p})$ is converging to $\Lambda_{\chi'_2}(x)$. We can verify also that $\Lambda_{\tilde{T'_2}}(x_{q,p})$ is weakly converging to $\Lambda_{\tilde{T'_2}}(x)$. 
\newline
As, for any $t$ in $\mathbb{R}$, and $y$ in $\pi(M'_0\cap M_2)'$, we have :
\[\delta'^{it}y\delta'^{-it}=\sigma_t^{\chi'_2\circ j}\circ\sigma_{-t}^{\chi'_2}(y)\]
we see that, for all such elements $x_{q,p}$, and $z$ in $\mathbb{C}$, $\delta'^{iz}x_{q,p}\delta'^{-iz}$ is bounded and belongs to $T_{\chi'_2, \tilde{T'_2}}$; in particular $\delta'^{-1/2}x_{q, p}\delta'^{1/2}$ belongs to $\gM_{\chi'_2}\cap\gM_{\tilde{T'_2}}$, and is analytic with respect to both $\chi'_2$ and $\chi'_2\circ j$.  Using then the operators $e_n$ introduced in ([V3], 1.1), which are analytic to both $\chi'_2$ and $\chi'_2\circ j$ and converging to $1$ when $n$ goes to infinity, we get that $e_nx_{q,p}$ belongs to $\gN_{\chi'_2}\cap\gN_{\tilde{T'_2}}$; on the other hand, as
\[e_nx_{q,p}\delta'^{1/2}=(e_n\delta'^{1/2})\delta'^{-1/2}x_{q, p}\delta'^{1/2}\]
belongs to $\gN_{\chi'_2}$, we see, by ([V3], 3.3), that $e_nx_{q,p}$ belongs to $\gN_{\chi'_2\circ j}$, and, therefore, to $T_{\chi'_2, \tilde{T'_2}}^{\chi'_2\circ j}$, from which we then get all the results claimed. 
\end{proof}

\subsection{Proposition}
\label{proptau}
{\it As usual, we shall identify $H_{\chi'_2}$ with $H_{\chi_2}$ by writing, for all $x$ in $\gN_{\chi_2}$ :
\[\Lambda_{\chi'_2}(J_{\chi_2}xJ_{\chi_2})=J_{\chi_2}\Lambda_{\chi_2}(x)\]
(i) If $a'$ belongs to $\gN_{\tilde{T'_2}}\cap\gN_{\chi'_2}$, $J_{\chi_2}\Lambda_{\chi'_2}(a')$ belongs to $D(\gH_s, \chi^o)$, and $R^{s, \chi^o}(J_{\chi_2}\Lambda_{\chi'_2}(a'))=\Lambda_{\tilde{T_2}}(J_{\chi_2}a'J_{\chi_2})$; and $\Lambda_{\chi'_2}(a')$ belongs to $D({}_r\gH, \chi)$ and $R^{r, \chi}(\Lambda_{\chi'_2}(a'))=J_{\chi_2}\Lambda_{\tilde{T_2}}(J_{\chi_2}a'J_{\chi_2})J_\chi$. 
\newline
(ii) We have then, for all $y$ in $\gN_{\chi'_2}$, $\xi$ in $D(\gH_s)\cap D({}_{\hat{r}}\gH)$, $\eta$ in $D({}_{\hat{r}}\gH)$ :
\[\Lambda_{\chi'_2}((id{}_s\underset{\chi}{*}{}_{\hat{r}}\omega_{\xi, \eta})(\Gamma (y))=(id*\omega_{\xi, \eta})(W)\Lambda_{\chi'_2}(y)\]
(iii) If $y$, $z$, $z'$ belong to $\gN_{\tilde{T'_2}}\cap\gN_{\chi'_2}$, we have :
\[(\omega_{J_{\chi_2}\Lambda_{\chi'_2}(z), J_{\chi_2}\Lambda_{\chi'_2}(z')}{}_s\underset{\chi}{*}{}_{\hat{r}}id)\Gamma (y)=(\omega_{\Lambda_{\chi'_2}(y), J_{\chi_2}z^*\Lambda_{\chi'_2}(z')}*id)(W)\]
(iv) for $a$, $b$ in $\gN_{\tilde{T'_2}}\cap\gN_{\chi'_2}$, we have :
\[(\omega_{J_{\chi_2}\Lambda_{\chi'_2}(a), J_{\chi_2}\Lambda_{\chi'_2}(a)}{}_s\underset{\chi}{*}{}_{\hat{r}}id)\Gamma (b^*b)=(\omega_{\Lambda_{\chi'_2}(b^*b), J_{\chi_2}\Lambda_{\chi'_2}(a^*a)}*id)(W)\]
(v) for $a$, $b$ in $\gN_{\tilde{T'_2}}\cap\gN_{\chi'_2}$, we have :}
\[j((\omega_{J_{\chi_2}\Lambda_{\chi'_2}(a), J_{\chi_2}\Lambda_{\chi'_2}(a)}{}_s\underset{\chi}{*}{}_{\hat{r}}id)\Gamma (b^*b))=
(\omega_{J_{\chi_2}\Lambda_{\chi'_2}(b), J_{\chi_2}\Lambda_{\chi'_2}(b)}{}_s\underset{\chi}{*}{}_{\hat{r}}id)\Gamma (a^*a)\]

\begin{proof}
Using \ref{prop}, we get (i) by direct calculations. 
\newline
Using \ref{propbeta} and \ref{thW}, we get :
\begin{eqnarray*}
\Lambda_{\chi'_2}((id*\omega_{\xi, \eta})(\Gamma (y))
&=&J_{\chi_2}\Lambda_{\chi_2}(J_{\chi_2}(id*\omega_{\xi, \eta})(W(y{}_r\underset{\chi^o}{\otimes}{}_s 1)W^*)J_{\chi_2})\\
&=&J_{\chi_2}\Lambda_{\chi_2}((id*\omega_{\hat{J}\xi, \hat{J}\eta})(W^*(J_{\chi_2}yJ_{\chi_2}{}_s\underset{\chi}{\otimes}{}_r 1)W))\\
&=&J_{\chi_2}(id*\omega_{\hat{J}\xi, \hat{J}\eta})(W^*)\Lambda_{\chi_2}(J_{\chi_2}yJ_{\chi_2})\\
&=&(id*\omega_{\xi, \eta})(W)\Lambda_{\chi'_2}(y)
\end{eqnarray*}
which gives (ii); moreover, we have :
\begin{multline*}
(\Gamma (y)(J_{\chi_2}\Lambda_{\chi'_2}(z){}_s\underset{\chi}{\otimes}{}_{\hat{r}}\xi )|J_{\chi_2}\Lambda_{\chi'_2}(z'){}_s\underset{\chi}{\otimes}{}_{\hat{r}}\eta )\\
=((id*\omega_{\xi, \eta})\Gamma (y)J_{\chi_2}\Lambda_{\chi'_2}(z)|J_{\chi_2}\Lambda_{\chi'_2}(z'))\\
=(J_{\chi_2}z^*J_{\chi_2}\Lambda_{\chi'_2}((id*\omega_{\xi, \eta})(\Gamma (y))|J_{\chi_2}\Lambda_{\chi'_2}(z'))
\end{multline*}
which, using (ii), is equal to :
\begin{multline*}
(\Lambda_{\chi'_2}((id*\omega_{\xi, \eta})(\Gamma (y))|J_{\chi_2}z\Lambda_{\chi'_2}(z'))\\
=((id*\omega_{\xi, \eta})(W)\Lambda_{\chi'_2}(y)|J_{\chi_2}z\Lambda_{\chi'_2}(z'))
\end{multline*}
from which we get (iii) and (iv). 
\newline
Using (iv) and \ref{thW}, we have :
\begin{eqnarray*}
j((\omega_{J_{\chi_2}\Lambda_{\chi'_2}(a), J_{\chi_2}\Lambda_{\chi'_2}(a)}*id)\Gamma (b^*b))
&=&j((\omega_{\Lambda_{\chi'_2}(b^*b), J_{\chi_2}\Lambda_{\chi'_2}(a^*a)}*id)(W))\\
&=&\hat{J}(\omega_{\Lambda_{\chi'_2}(b^*b), J_{\chi_2}\Lambda_{\chi'_2}(a^*a)}*id)(W)^*\hat{J}\\
&=&\hat{J}(\omega_{J_{\chi_2}\Lambda_{\chi'_2}(a^*a), \Lambda_{\chi'_2}(b^*b)}*id)(W^*)\hat{J}\\
&=&(\omega_{\Lambda_{\chi'_2}(a^*a), J_{\chi_2}\Lambda_{\chi'_2}(b^*b)}*id)(W)\\
&=&(\omega_{J_{\chi_2}\Lambda_{\chi'_2}(b), J_{\chi_2}\Lambda_{\chi'_2}(b)}*id)\Gamma (a^*a)
\end{eqnarray*}
which is (v).  \end{proof}

\subsection{Proposition}
\label{proptau2}
{\it For all $t$ in $\mathbb{R}$ and $x$ in $\pi(M'_ˆ\cap M_2)'$, we have :}
\[\Gamma\tau_t(x)=(\sigma_t^{\chi'_2\circ j}{}_s\underset{\chi}{*}{}_{\hat{r}}\sigma_{-t}^{\chi'_2})\Gamma (x)\]
\[\Gamma\sigma_t^{\chi'_2\circ j}\sigma_{-t}^{\chi'_2}(x)=
(\sigma_t^{\chi'_2\circ j}\sigma_{-t}^{\chi'_2}{}_s\underset{\chi}{*}{}_{\hat{r}}\sigma_t^{\chi'_2\circ j}\sigma_{-t}^{\chi'_2})\Gamma (x)\]

\begin{proof}
Let $x$ be in $\pi(M'_0\cap M_2)'$ and $b$ in $T_{\chi'_2, \tilde{T_2}'}$; using \ref{modulus}, we have :
\begin{eqnarray*}
\omega_{J_{\chi_2}\Lambda_{\chi'_2}(b)}\circ\sigma_t^{\chi'_2\circ j}(x)
&=&
\chi'_2(\sigma_{-i/2}^{\chi'_2}(b^*)^*\sigma_t^{\chi'_2\circ j}(x)\sigma_{-i/2}^{\chi'_2}(b^*))\\
&=&\chi'_2\circ\sigma_t^{\chi'_2\circ j}((\sigma_{-i/2}^{\chi'_2}\sigma_{-t}^{\chi'_2\circ j}(b^*))^*x(\sigma_{-i/2}^{\chi'_2}\sigma_{-t}^{\chi'_2\circ j}(b^*)))\\
&=&\chi'_2((\sigma_{-i/2}^{\chi'_2}\sigma_{-t}^{\chi'_2\circ j}(b^*))^*\pi(\lambda)^{t/2}x\pi(\lambda)^{t/2}\sigma_{-i/2}^{\chi'_2}\sigma_{-t}^{\chi'_2\circ j}(b^*))\\
&=&\omega_{J_{\chi'_2}\Lambda_{\chi'_2}(\pi(\lambda)^{t/2}\sigma_{-t}^{\chi'_2\circ j}(b))}(x)
\end{eqnarray*}
The same way, we shall get $\omega_{J_{\chi_2}\Lambda_{\chi'_2}(b)}\circ\tau_t=
\omega_{J_{\chi'_2}\Lambda_{\chi'_2}(\pi(\lambda)^{t/2}\tau_{-t}(b))}$.
Then :
\[(\omega_{J_{\chi_2}\Lambda_{\chi'_2}(b)}*id)(\sigma_t^{\chi'_2\circ j}{}_s\underset{\chi}{*}{}_{\hat{r}}\sigma_{-t}^{\chi'_2})\Gamma(\tau_{-t}(a^*a))\]
is equal to :
\begin{multline*}
\sigma_{-t}^{\chi'_2}(\omega_{J_{\chi'_2}\Lambda_{\chi'_2}(\pi(\lambda)^{t/2}\sigma_{-t}^{\chi'_2\circ j}(b))}*id)\Gamma(\tau_{-t}(a^*a))\\
=j\sigma_t^{\chi'_2\circ j}j((\omega_{J_{\chi'_2}\Lambda_{\chi'_2}(\pi(\lambda)^{t/2}\sigma_{-t}^{\chi'_2\circ j}(b))}*id)\Gamma(\tau_{-t}(a^*a))
\end{multline*}
Using \ref{proptau} and the property of $\pi(\lambda)$, it is equal to :
\begin{multline*}
j\sigma_t^{\chi'_2\circ j}(\omega_{J_{\chi_2}\Lambda_{\chi'_2}(\tau_{-t}(a))}*id)\Gamma(\pi(\lambda)^{t/2}\sigma_{-t}^{\chi'_2\circ j}(b^*b)\pi(\lambda)^{t/2})\\
=j\sigma_t^{\chi'_2\circ j}(\omega_{J_{\chi_2}\Lambda_{\chi'_2}(\pi(\lambda)^{t/2}\tau_{-t}(a))}*id)\Gamma(\sigma_{-t}^{\chi'_2\circ j}(b^*b))\\
=j\sigma_t^{\chi'_2\circ j}((\omega_{J_{\chi_2}\Lambda_{\chi'_2}(a)}\circ\tau_t)*id)\Gamma(\sigma_{-t}^{\chi'_2\circ j}(b^*b))\
\end{multline*}
Using \ref{propsigma} and \ref{proptau}(v) again, we see it is equal to :
\[j(\omega_{J_{\chi_2}\Lambda_{\chi'_2}(a)}*id)\Gamma(b^*b)
=(\omega_{J_{\chi_2}\Lambda_{\chi'_2}(b)}*id)\Gamma(a^*a)\]
from which we get the first equality. Then, we get, using \ref{propsigma}, \ref{thtau3} and the first relation, that : 
\begin{multline*}
\Gamma\sigma_t^{\chi'_2\circ j}\sigma_{-t}^{\chi'_2}=(\tau_t*\sigma_t^{\chi'_2\circ j})(\sigma_{-t}^{\chi'_2}*\tau_t)\Gamma\\=(\sigma_{-t}^{\chi'_2}*\sigma_t^{\chi'_2\circ j})\Gamma\tau_t=
(\sigma_t^{\chi'_2\circ j}\sigma_{-t}^{\chi'_2}{}_s\underset{\chi}{*}{}_{\hat{r}}\sigma_t^{\chi'_2\circ j}\sigma_{-t}^{\chi'_2})\Gamma
\end{multline*}

\end{proof}

\subsection{Proposition}
\label{propdelta}
{\it It is possible to define a one-parameter unitary group of unitaries, denoted $\delta'^{it}{}_s\underset{\chi}{\otimes}{}_{\hat{r}}\delta'^{it}$, with natural values on elementary tensors; moreover, with have, for all $t$ in $\mathbb{R}$, and $x$ in $\pi(M'_0\cap M_2)'$ :
\[(\delta'^{it}{}_s\underset{\chi}{\otimes}{}_{\hat{r}}\delta'^{it})\Gamma (x)(\delta'^{-it}{}_s\underset{\chi}{\otimes}{}_{\hat{r}}\delta'^{it})=\Gamma(\delta'^{it}x\delta'^{-it})\]
and, the two one-parameter automorphism groups of unitaries $\Gamma(\delta'^{is})$ and $\delta'^{it}{}_s\underset{\chi}{\otimes}{}_{\hat{r}}\delta'^{it}$ commute.}

\begin{proof}
For all $t$ in $\mathbb{R}$, and $x$ in $M'_0\cap M_1$, we have :
\begin{multline*}
\sigma_t^{\chi'_2}(s(x))=\Delta_{\chi_2}^{-it}J_{\chi_2}x^*J_{\chi_2}\Delta_{\chi_2}^{it}\\=J_{\chi_2}\Delta_{\chi_2}^{-it}x^*\Delta_{\chi_2}^{it}J_{\chi_2}=s(\sigma_{-t}^{\chi_2}(x))=s(\sigma_{-t}^{\chi}(x))
\end{multline*}
and we have seen in \ref{propsigma2} that $\sigma_t^{\chi'_2}(\hat{r}(x))=\hat{r}(\sigma_{-t}^{T_1}(x))$.
So, we have :
\begin{multline*}
\sigma_t^{\chi'_2\circ j}(s(x))=j\circ\sigma_{-t}^{\chi'_2}\circ j(s(x))\\=
j\circ\sigma_{-t}^{\chi'_2}(\hat{r}(x))=j\circ\hat{r}(\sigma_{t}^{T_1}(x))=s(\sigma_{t}^{T_1}(x))
\end{multline*}
and, therefore :
\[\delta'^{it}s(x)\delta'^{it}=\sigma_t^{\chi'_2\circ j}\sigma_{-t}^{\chi'_2}(s(x))=
\sigma_t^{\chi'_2\circ j}(s(\sigma_{t}^{\chi}(x))=s(\sigma_{t}^{T_1}\sigma_{t}^{\chi}(x))\]
and an analogous calculation for $\delta'^{it}\hat{r}(x)\delta'^{it}$ will lead to the possibility of defining $\delta'^{it}{}_s\underset{\chi}{\otimes}{}_{\hat{r}}\delta'^{it}$. 
\newline
On the other hand, by \ref{modulus}, that, for all $x$ in $\pi(M'_0\cap M_2)'$, for all $t$ in $\mathbb{R}$, we have $\delta'^{it}x\delta'^{-it}=\sigma_t^{\chi'_2\circ j}\sigma_{-t}^{\chi'_2}(x)$. 
\newline
Therefore, the first formula is a corollary of \ref{proptau2}. If, in this formula, we take $x=\delta'^{is}$, we obtain the second result. 
\end{proof}

\subsection{Lemma}
\label{lemmodulus}
 {\it (i) Let $b$ be in $\gN_{\tilde{T'_2}}\cap\gN_{\chi'_2}\cap\gN_{\chi'_2\circ j}$, and let $X$ positive affiliated to $\pi(M'_0\cap M_2)'$ be such that $\delta'^{-1/2}X\delta'^{-1/2}$ is bounded; then, the element of the extended positive part  $(id{}_s\underset{\chi}{*}{}_{\hat{r}}\chi'_2)\Gamma (X)$ is such that :}
\begin{multline*}
<(id{}_s\underset{\chi}{*}{}_{\hat{r}}\chi'_2)\Gamma (X), \omega_{J_{\chi_2}\Lambda_{\chi'_2}(b)}>\\=
<j\circ\tilde{T'_2}\circ j (\delta'^{-1/2}X\delta'^{-1/2}), \omega_{\delta'^{1/2}J_{\chi_2}\Lambda_{\chi'_2}(b)}>
\end{multline*}
{\it If $X$ is bounded, such that $\delta'^{-1/2}X\delta'^{-1/2}$ is bounded and in $\gM_{j\circ\tilde{T'_2}\circ j}^+$, then $(\omega_{J_{\chi_2}\Lambda_{\chi'_2}(b)}{}_s\underset{\chi}{*}{}_{\hat{r}}id)\Gamma (X)$ belongs to  $\gM_{\tilde{T'_2}}^+\cap\gM_{\chi'_2}^+$. 
\newline
(ii) If $Y$ is in $\gM_{j\circ\tilde{T'_2}\circ j}^+$, we have :}
\[\delta'^{1/2}(j\circ\tilde{T'_2}\circ j)(Y)\delta'^{1/2}=(id{}_s\underset{\chi}{*}{}_{\hat{r}}\chi'_2)\Gamma (\delta'^{1/2}Y\delta'^{1/2})\]

\begin{proof}
Let us suppose first that $a$ and $b$ belong to $\gN_{\tilde{T'_2}}\cap\gN_{\chi'_2}$; using \ref{proptau}(v), we have :
\[j(\omega_{J_{\chi_2}\Lambda_{\chi'_2}(a)}{}_s\underset{\chi}{*}{}_{\hat{r}}id)\Gamma (b^*b)=
(\omega_{J_{\chi_2}\Lambda_{\chi'_2}(b)}{}_s\underset{\chi}{*}{}_{\hat{r}}id)\Gamma (a^*a)\]
and, therefore, applying the weight $\chi'_2$, and using \ref{thinv}, we obtain :
\[<(id{}_s\underset{\chi}{*}{}_{\hat{r}}\chi'_2)\Gamma (a^*a), \omega_{J_{\chi_2}\Lambda_{\chi'_2}(b)}>=
<j\circ\tilde{T'_2}\circ j (b^*b), \omega_{J_{\chi_2}\Lambda_{\chi'_2}(a)}>\]
Let us suppose now that $a$ and $b$ belong to $\gN_{\tilde{T'_2}}\cap\gN_{\chi'_2}\cap\gN_{\chi'_2\circ j}$ (which is dense by \ref{cormodulus2}); we get, using [V3], that $J_{\chi'_2\circ j}\Lambda_{\chi'_2\circ j}(a)$  is in the domain of $\delta'^{-1/2}$, and that $\delta^{-1/2}J_{\chi'_2\circ j}\Lambda_{\chi'_2\circ j}(a)$ can be indentified with $\pi(\lambda)^{i/4}J_{\chi_2}\Lambda_{\chi'_2}(a)$. So we get that :
\[<(id{}_s\underset{\chi}{*}{}_{\hat{r}}\chi'_2)\Gamma (a^*a), \omega_{J_{\chi_2}\Lambda_{\chi'_2}(b)}>=
<j\circ\tilde{T'_2}\circ j (b^*b), \omega_{\delta'^{-1/2}J_{\chi'_2\circ j}\Lambda_{\chi'_2\circ j}(a)}>\]
Let us suppose now that $a\delta'^{-1/2}$ is bounded (we shall then denote again $a\delta'^{-1/2}$ its closure, which, using again [V3] belongs to $\gN_{\chi'_2\circ j}$). We then get that :
\[<(id{}_s\underset{\chi}{*}{}_{\hat{r}}\chi'_2)\Gamma (a^*a), \omega_{J_{\chi_2}\Lambda_{\chi'_2}(b)}>=
<j\circ\tilde{T'_2}\circ j (b^*b), \omega_{J_{\chi'_2\circ j}\Lambda_{\chi'_2}(a\delta'^{-1/2})}>\]
And, using \ref{propT}, it is equal to :
\[<j\circ\tilde{T'_2}\circ j (\delta'^{-1/2}a^*a\delta'^{-1/2}), \omega_{J_{\chi'_2\circ j}\Lambda_{\chi'_2\circ j}(b)}>\]
or, using again these identifications, to :
\[<j\circ\tilde{T'_2}\circ j (\delta'^{-1/2}a^*a\delta'^{-1/2}), \omega_{\delta'^{1/2}J_{\chi_2}\Lambda_{\chi'_2}(b)}>\]
If $X$ is positive such that $\delta'^{-1/2}X\delta'^{-1/2}$ is bounded, we may consider $X$ as the upper limit of elements of the type $a_i^*a_i$, where the operators $a_i$ belong to the dense left ideal $\gN_{\tilde{T'_2}}\cap\gN_{\chi'_2}\cap\gN_{\chi'_2\circ j}$; then every $a_i\delta'^{-1/2}$ is bounded, and
we get the first formula of (i) by increasing limits. 
\newline
The second formula of (i) is just a corollary, using \ref{fiber}. 
The proof of (ii) is easy, using the fact that we are in an essential domain of $\delta'^{1/2}$, thanks to \ref{cormodulus2}. 
\end{proof}

\subsection{Theorem}
\label{thmodulus}
{\it We have $\Gamma(\delta')=\delta'{}_s\underset{\chi}{\otimes}{}_{\hat{r}}\delta'$. }

\begin{proof}
Applying $\Gamma$ to \ref{cormodulus}(ii), we get, for all $Y$ in $\gM_{j\circ\tilde{T'_2}\circ j}^+$:
\[\Gamma(\delta'^{1/2})((j\circ\tilde{T'_2}\circ j)(Y)_s\underset{\chi}{\otimes}{}_{\hat{r}}1)\Gamma(\delta'^{1/2})
=\Gamma((id{}_s\underset{\chi}{*}{}_{\hat{r}}\chi'_2)\Gamma (\delta'^{1/2}Y\delta'^{1/2}))\]
which is equal to :
\[((id{}_s\underset{\chi}{*}{}_{\hat{r}}id{}_s\underset{\chi}{*}{}_{\hat{r}}\chi'_2)(\Gamma_s\underset{\chi}{*}{}_{\hat{r}}id)\Gamma (\delta'^{1/2}Y\delta'^{1/2}))\]
and, therefore, to :
\[((id{}_s\underset{\chi}{*}{}_{\hat{r}}id{}_s\underset{\chi}{*}{}_{\hat{r}}\chi'_2)(id_s\underset{\chi}{*}{}_{\hat{r}}\Gamma)\Gamma (\delta'^{1/2}Y\delta'^{1/2}))\]
Let now $b$ be in the set $T_{\chi'_2, \tilde{T'_2}}^{\chi'_2\circ j}$ defined in \ref{cormodulus2}, and let define $Z$ by :
\[Z=(\omega_{J_{\chi_2}\Lambda_{\chi'_2}(b)}{}_s\underset{\chi}{*}{}_{\hat{r}}id)\Gamma(\delta'^{1/2}Y\delta'^{1/2})\]
 Then, we get, using \ref{propdelta}, that :
\begin{multline*}
\delta'^{-1/2}Z\delta'^{-1/2}=
\delta'^{-1/2}(\omega_{J_{\chi_2}\Lambda_{\chi'_2}(b)}{}_s\underset{\chi}{*}{}_{\hat{r}}id)\Gamma(\delta'^{1/2}Y\delta'^{1/2})\delta'^{-1/2}\\
=(\omega_{\delta'^{1/2}J_{\chi_2}\Lambda_{\chi'_2}(b)}{}_s\underset{\chi}{*}{}_{\hat{r}}id)\Gamma(Y)
\end{multline*}
which is bounded, by \ref{cormodulus2}.  So, we can use \ref{lemmodulus}(i) and we get, for $b'$ in $\gN_{\chi'_2}\cap\gN_{\tilde{T'_2}}$,  that :
\[<\Gamma(\delta'^{1/2})(j\circ\tilde{T'_2}\circ j)(Y)_s\underset{\chi}{\otimes}{}_{\hat{r}}1)\Gamma(\delta'^{1/2}), \omega_{J_{\chi_2}\Lambda_{\chi'_2}(b){}_s\underset{\chi}{\otimes}{}_{\hat{r}}J_{\chi_2}\Lambda_{\chi'_2}(b')}>\]
is equal to :
\[<(id{}_s\underset{\chi}{*}{}_{\hat{r}}\chi'_2)\Gamma((\omega_{J_{\chi_2}\Lambda_{\chi'_2}(b)}{}_s\underset{\chi}{*}{}_{\hat{r}}id)\Gamma(\delta'^{1/2}Y\delta'^{1/2})), \omega_{J_{\chi_2}\Lambda_{\chi'_2}(b')}>\]
which is equal to :
\[<j\circ\tilde{T'_2}\circ j (\delta'^{-1/2}Z\delta'^{-1/2}), \omega_{\delta'^{1/2}J_{\chi_2}\Lambda_{\chi'_2}(b')}>\]
or to :
\[<(j\circ\tilde{T'_2}\circ j (Y)_s\underset{\chi}{\otimes}{}_{\hat{r}}1), \omega_{\delta'^{1/2}J_{\chi_2}\Lambda_{\chi'_2}(b){}_s\underset{\chi}{\otimes}{}_{\hat{r}}\delta'^{1/2}J_{\chi_2}\Lambda_{\chi'_2}(b')}>\]
from which we infer, by increasing limits, that :
\[<\Gamma(\delta'), \omega_{J_{\chi_2}\Lambda_{\chi'_2}(b){}_s\underset{\chi}{\otimes}{}_{\hat{r}}J_{\chi_2}\Lambda_{\chi'_2}(b')}>=\|\delta'^{1/2}J_{\chi_2}\Lambda_{\chi'_2}(b){}_s\underset{\chi}{\otimes}{}_{\hat{r}}\delta'^{1/2}J_{\chi_2}\Lambda_{\chi'_2}(b')\|^2\]
which, by \ref{cormodulus2} again, finishes the proof.  \end{proof}

\section{A density theorem}
\label{density}
In that chapter, we prove that there are sufficiently enough operators which are both bounded under the left-invariant operator-valued weight and the right-invariant operator-valued weight (\ref{thdens1}); this fact was proved before only in the case when the basis was semi-finite. This allows, as a corollary, to found bounded elements for both $s$ and $\hat{r}$ (\ref{cordens}), which will be usefull for duality. 

\subsection{Lemma}
\label{lemdens1}
{\it (i) Let $y$, $z$, $z'$ in $\gN_{\tilde{T'_2}}\cap\gN_{\chi'_2}$; then, we have :
\begin{multline*}
(\omega_{\Lambda_{\chi'_2}(y), J_{\chi_2}z^*\Lambda_{\chi'_2}(z')}*id)(W)^*
((\omega_{\Lambda_{\chi'_2}(y), J_{\chi_2}z^*\Lambda_{\chi'_2}(z')}*id)(W))\\
\leq \|\tilde{T'_2}(z'^*z')\|
(\omega_{J_{\chi_2}\Lambda_{\chi'_2}(z)}{}_s\underset{\chi}{*}{}_{\hat{r}}id)\Gamma (y^*y)
\end{multline*}
(ii) Let $y$, $z'$ in $\gN_{\tilde{T'_2}}\cap\gN_{\chi'_2}$, and $z$ in $\pi(M'_0\cap M_2)'$; then, we have :}
\begin{multline*}
j((\omega_{\Lambda_{\chi'_2}(y), J_{\chi_2}z^*\Lambda_{\chi'_2}(z')}*id)(W)^*
((\omega_{\Lambda_{\chi'_2}(y), J_{\chi_2}z^*\Lambda_{\chi'_2}(z')}*id)(W)))\\
\leq \|\tilde{T'_2}(z'^*z')\|
(\omega_{J_{\chi_2}\Lambda_{\chi'_2}(y)}{}_s\underset{\chi}{*}{}_{\hat{r}}id)\Gamma (z^*z)
\end{multline*}

\begin{proof}
Using \ref{proptau}(iii), we get that :
\[(\omega_{\Lambda_{\chi'_2}(y), J_{\chi_2}z^*\Lambda_{\chi'_2}(z')}*id)(W)^*
(\omega_{\Lambda_{\chi'_2}(y), J_{\chi_2}z^*\Lambda_{\chi'_2}(z')}*id)(W)\]
is equal to :
\begin{multline*}
(\omega_{J_{\chi_2}\Lambda_{\chi'_2}(z), J_{\chi_2}\Lambda_{\chi'_2}(z')}{}_s\underset{\chi}{*}{}_{\hat{r}}id)\Gamma(y)^*
(\omega_{J_{\chi_2}\Lambda_{\chi'_2}(z), J_{\chi_2}\Lambda_{\chi'_2}(z')}{}_s\underset{\chi}{*}{}_{\hat{r}}id)\Gamma(y)\\
=\lambda^{s, \hat{r}*}_{J_{\chi_2}\Lambda_{\chi'_2}(z)}\Gamma(y^*)
\lambda^{s, \hat{r}}_{J_{\chi_2}\Lambda_{\chi'_2}(z')}\lambda^{s, \hat{r}*}_{J_{\chi_2}\Lambda_{\chi'_2}(z')}
\Gamma(y)\lambda^{s, \hat{r}}_{J_{\chi_2}\Lambda_{\chi'_2}(z)}
\end{multline*}
which is less than :
\[\|R^{s, \chi^o}(J_{\chi_2}\Lambda_{\chi'_2}(z'))\|^2
(\omega_{J_{\chi_2}\Lambda_{\chi'_2}(z)}{}_s\underset{\chi}{*}{}_{\hat{r}}id)\Gamma (y^*y)\]
from which we get (i), using \ref{proptau}(i). 
\newline
Moreover, we obtain, using \ref{proptau}(v), that :
\begin{multline*}
j((\omega_{\Lambda_{\chi'_2}(y), J_{\chi_2}z^*\Lambda_{\chi'_2}(z')}*id)(W)^*
(\omega_{\Lambda_{\chi'_2}(y), J_{\chi_2}z^*\Lambda_{\chi'_2}(z')}*id)(W))\\
\leq \|\tilde{T'_2}(z'^*z')\|
j((\omega_{J_{\chi_2}\Lambda_{\chi'_2}(z)}{}_s\underset{\chi}{*}{}_{\hat{r}}id)\Gamma (y^*y))\\
=\|\tilde{T'_2}(z'^*z')\|(\omega_{J_{\chi_2}\Lambda_{\chi'_2}(y)}{}_s\underset{\chi}{*}{}_{\hat{r}}id)\Gamma (z^*z)
\end{multline*}
Let us suppose now that $z$ belongs to $\pi(M'_0\cap M_2)'$; using Kaplansky's theorem, there exist a family $z_i$ in $\gN_{\tilde{T'_2}}\cap\gN_{\chi'_2}$, weakly converging to $z$, with $\|z_i\|\leq\|z\|$; from which we infer that $R^{s, \chi^o}(J_{\chi_2}z_i^*\Lambda_{\chi'_2}(z'))$ is weakly converging to $R^{s, \chi^o}(J_{\chi_2}z^*\Lambda_{\chi'_2}(z'))$, with :
\[\|R^{s, \chi^o}(J_{\chi_2}z_i^*\Lambda_{\chi'_2}(z'))\|\leq\|R^{s, \chi^o}(J_{\chi_2}z^*\Lambda_{\chi'_2}(z'))\|\]
and, therefore, that 
$(\omega_{\Lambda_{\chi'_2}(y), J_{\chi_2}z_i^*\Lambda_{\chi'_2}(z')}*id)(W)$ is weakly converging to 
$(\omega_{\Lambda_{\chi'_2}(y), J_{\chi_2}z^*\Lambda_{\chi'_2}(z')}*id)(W)$, with :
\[\|(\omega_{\Lambda_{\chi'_2}(y), J_{\chi_2}z_i^*\Lambda_{\chi'_2}(z')}*id)(W)\|\leq\|(\omega_{\Lambda_{\chi'_2}(y), J_{\chi_2}z^*\Lambda_{\chi'_2}(z')}*id)(W)\|\]
Therefore, we obtain, for such a $z$, that we still have :
\begin{multline*}
j((\omega_{\Lambda_{\chi'_2}(y), J_{\chi_2}z^*\Lambda_{\chi'_2}(z')}*id)(W)^*
(\omega_{\Lambda_{\chi'_2}(y), J_{\chi_2}z^*\Lambda_{\chi'_2}(z')}*id)(W))\\
\leq \|\tilde{T'_2}(z'^*z')\|(\omega_{J_{\chi_2}\Lambda_{\chi'_2}(y)}{}_s\underset{\chi}{*}{}_{\hat{r}}id)\Gamma (z^*z)
\end{multline*}
which finishes the proof. \end{proof}

\subsection{Proposition}
\label{propdens1}
{\it If $z$ belongs to $\gN_{j\circ\tilde{T'_2}\circ j}$, $y$ and $z'$ to $\gN_{\tilde{T'_2}}\cap\gN_{\chi'_2}$, then 
$(\omega_{\Lambda_{\chi'_2}(y), J_{\chi_2}z^*\Lambda_{\chi'_2}(z')}*id)(W)$ belongs to $\gN_{\tilde{T'_2}}\cap\gN_{\chi'_2}$.}

\begin{proof}
Using \ref{lemdens1} and \ref{thinv}(ii), we get that :
\[\chi'_2((\omega_{\Lambda_{\chi'_2}(y), J_{\chi_2}z^*\Lambda_{\chi'_2}(z')}*id)(W)^*(\omega_{\Lambda_{\chi'_2}(y), J_{\chi_2}z^*\Lambda_{\chi'_2}(z')}*id)(W))\]
is less or equal to :
\[\|T'_2(z'^*z)\|(j\circ\tilde{T'_2}\circ j(z^*z)J_{\chi_2}\Lambda_{\chi'_2}(y)|J_{\chi_2}\Lambda_{\chi'_2}(y))\]
which proves that $(\omega_{\Lambda_{\chi'_2}(y), J_{\chi_2}z^*\Lambda_{\chi'_2}(z')}*id)(W)$ belongs to $\gN_{\chi'_2}$.
\newline
With same arguments, we get that :
\[j\circ\tilde{T'_2}((\omega_{\Lambda_{\chi'_2}(y), J_{\chi_2}z^*\Lambda_{\chi'_2}(z')}*id)(W)^*(\omega_{\Lambda_{\chi'_2}(y), J_{\chi_2}z^*\Lambda_{\chi'_2}(z')}*id)(W))\]
is majorized by $\|T'_2(z'^*z)\|\|j\circ\tilde{T'_2}\circ j(z^*z)\|\|T'_2(y^*y)\|$, which proves that $(\omega_{\Lambda_{\chi'_2}(y), J_{\chi_2}z^*\Lambda_{\chi'_2}(z')}*id)(W)$ belongs to $\gN_{\tilde{T'_2}}$. \end{proof}

\subsection{Lemma}
\label{lemdens2}
{\it Let $y_1$, $y_2$, $z$ in $\gN_{\tilde{T'_2}}\cap\gN_{\chi'_2}$; then we have :}
\begin{multline*}
j((\omega_{y_1^*\Lambda_{\chi'_2}(y_2), J_{\chi_2}\Lambda_{\chi'_2}(z)}*id)(W)^*)^*
j((\omega_{y_1^*\Lambda_{\chi'_2}(y_2), J_{\chi_2}\Lambda_{\chi'_2}(z)}*id)(W)^*)\\
\leq\|\tilde{T'_2}(y_1^*y_1)\|(\omega_{J_{\chi_2}\Lambda_{\chi'_2}(y_2)}{}_s\underset{\chi}{*}{}_{\hat{r}}id)\Gamma (zz^*)
\end{multline*}
\begin{proof}
If $y$, $z$ are in $\gN_{\tilde{T'_2}}$, we have, using \ref{thtau2} and \ref{thW} :
\begin{eqnarray*}
j((\omega_{\Lambda_{\chi'_2}(y), J_{\chi_2}\Lambda_{\chi'_2}(z)}*id)(W)^*)
&=&\hat{J}(\omega_{\Lambda_{\chi'_2}(y), J_{\chi_2}\Lambda_{\chi'_2}(z)}*id)(W)\hat{J}\\
&=&(\omega_{J_{\chi_2}\Lambda_{\chi'_2}(y), \Lambda_{\chi'_2}(z)}*id)(W^*)\\
&=&(\omega_{\Lambda_{\chi'_2}(z), J_{\chi_2}\Lambda_{\chi'_2}(y)}*id)(W)^*
\end{eqnarray*}
and, we get that :
\[j((\omega_{y_1^*\Lambda_{\chi'_2}(y_2), J_{\chi_2}\Lambda_{\chi'_2}(z)}*id)(W)^*)^*
j((\omega_{y_1^*\Lambda_{\chi'_2}(y_2), J_{\chi_2}\Lambda_{\chi'_2}(z)}*id)(W)^*)\]
is equal to :
\[(\omega_{\Lambda_{\chi'_2}(z), J_{\chi_2}y_1^*\Lambda_{\chi'_2}(y_2)}*id)(W)
(\omega_{\Lambda_{\chi'_2}(z), J_{\chi_2}y_1^*\Lambda_{\chi'_2}(y_2)}*id)(W)^*\]
and, thanks to \ref{proptau}(iii), is equal to :
\[(\omega_{J_{\chi_2}\Lambda_{\chi'_2}(y_1), J_{\chi_2}\Lambda_{\chi'_2}(y_2)}{}_s\underset{\chi}{*}{}_{\hat{r}}id)\Gamma(z)
(\omega_{J_{\chi_2}\Lambda_{\chi'_2}(y_1), J_{\chi_2}\Lambda_{\chi'_2}(y_2)}{}_s\underset{\chi}{*}{}_{\hat{r}}id)\Gamma(z))^*\]
which can be written as :
\[\lambda^{s, \hat{r}*}_{J_{\chi_2}\Lambda_{\chi'_2}(y_2)}\Gamma (z)
\lambda^{s, \hat{r}}_{J_{\chi_2}\Lambda_{\chi'_2}(y_1)}
\lambda^{s, \hat{r}*}_{J_{\chi_2}\Lambda_{\chi'_2}(y_1)}
\Gamma(z^*)\lambda^{s, \hat{r}}_{J_{\chi_2}\Lambda_{\chi'_2}(y_2)}\]
which is less than :
\[\|R^{s, \chi^o}(J_{\chi_2}\Lambda_{\chi'_2}(y_1))\|^2
(\omega_{J_{\chi_2}\Lambda_{\chi'_2}(y_2)}{}_s\underset{\chi}{*}{}_{\hat{r}}id)\Gamma (zz^*)\]
which gives the result, by \ref{proptau}(i). \end{proof}

\subsection{Proposition}
\label{propdens2}
{\it Let $y_1$, $z'$ in $\gN_{\tilde{T'_2}}\cap\gN_{\chi'_2}$, $y_2$ in $\gN_{\tilde{T'_2}}\cap\gN_{\chi'_2}\cap\gN_{\chi'_2\circ j}$, $z$ in $j(T_{\chi'_2, \tilde{T'_2}}^{\chi'_2\circ j})^*$ (where $T_{\chi'_2, \tilde{T'_2}}^{\chi'_2\circ j}$ has been defined in \ref{cormodulus2}), and $e_n$ the analytic elements associated to the Radon-Nykodym derivative $\delta'$, defined in [V3]); then the operators 
$(\omega_{y_1^*\Lambda_{\chi'_2}(y_2), J_{\chi_2}z^*e_n^*\Lambda_{\chi'_2}(z')}*id)(W)$ belong to $\gN_{\tilde{T'_2}}\cap\gN_{\chi'_2}\cap\gN_{j\circ\tilde{T'_2}\circ j}\cap\gN_{\chi'_2\circ j}$. }

\begin{proof}
Let us write $X=(\omega_{y_1^*\Lambda_{\chi'_2}(y_2), J_{\chi_2}z^*e_n^*\Lambda_{\chi'_2}(z')}*id)(W)$. 
\newline
As $y_1^*y_2$ belongs to $\gN_{\tilde{T'_2}}\cap\gN_{\chi'_2}$, and $z$ belongs to $j(\gN_{\tilde{T'_2}}^*)=\gN_{j\circ\tilde{T'_2}\circ j}$, and, therefore, $e_nz$ belongs to $\gN_{j\circ\tilde{T'_2}\circ j}$, we get, using \ref{propdens1}, that $X$ belongs to $\gN_{\tilde{T'_2}}\cap\gN_{\chi'_2}$. 
\newline
On the other hand, as $y_1$, $y_2$, $z^*e_n^*z'$ belong to $\gN_{\tilde{T'_2}}\cap\gN_{\chi'_2}$, we can use \ref{lemdens2}, and we get that :
\begin{eqnarray*}
j(X^*)^*j(X^*)
&\leq&\|\tilde{T'_2}(y_1^*y_1)\|(\omega_{J_{\chi_2}\Lambda_{\chi'_2}(y_2)}{}_s\underset{\chi}{*}{}_{\hat{r}}id)\Gamma (z^*e_n^*z'z'^*e_nz)\\
&\leq&\|\tilde{T'_2}(y_1^*y_1)\|\|z'\|^2(\omega_{J_{\chi_2}\Lambda_{\chi'_2}(y_2)}{}_s\underset{\chi}{*}{}_{\hat{r}}id)\Gamma (z^*e_n^*e_nz)
\end{eqnarray*}
Let us apply $\tilde{T'_2}$ to this inequality; we get that :
\[\tilde{T'_2}(j(X^*)^*j(X^*))
\leq\|\tilde{T'_2}(y_1^*y_1)\|\|z'\|^2 \tilde{T'_2}((\omega_{J_{\chi_2}\Lambda_{\chi'_2}(y_2)}{}_s\underset{\chi}{*}{}_{\hat{r}}id)\Gamma (z^*e_n^*e_nz))\]
which is equal, thanks to \ref{lemmodulus}(i), to :
\[\|\tilde{T'_2}(y_1^*y_1)\|\|z'\|^2< j\circ\tilde{T'_2}\circ j(\delta'^{-1/2}z^*e_n^*e_nz\delta'^{-1/2}), \omega_{\delta'^{-1/2}J_{\chi_2}\Lambda_{\chi'_2}(y_2)}>\]
With the hypothesis, we get that $\delta'^{1/2}z\delta'^{-1/2}$ belongs to $\gN_{j\circ\tilde{T'_2}\circ j}$, and, therefore $e_nz\delta'^{-1/2}=(e_n\delta^{-1/2})\delta'^{1/2}z\delta'^{-1/2}$ belongs also to $\gN_{j\circ\tilde{T'_2}\circ j}$. We also get that $J_{\chi_2}\Lambda_{\chi'_2}(y_2)$ belongs to $D(\delta^{-1/2})$, which proves that $j(X^*)$ belongs to $\gN_{\tilde{T'_2}}$, and, therefore, that $X$ belongs to $\gN_{j\circ\tilde{T'_2}\circ j}$.  We prove by similar calculations that $j(X^*)$ belongs to $\gN_{\chi'_2}$, and, therefore, that $X$ belongs to $\gN_{\chi'_2\circ j}$. \end{proof}

\subsection{Theorem}
\label{thdens1}
{\it The left ideal $\gN_{\tilde{T'_2}}\cap\gN_{\chi'_2}\cap\gN_{j\circ\tilde{T'_2}\circ j}\cap\gN_{\chi'_2\circ j}$ is dense in $\pi(M'_0\cap M_2)'$, and the linear space $\Lambda_{\chi'_2}(\gN_{\tilde{T'_2}}\cap\gN_{\chi'_2}\cap\gN_{j\circ\tilde{T'_2}\circ j}\cap\gN_{\chi'_2\circ j})$ is dense in $\gH$.}

\begin{proof}
Let $y$ be in $\gN_{\tilde{T'_2}}\cap\gN_{\chi'_2}\cap\gN_{\chi'_2\circ j}$, and $z$ in $\gN_{\tilde{T'_2}}\cap\gN_{\chi'_2}$. Taking, by Kaplansky's theorem, a bounded family $e_i$ in $\gN_{\tilde{T'_2}}\cap\gN_{\chi'_2}$ strongly converging to $1$, we get that $R^{r, \chi}(e_i^*\Lambda_{\chi'_2}(y))$ is weakly converging to $R^{r, \chi}(\Lambda_{\chi'_2}(y))$; taking also a bounded family $f_k$ in $j(T_{\chi'_2, \tilde{T'_2}}^{\chi'_2\circ j})^*$ strongly converging to $1$, we get that $R^{s, \chi^o}(J_{\chi_2}f_k^*e_n^*\Lambda_{\chi'_2}(z))$ is weakly converging, when $n$ and $k$ go to infinity, to $R^{r, \chi^o}(\Lambda_{\chi'_2}(z))$. 
\newline
Therefore, using \ref{propdens2}, we get that $(\omega_{\Lambda_{\chi'_2}(y), J_{\chi_2}\Lambda_{\chi'_2}(z)}*id)(W)$ belongs to the weak closure of $\gN_{\tilde{T'_2}}\cap\gN_{\chi'_2}\cap\gN_{j\circ\tilde{T'_2}\circ j}\cap\gN_{\chi'_2\circ j}$. 
\newline
Using then \ref{cormodulus2}, we get that, for any $x$  in $T_{\chi'_2, \tilde{T'_2}}$, there exists $y_i$  in $\gN_{\tilde{T'_2}}\cap\gN_{\chi'_2}\cap\gN_{\chi'_2\circ j}$ such that $\Lambda_{\tilde{T'_2}}(y_i)$ is weakly converging to $\Lambda_{\tilde{T'_2}}(x)$, or, equivalently, that $R^{r, \chi}(\Lambda_{\chi'_2}(y_i))$ is weakly converging to $R^{r, \chi}(\Lambda_{\chi'_2}(x))$. Therefore, we get that $(\omega_{\Lambda_{\chi'_2}(x), J_{\chi_2}\Lambda_{\chi'_2}(z)}*id)(W)$ belongs to the weak closure of $\gN_{\tilde{T'_2}}\cap\gN_{\chi'_2}\cap\gN_{j\circ\tilde{T'_2}\circ j}\cap\gN_{\chi'_2\circ j}$.
\newline
Using \ref{prop}(ii), this result remains true for $x$ in $\gN_{\tilde{T'_2}}\cap\gN_{\chi'_2}\cap\gN_{\tilde{T'_2}}^*\cap\gN_{\chi'_2}^*$.
\newline
If now $x$ belongs to $\gN_{\tilde{T'_2}}\cap\gN_{\chi'_2}$, and $h_i$ is a bounded family in $\gN_{\tilde{T'_2}}\cap\gN_{\chi'_2}$, as $\Lambda_{\tilde{T'_2}}(h_i^*x)=h_i^*\Lambda_{\tilde{T'_2}}(x)$ is weakly converging to $\Lambda_{\tilde{T'_2}}(x)$, we finaly obtain that, for any $x$, $z$ in $\gN_{\tilde{T'_2}}\cap\gN_{\chi'_2}$, the operator $(\omega_{\Lambda_{\chi'_2}(x), J_{\chi_2}\Lambda_{\chi'_2}(z)}*id)(W)$ belongs to the weak closure of $\gN_{\tilde{T'_2}}\cap\gN_{\chi'_2}\cap\gN_{j\circ\tilde{T'_2}\circ j}\cap\gN_{\chi'_2\circ j}$. 
 \newline
 Using then \ref{prop}(i), we get that, for all $\xi$ in $D(\gH_r, \chi^o)$, and $\eta$ in $D({}_s\gH, \chi)$, the operator $(\omega_{\xi, \eta}*id)(W)$ belongs to the weak closure of $\gN_{\tilde{T'_2}}\cap\gN_{\chi'_2}\cap\gN_{j\circ\tilde{T'_2}\circ j}\cap\gN_{\chi'_2\circ j}$.  Which proves the density of $\gN_{\tilde{T'_2}}\cap\gN_{\chi'_2}\cap\gN_{j\circ\tilde{T'_2}\circ j}\cap\gN_{\chi'_2\circ j}$ in $\pi(M'_0\cap M_2)'$, by \ref{exW}. 
\newline
Let $g_n$ an increasing sequence of positive elements of  $\gM_{\tilde{T'_2}}\cap\gM_{j\circ\tilde{T'_2}\circ j}\cap\gM_{\chi'_2}\cap\gM_{\chi'_2\circ j}$ strongly converging to $1$; the operators 
\[h_n=\sqrt{\frac{1}{\pi}}\int_{-\infty}^{\infty}e^{-t^2}\sigma_t^{\chi'_2}(g_n)dt\]
are in $\gM_{\tilde{T'_2}}\cap\gM_{\chi'_2}$, analytic with respect to $\chi'_2$, and, for any $z$ in $\Bbb{C}$, $\sigma_z^{\chi'_2}(h_n)$ is a bounded sequence  strongly converging to $1$. 
\newline
Let now $\lambda=\int_0^{\infty}tde_t$ be the scaling operator, affiliated to $Z(M_0)\cap Z(M_1)$, and let us write $h'_n=\pi(\int_{1/n}^nde_t)h_n$; these operators are in $\gN_{\tilde{T'_2}}\cap\gN_{\chi'_2}$, analytic with respect to $\chi'_2$, and, for any $z$ in $\Bbb{C}$, $\sigma_z^{\chi'_2}(h_n)$ is a bounded sequence strongly converging to $1$. 
Moreover, using \ref{modulus} and \ref{cormodulus}, we get that the operators $h'_n$ belong also to $\gN_{j\circ\tilde{T'_2}\circ j}\cap\gN_{\chi'_2\circ j}$. 
\newline
Let now $x$ be in $\gN_{\chi'_2}$; we get that $xh'_n\in\gN_{\tilde{T'_2}}\cap\gN_{\chi'_2}\cap\gN_{j\circ\tilde{T'_2}\circ j}\cap\gN_{\chi'_2\circ j}$, and that 
\[\Lambda_{\chi'_2}(xh'_n)=J_{\chi_2}\sigma_{-i/2}^{\chi'_2}(h'_n)J_{\chi_2}\Lambda_{\chi'_2}(x)\rightarrow_n\Lambda_{\chi'_2}(x)\]
which finishes the proof.  \end{proof}

\subsection{Theorem}
\label{thdens2}
{\it Let $T_{\chi'_2, \tilde{T'_2}, \chi'_2\circ j, j\circ\tilde{T'_2}\circ j}$ be the subset of elements $x$ in $\gN_{\tilde{T'_2}}\cap\gN_{\chi'_2}\cap\gN_{j\circ\tilde{T'_2}\circ j}\cap\gN_{\chi'_2\circ j}$, analytic with respect to both $\chi'_2$ and $\chi'_2\circ j$, and such that, for all $z$, $z'$ in $\Bbb{C}$, $\sigma_z^{\chi'_2}\circ\sigma_{z'}^{\chi'_2\circ j}(x)$ belongs to $\gN_{\tilde{T'_2}}\cap\gN_{\chi'_2}\cap\gN_{j\circ\tilde{T'_2}\circ j}\cap\gN_{\chi'_2\circ j}$. Then $T_{\chi'_2, \tilde{T'_2}, \chi'_2\circ j, j\circ\tilde{T'_2}\circ j}$ is dense in $\pi(M'_0\cap M_2)'$, and $\Lambda_{\chi'_2}(T_{\chi'_2, \tilde{T'_2}, \chi'_2\circ j, j\circ\tilde{T'_2}\circ j})$ is dense in $\gH$. }

\begin{proof}
Let $x$ be a positive operator in $\gM_{\chi'_2}\cap\gM_{\tilde{T'_2}}\cap\gM_{\chi'_2\circ j}\cap\gM_{j\circ\tilde{T'_2}\circ j}$; let now $\lambda=\int_0^{\infty}tde_t$ be the scaling operator, affiliated to $Z(M_0)\cap Z(M_1)$, and let us define :
\[x_n=\pi(\int_{1/n}^nde_t)\frac{n}{\pi}\int_{-\infty}^{\infty}\int_{-\infty}^{\infty}e^{-n(t^2+s^2)}\sigma_t^{\chi'_2}\sigma_s^{\chi'_2\circ j}(x)dsdt\]
By similar arguments as in \ref{lemma3}, we can show that $x_n$ is analytic both with respect to $\chi'_2$ and $\chi'_2\circ j$. 
\newline
Using \ref{modulus}(iii) and \ref{cormodulus} and the same technics as in ([EN], 10.12), we see that, thanks to $\pi(\int_{1/n}^nde_t)$, the operators  
$\sigma_z^{\chi'_2}(x_n)$ and $\sigma_z^{\chi'_2\circ j}(x_n)$ are linear combinations of positive elements in $\gM_{\chi'_2}\cap\gM_{\tilde{T'_2}}\cap\gM_{\chi'_2\circ j}\cap\gM_{j\circ\tilde{T'_2}\circ j}$. 
\end{proof}

\subsection{Corollary}
\label{cordens}
{\it There exist a dense linear subspace $E$ of $\gN_{\chi'_2\circ j}$ such that $\Lambda_{\chi'_2\circ j}(E)$ is dense in $L^2(\pi(M'_0\cap M_2)', \chi'_2\circ j)$, and :}
\[J_{\chi'_2\circ j}\Lambda_{\chi'_2\circ j}(E)\subset D({}_{\hat{r}}L^2(\pi(M'_0\cap M_2)', \chi'_2\circ j))\cap D(L^2(\pi(M'_0\cap M_2)', \chi'_2\circ j)_s)\]

\begin{proof}
Let $E$ be the linear subspace generated by the elements of the form $e_nx$, where $e_n$ are the analytic elements associated to the Radon-Nykodym derivative $\delta'$, defined in [V3], and $x$ belongs to the subset $T_{\chi'_2, \tilde{T'_2}, \chi'_2\circ j, j\circ\tilde{T'_2}\circ j}$ defined in \ref{thdens2}; it is clear that $E$ is a subset of $\gN_{\chi'_2\circ j}$, dense in $\pi(M'_0\cap M_2)'$, and that $\Lambda_{\chi'_2\circ j}(E)$ is dense in $L^2(\pi(M'_0\cap M_2)', \chi'_2\circ j)$. 
\newline
As $E\subset\gN_{\chi'_2\circ j}\cap\gN_{j\circ\tilde{T'_2}\circ j}$, we get, using \ref{prop}(i), that 
\[J_{\chi'_2\circ j}\Lambda_{\chi'_2\circ j}(E)\subset D({}_{\hat{r}}L^2(\pi(M'_0\cap M_2)', \chi'_2\circ j))\]
\newline
Using ([V3]), we get that, if $x$ is in $T_{\chi'_2, \tilde{T'_2}, \chi'_2\circ j, j\circ\tilde{T'_2}\circ j}$, we can identify the vector $J_{\chi'_2\circ j}\Lambda_{\chi'_2\circ j}(e_nx)$ with $\delta'^{-1/2}J_{\chi_2}\Lambda_{\chi'_2}(e_nx)$. As we have :
\[e_nx\delta'^{-1/2}=(e_n\delta'^{-1/2})(\delta'^{1/2}x\delta'^{-1/2})\]
and, by \ref{thdens2}, $\delta'^{1/2}x\delta'^{-1/2}$ is a bounded operator in $\gN_{\tilde{T'_2}}$, so is $e_nx\delta'^{-1/2}$, and, therefore, we have :
\[\delta'^{-1/2}J_{\chi_2}\Lambda_{\chi'_2}(e_nx)=\pi(\lambda)^{1/4}J_{\chi_2}\Lambda_{\chi'_2}(e_nx\delta'^{-1/2})\subset J_{\chi_2}\Lambda_{\chi'_2}(\gN_{\chi'_2}\cap\gN_{\tilde{T'_2}})\]
and we get, by \ref{proptau}(i), that $J_{\chi'_2\circ j}\Lambda_{\chi'_2\circ j}(e_nx)$ belongs to $D(L^2(\pi(M'_0\cap M_2)', \chi'_2)_s)$, and, by linearity, we get the result. \end{proof}

\section{Hopf-bimodule structures on the relative commutants}
\label{relcom}
In this section, we put all the structures obtained on the relative commutant algebra $M'_0\cap M_2$.  We obtain this way on this algebra a canonical structure of Hopf-bimodule, with a co-inverse, a left-invariant operator-valued weight, and a right-invariant operator-valued weight (\ref{thleft}). We obtain also, thanks to chapters \ref{modular2} and \ref{density}, modular and analytical properties of this object. 

\subsection{Lemma}
\label{jj}
{\it Let $\mathcal J$ be the canonical anti-isomorphism from $\pi(M'_0\cap M_2)$ to $\pi(M'_0\cap M_2)'$ given, for $x$ in  by $M'_0\cap M_2$, by $\mathcal J(\pi (x))=J_{\chi_2}\pi(x)^*J_{\chi_2}$. Then, we get that $j\circ\mathcal J(\pi(x))=\mathcal J(\pi(j_1(x))$. }

\begin{proof}
It has been proved in ([E2], Prop. 7.7) \end{proof}

\subsection{Theorem}
\label{thleft}
{\it Let $M_0\subset M_1$ be a depth 2 inclusion of $\sigma$-finite von Neumann algebras, equipped with a regular normal semi-finite faithful operator-valued weight $T_1$. Let us suppose there exists on $M'_0\cap M_1$ a normal faithful semi-finite weight $\chi$ invariant under the modular automorphism group of $T_1$. Then, }
\newline {\it (i) there exists an application $\tilde{\Gamma}$ from $M'_0\cap M_2$ to 
\[(M'_0\cap M_2)_{j_1}\underset{M'_0\cap M_1}{*}{}_{id}(M'_0\cap M_2)\]
such that $(M'_0\cap M_1, M'_0\cap M_2, id, j_1, \tilde{\Gamma})$ is a Hopf-bimodule, (where $id$ means here the injection of $M'_0\cap M_1$ into $M'_0\cap M_2$, and $j_1$ means here the restriction of $j_1$ to $M'_0\cap M_1$, considered then as an anti-representation of $M'_0\cap M_1$ into $M'_0\cap M_2$). Moreover, the anti-automorphism $j_1$ of $M'_0\cap M_2$ is a co-inverse for this Hopf-bimodule structure. }
\newline {\it (ii) $\tilde{T_2}$ is a left-invariant operator-valued weight from $M'_0\cap M_2$ to $M'_0\cap M_1$, and, therefore, $j_1\circ\tilde{T_2}\circ j_1$ is a right-invariant operator-valued weight from $M'_0\cap M_2$ to $M'_1\cap M_2$. } 
\newline {\it (iii) Let $\chi_2$ be the weight $\chi\circ \tilde{T_2}$; there exist a one-parameter automorphism group $\tilde{\tau_t}$ of $M'_0\cap M_2$, commuting with the modular automorphism group  $\sigma_s^{\chi_2}$, such that, for all $t$ in $\mathbb{R}$, we have :
\[\tilde{\Gamma}\circ\sigma_t^{\chi_2}=(\tilde{\tau_t}_{j_1}\underset{\chi}{*}{}_{id}\sigma_t^{\chi_2})\circ\tilde{\Gamma}\]
Moreover, we have $j_1\circ\tilde{\tau_t}=\tilde{\tau_t}\circ j_1$. }

\begin{proof}
For simplification, let us identify $M'_0\cap M_2$ with $\pi(M'_0\cap M_2)$; then $j\circ\mathcal J=\mathcal J\circ j_1$ is an isomorphism from $M'_0\cap M_2$ onto $\pi(M'_0\cap M_2)'$. We just then use this isomorphism to get all the results, from those obtained in chapter \ref{bW}.  \end{proof}

\subsection{Theorem}
\label{thleft2}
{\it Let $M_0\subset M_1$ be a depth 2 inclusion of $\sigma$-finite von Neumann algebras, equipped with a regular normal semi-finite faithful operator-valued weight $T_1$. Let us suppose there exists on $M'_0\cap M_1$ a normal faithful semi-finite weight $\chi$ invariant under the modular automorphism group of $T_1$. Then the Hopf-bimodule structure $(M'_0\cap M_1, M'_0\cap M_2, id, j_1, \tilde{\Gamma})$ with co-inverse $j_1$, left-invariant operator valued weight $\tilde{T_2}$ (and right-invariant $j_1\circ\tilde{T_2}\circ j_1$) bears the following properties :
\newline (i) The modular automorphism groups $\sigma_t^{\chi_2}$ and $\sigma_s^{\chi_2\circ j}$ commute, and there exist a positive invertible operator $\lambda$ affiliated to $Z(M_0)\cap Z(M_1)$, and a positive invertible operator $\delta$ affiliated to $M'_0\cap M_2$ such that, for all $s$, $t$ in $\mathbb{R}$, and  :
\[\sigma_t^{\chi_2}(\delta^{is})=\lambda^{ist}(\delta^{is})\]
\[(D\chi_2 : D\chi_2\circ j_1)_t=\lambda^{it^2/2}\delta^{it}\]
\[\chi_2\circ\sigma_t^{\chi_2\circ j_1}(x)=\chi_2\circ\tilde{\tau_t}(x)=\chi_2(\lambda^{t/2}x\lambda^{t/2})\]
Moreover, it is possible to define on elementary tensors a one-parameter group of unitaries $\delta^{it}{}_{j_1}\underset{\chi}{\otimes}{}_{id}\delta^{it}$, whose generator satisfies :}
\[\Gamma(\delta)=\delta{}_{j_1}\underset{\chi}{\otimes}{}_{id}\delta\]
\newline {\it (ii) There exists a dense subspace $E$ of $\gN_{\chi_2}$ such that the subspace $\Lambda_{\chi_2}(E)$ is dense in $\gH$, and included in $D({}_{id}\gH)\cap D(\gH_{j_1})$. }

\begin{proof}
We use the same isomorphism as in \ref{thleft} to obtain these results from chapters \ref{modular2} and \ref{density}.  \end{proof}

\section{Lesieur's measured quantum groupoids}
\label{adapted}
In this chapter, we restrict to the case where the weight $\chi$ defined on the basis has the same modular automorphism group than $T_1$, i.e. when we have, for all $t$ in $\mathbb{R}$, $\sigma_t^{\chi}=\sigma_t^{T_1}$. 

\subsection{Definition}
\label{defadapted}
Let $\chi$ be a normal semi-finite faithful weight on $M'_0\cap M_1$; following [E2], we shall say that $\chi$ is adapted to $T_1$ if, for all $t$ in $\mathbb{R}$, we have $\sigma_t^{\chi}=\sigma_t^{T_1}$; it is then clear that $\chi$ is invariant under $\sigma_t^{T_1}$, and we may apply all the preceeding results. It must be noticed that the weight $\chi^2=\chi\circ j_1$ on $M'_1\cap M_2$ is not adapted to $T_2$ : the calculations made in \ref{inv} show that we have, for all $t$ in $\mathbb{R}$ and $t$ in $M'_1\cap M_2$, $\sigma_t^{\chi^2}(y)=\sigma_{-t}^{T_2}(y)$. We shall say then that $\chi^2$ is anti-adapted to $T_2$. 

\subsection{Theorem}
\label{thadapted}
{\it Let $M_0\subset M_1$ be a depth 2 inclusion of $\sigma$-finite von Neumann algebras, equipped with a regular normal semi-finite faithful operator-valued weight $T_1$. Let us suppose there exists on $M'_0\cap M_1$ an adapted normal faithful semi-finite weight $\chi$, as defined in \ref{defadapted}; then the Hopf-bimodule constructed in \ref{thleft} by defining $\tilde{\Gamma}'$ from $M'_1\cap M_3$ to 
\[(M'_1\cap M_3)_{j_2}\underset{M'_1\cap M_2}{*}{}_{id}(M'_1\cap M_3)\]
and taking the left invariant operator-valued weight $\tilde{T_3}=T_{3|M'_1\cap M_3}$ and the right-invariant operator-valued weight $j_2\circ\tilde{T_3}\circ j_2$, is a "measured quantum groupoid" in the sense of Lesieur ([L]).}

\begin{proof}
Let us apply \ref{corj} to the inclusion $M_1\subset M_2$; we obtain, that, for all $x$ in $M'_2\cap M_3$, and $t$ in $\mathbb{R}$, we have $\sigma_t^{\tilde{T_3}}(x)=\sigma_t^{T_3}(x)$; moreover, using both times the remark done in \ref{defadapted}, we see that the weight $\chi^3=\chi^2\circ j_2$ is adapted to $T_3$, and therefore, we have $\sigma_t^{T_3}(x)=\sigma_t^{\chi^3}(x)$. 
\newline
So, for any $y$ in $M'_1\cap M_2$, we have :
\[\sigma_t^{\tilde{T_3}}(j_2(y))=\sigma_t^{T_3}((j_2(y))=\sigma_t^{\chi^3}((j_2(y))=j_2(\sigma_{-t}^{\chi^2}(y))\]
which means that the weight $\chi^2$ is adapted to $\tilde{T_3}$, in the sense of Lesieur. 
\end{proof}

\section{Bibliography}
[AR] C. Anantharaman-Delaroche et J. Renault : Amenable Groupo-ids; Monographies de l'Enseignement Math\'ematique, {\bf 36}. L'Enseigne-ment
Math\'ematique, Gen\`eve, 2000. 196 pp
\newline\indent
[BS] S. Baaj et G. Skandalis : Unitaires multiplicatifs et dualit\'{e} pour les produits  crois\'{e}s
de
$\mathbb{C}^*$-alg\`{e}bres, {\it Ann. Sci. ENS}, {\bf 26} (1993), 425-488.
\newline\indent
[BSz1] G. B\"{o}hm and K. Szlach\'{a}nyi : A Coassociative $\mathbb{C}^*$-Quantum group with Non
Integral Dimensions, {\it Lett. Math. Phys.}, {\bf 38} (1996), 437-456.
\newline\indent
[BSz2] G. B\"{o}hm and K. Szlach\'{a}nyi : Weak $\mathbb{C}^*$-Hopf Algebras : the coassociative
symmetry of non-integral dimensions, in Quantum Groups and Quantum spaces {\it Banach Center
Publications}, {\bf 40} (1997), 9-19.
\newline\indent
[C1] A. Connes: On the spatial theory of von Neumann algebras, 
{\it J. Funct. Analysis}, {\bf 35} (1980), 153-164.
\newline\indent
[C2] A. Connes: Non commutative Geometry, Academic Press, 1994
\newline\indent
[E1] M. Enock : Inclusions irr\'eductibles de facteurs et
unitaires multiplicatifs II, {\it J. Funct. Analysis}, {\bf 154} (1998), 67-109.
\newline\indent
[E2] M. Enock : Inclusions of von Neumann algebras and quantum groupo\"{\i}ds II, {\it J. Funct. Analalysis}, {\bf 178} (2000), 156-225.
\newline\indent
[EN] M. Enock, R. Nest : Inclusions of factors,
multiplicative unitaries and Kac algebras, {\it J. Funct. Analysis}, {\bf137} (1996),
466-543.
\newline\indent
[ES] M. Enock, J.-M. Schwartz : Kac algebras and Duality of locally compact Groups,
Springer-Verlag, Berlin, 1989.
\newline\indent
[EV] M. Enock, J.-M. Vallin : Inclusions of von Neumann algebras and quantum groupo\"{\i}ds,
{\it J. Funct. Analalysis}, {\bf 172} (2000), 249-300.
\newline\indent
[GHJ] F.M. Goodman, P. de la Harpe, V.F.R. Jones: Coxeter
Graphs and Towers of Algebras, MSRI Publ. 14 , 1989, Springer.
\newline\indent
[J] V. Jones: Index for subfactors, {\it Invent. Math.},  {\bf 72} (1983), 1-25.
\newline\indent
[KV1] J. Kustermans and S. Vaes, Locally compact quantum groups, {\it Ann. Sci. ENS}, {\bf 33} (2000), 837-934.
\newline\indent
[KV2] J. Kustermans and S. vaes, Locally compact quantum groups in the von Neumann algebraic setting, {\it Math. Scand.}, {\bf 92} (2003), 68-92. 
\newline\indent
[L] F. Lesieur : thesis, University of Orleans, available at :
\newline
http://tel.ccsd.cnrs.fr/documents/archives0/00/00/55/05
\newline\indent
[MN] T. Masuda and Y. Nakagami : A von Neumann Algebra framework for
the duality of the quantum groups; {\it Publ. RIMS Kyoto}, {\bf 30} (1994), 799-850.
\newline\indent
[NV1] D. Nikshych, L. Va\u{\i}nerman : Algebraic versions of a finite dimensional quantum
groupoid, in Lecture Notes in Pure and Applied Mathematics, Marcel Dekker, 2000.
\newline\indent
[NV2] D. Nikshych, L. Va\u{\i}nerman : A characterization of depth 2 subfactors of $II_1$
factors, {\it J. Funct. Analysis}, {\bf 171} (2000), 278-307.
\newline\indent
[PP] M. Pimsner and S. Popa : Iterating the basic
construction, {\it Transactions AMS}, 
 {\bf 310}-1 (1988), 127-133.
\newline\indent
[R1] J. Renault : A Groupoid Approach to $\mathbb{C}^*$-Algebras, {\it Lecture Notes in
Math.} {\bf 793}, Springer-Verlag
\newline\indent
[R2] J. Renault : The Fourier algebra of a measured groupoid and its multipliers, {\it J. Funct. Analysis}, {\bf 145} (1997), 455-490.
\newline\indent
[S1] J.-L. Sauvageot : Produit tensoriel de $Z$-modules et applications, in Operator
Algebras and their Connections with Topology and Ergodic Theory, Proceedings Bu\c{s}teni, Romania,
1983, {\it Lecture Notes in
Math.} {\bf 1132}, Springer-Verlag, 468-485.
\newline\indent
[S2] J.-L. Sauvageot : Sur le produit tensoriel
relatif d'espaces de Hilbert,  {\it J. Operator Theory}, {\bf 9} (1983), 237-352.
\newline\indent
[St] \c{S}. Str\u{a}til\u{a} : Modular theory in
operator algebras, Abacus Press,  Turnbridge Wells, England, 1981.
\newline\indent
[Sz] K. Szlach\'{a}nyi : Weak Hopf algebras, in Operators Algebras and Quantum Field Theory,
S. Doplicher, R. Longo, J.E. Roberts, L. Zsido editors, International Press, 1996.
\newline\indent
[T] M. Takesaki : Theory of Operator Algebras II, Springer, Berlin, 2003. 
\newline\indent
[V1] S. Vaes : Locally compact quantum groups, thesis, Katholieke Universitet Leuven, 2001. 
\newline\indent
[V2] S. Vaes : The unitary implementation of a locally compact quantum group action, {\it J. Funct. Analysis}, {\bf 180} (2001), 426-480.
\newline\indent
[V3] S. Vaes : A Radon-Nikodym theorem for von Neumann algebras, {\it J. Operator Theory}, {\bf 180} (2001), 477-489. 
\newline\indent
[Val1] J.-M. Vallin : Bimodules de Hopf
et Poids op\'eratoriels de Haar, {\it J. Operator theory}, {\bf 35} (1996), 39-65
\newline\indent
[Val2] J.-M. Vallin : Unitaire pseudo-multiplicatif associ\'e \`a un groupo\"{\i}de;
applications \`a la moyennabilit\'e,  {\it J. Operator theory}, {\bf 44} (2000), 347-368.
\newline\indent
[Val3] J.-M. Vallin : Groupo\"{\i}des quantiques finis, {\it J. Algebra}, {\bf 239} (2001), 215-261.
\newline\indent
[Val4] J.-M. Vallin : Multiplicative partial isometries and finite quantum groupoids, in Locally Compact Quantum Groups and Groupoids, IRMA Lectures in Mathematics and Theoretical Physics 2, V. Turaev, L. Vainerman editors, de Gruyter, 2002. 
\newline\indent
[W1] S.L. Woronowicz : Tannaka-Krein duality
for compact matrix pseudogroups. Twisted
$SU(N)$ group. {\it Invent. Math.}, {\bf 93} (1988), 35-76.
\newline\indent
[W2] S.L. Woronowicz : Compact quantum group, in "Sym\'{e}tries quantiques" (Les Houches, 1995), North-Holland, Amsterdam (1998),
845-884.
\newline\indent
[W3] S.L. Woronowicz : From multiplicative unitaries to quantum groups, {\it Int. J. Math.}, {\bf 7}
(1996), 127-149. 
\newline\indent
[Y1] T. Yamanouchi : Duality for actions and coactions of
measured Groupoids on von Neumann Algebras, {\it Memoirs of the A.M.S.}, {\bf 101} (1993), 1-109.
\newline\indent
[Y2] T. Yamanouchi : Duality for generalized Kac algebras and a characterization on
finite groupoids algebras, {\it J. Algebra}, {\bf 163} (1994), 9-50.
\newline\indent

\end{document}